\documentclass[preprint,3p,12pt,sort&compress]{elsarticle}

\usepackage{natbib}

\usepackage{booktabs}
\usepackage{color}
\usepackage{amsmath}
\usepackage{amsfonts}
\usepackage{amssymb}
\usepackage{latexsym}
\usepackage{array}
\usepackage{amsthm}
\usepackage{graphics}
\usepackage{comment}
\usepackage{epsf,graphicx,subfigure}
\usepackage{epsfig}
\usepackage{epstopdf}
\usepackage{stmaryrd}
\usepackage{tikz}
\usepackage{hyperref}

\usepackage[normalem]{ulem}

\usetikzlibrary{mindmap}

\usetikzlibrary{arrows,positioning,shapes.geometric}

\usepackage{longtable}
\usepackage{rotating}
\usepackage{lscape}
\usepackage{adjustbox}

\usepackage{verbatim}
\usetikzlibrary{positioning,shadows,arrows}

\newtheorem{Theorem}{Theorem}
\newtheorem{Proposition}{Proposition}
\newtheorem{Lemma}{Lemma}
\newtheorem{Corollary}{Corollary}
\newtheorem{Assumption}{Assumption}
\newtheorem{Remark}{Remark}
\newtheorem{Definition}{Definition}
\newtheorem{Example}{Example}

\makeatletter
\def\ps@pprintTitle{%
  \let\@oddhead\@empty
  \let\@evenhead\@empty
  \def\@oddfoot{\reset@font\hfil\thepage\hfil}
  \let\@evenfoot\@oddfoot
}
\makeatother

\begin{document}

\begin{frontmatter}

\title{Probabilistic Analysis of Scalogram Ridges in Signal Processing}
\author[GRLiu]{Gi-Ren Liu}
\author[YCSheu]{Yuan-Chung Sheu}
\author[HTWu]{Hau-Tieng Wu\corref{mycorrespondingauthor}}
\cortext[mycorrespondingauthor]{Hau-Tieng Wu}
\ead{hauwu@cims.nyu.edu}

\address[GRLiu]{Department of Mathematics, National Cheng Kung University, Tainan, Taiwan}
\address[YCSheu]{Department of Applied Mathematics, National Yang Ming Chiao Tung University, Hsinchu, Taiwan}
\address[HTWu]{Courant Institute of Mathematical Sciences, New York University, New York, United States of America}

\begin{abstract}
While ridges in the scalogram, determined by the squared modulus of analytic wavelet transform (AWT), is a widely accepted concept and utilized in nonstationary time series analysis, their behavior in noisy environments remains underexplored.
Our objective is to provide a theoretical foundation for scalogram ridges by defining ridges as a potentially set-valued random process connecting local maxima of the scalogram along the scale axis and analyzing their properties when the signal fulfills the adaptive harmonic model and is contaminated by stationary Gaussian noise.
In addition to establishing several key properties of the AWT for random processes, we investigate the probabilistic characteristics of the resulting random ridge points in the scalogram. Specifically, we establish the uniqueness property of the ridge point at individual time instances and prove the upper hemicontinuity of the ridge random process. Furthermore, we derive bounds on the probability that the deviation between the ridges of noisy and clean signals exceeds a specified threshold, and these bounds depend on the signal-to-noise ratio. To achieve these ridge deviation results, we derive maximal inequalities for the complex modulus of nonstationary Gaussian processes, leveraging classical tools such as the Borell-TIS inequality and Dudley's theorem, which might be of independent interest.
\end{abstract}

\begin{keyword}
adaptive harmonic model; analytic wavelet transform;
Borell-TIS inequality;
circularly symmetric Gaussian process;
maximum inequality; ridge; ridge deviation.
\MSC[2010]: 42C40, 60G15, 60G60, 60G70, 60H30.

\end{keyword}

\end{frontmatter}

\section{Introduction}

\begin{figure}[hbt!]
\centering
\includegraphics[width=0.99\textwidth]{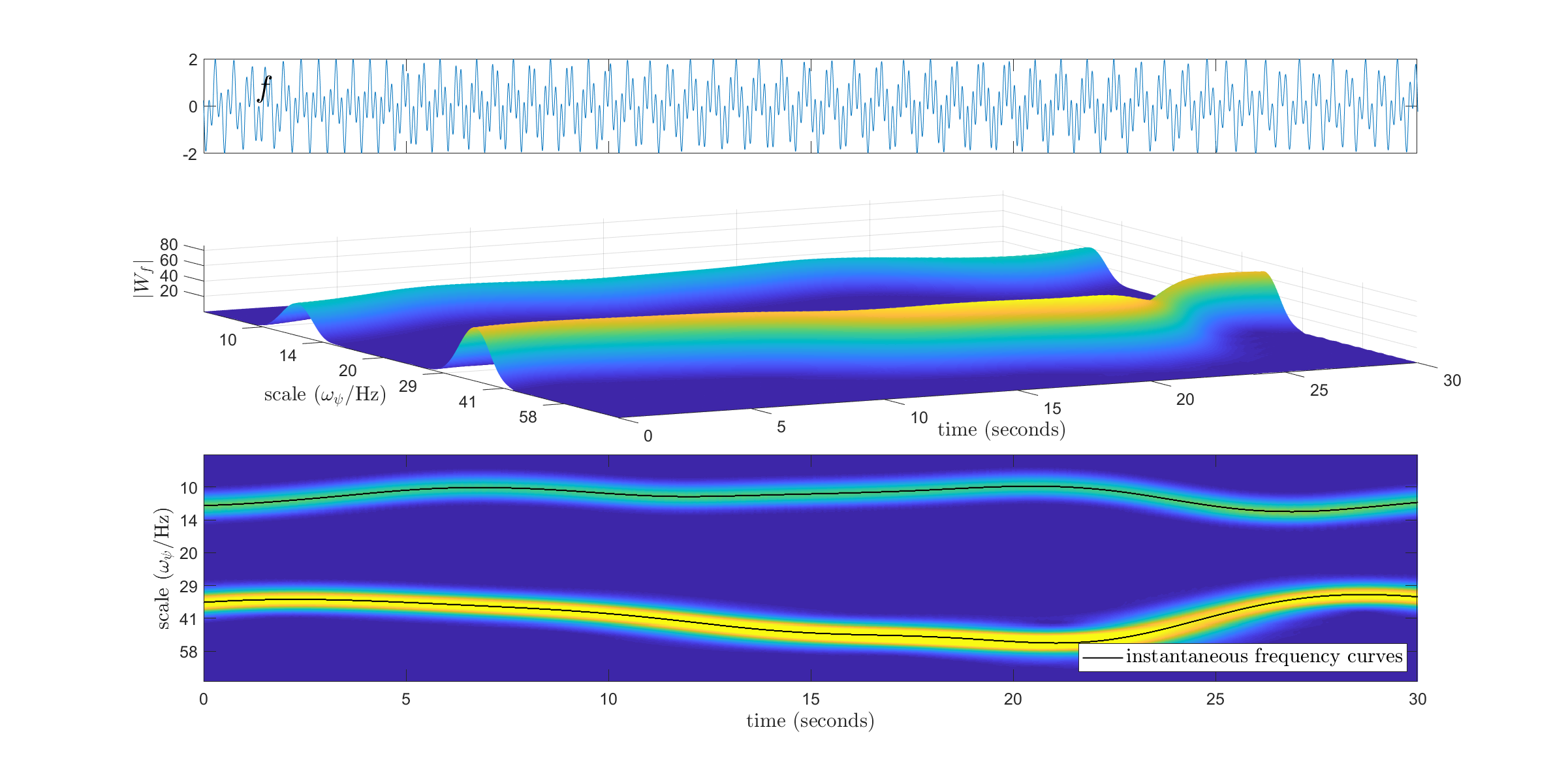}
\hspace{-1.5cm}
\includegraphics[width=0.99\textwidth]{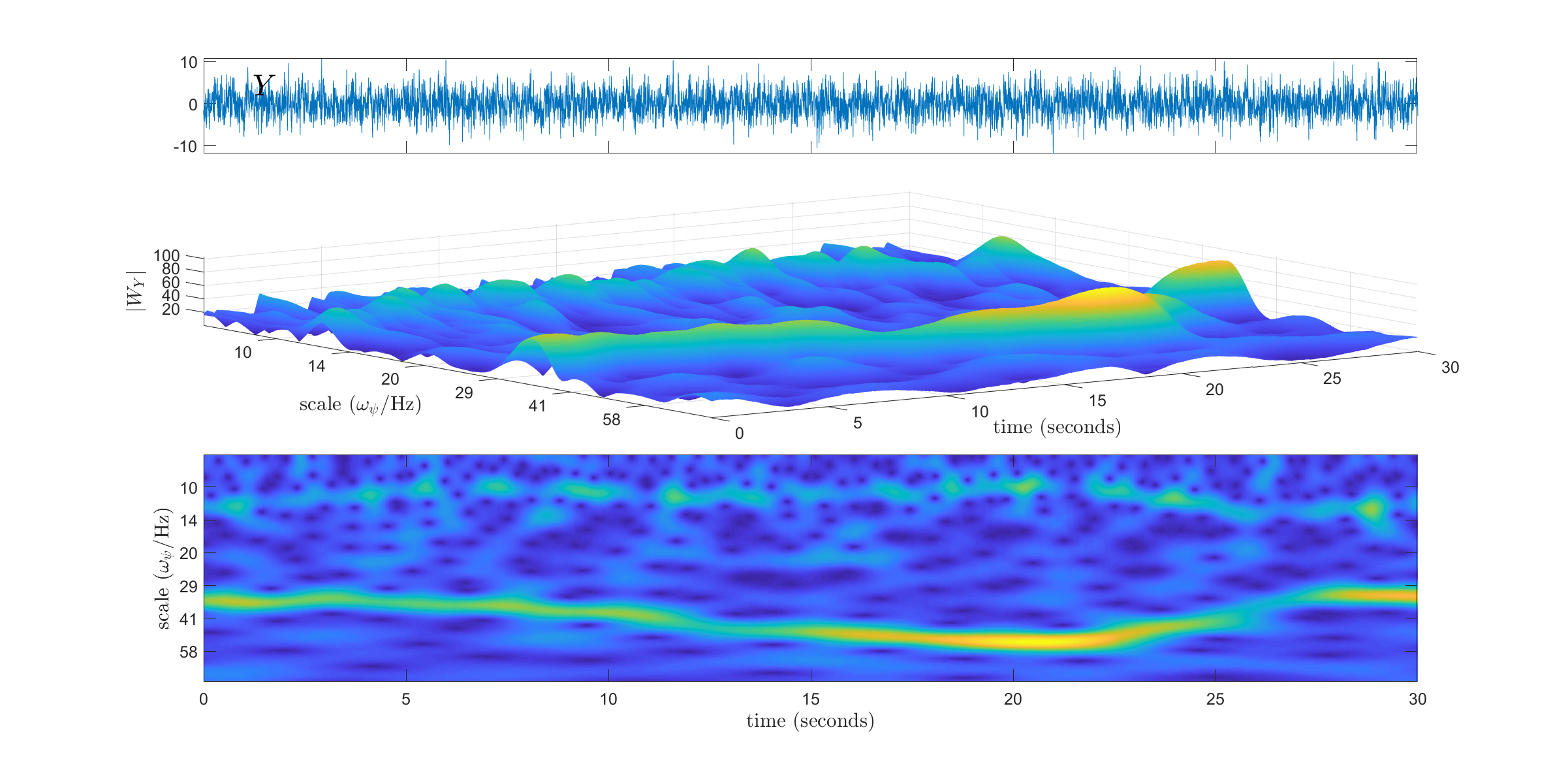}
\caption{The top panel illustrates a signal $f$, composed of two frequency-modulated sinusoids, along with its 3-dimensional and image representations in the time-scale domain.
Here, the mother wavelet used has center frequency $\omega_{\psi}= 80$ Hz, as defined in (\ref{def:centralfrequency_psi}).
Hence, the signal $f$'s IF lies approximately within the interval
$[\frac{\omega_{\psi}}{58},\frac{\omega_{\psi}}{10}]$  Hz.
The bottom panel presents the same for its noise-affected counterpart $Y$, with a signal-to-noise ratio of -8.98 dB. Here, $Y$ is the sum of $f$ and a sample path of the Gaussian process $\Phi$.
The y-axis represents the scale variable
$s$, which is dimensionless.
The tick values $10, 14, \ldots, 58$ are scale values rather than physical frequencies.
For a mother wavelet with central frequency
$\omega_{\psi}=80$ Hz, a scale value $s$ corresponds to the frequency $\omega_{\psi}/s$ Hz.
Thus, the tick values
$10, 14, \ldots,  58$ correspond to the frequencies
$80/10, 80/14,\ldots,80/58$ Hz, respectively.
 For the single-component case, see Figure \ref{fig:wavelet_potential_ridge_Surface plot:one} in the appendix.
 }\label{fig:wavelet_potential_ridge_Surface plot}
\end{figure}

Ridge analysis plays a vital role in connecting time-frequency analysis \cite{daubechies1992ten,flandrin1998time} and signal processing \cite{mallat1999wavelet}, particularly for capturing dynamic information.
Traditionally, ridges in time-frequency representations (TFRs) are interpreted ad hoc as sequences of local maxima along the frequency or scale axis, often visualized as continuous curves over mountain-shaped structures in the TFR. See Figure \ref{fig:wavelet_potential_ridge_Surface plot}. This continuity is central to definitions for adaptive harmonic signals \cite{lilly2010analytic} and underlies many ridge extraction algorithms \cite{carmona1999multiridge,carmona1997characterization,rankine2007if,Colominas_Meignen_Pham_2020,iatsenko2016extraction,Laurent_Meignen_2021,Legros_Fourer_2021,meignen2017demodulation,Ozkurt_Savaci_2005,su2024ridge,Zhu_Zhang_Gao_Li_2019}.
Ridges encode instantaneous frequency (IF), making them foundational for various applications \cite{flandrin1998time,mallat1999wavelet}.
Despite widespread use and the development of many algorithms, the statistical properties of ridges and algorithmic accuracy under noise remain underexplored.
While the term ``ridge'' also appears in other fields such as statistics and image analysis
\cite{hall1992ridge, eberly2012ridges}, its meaning in TFRs is distinct. In what follows, we focus exclusively on ridges in TFRs.

Before proceeding, we outline a general framework for ridge analysis that guides this study. Understanding a signal scientifically involves four steps. In Step \texttt{I}, identify the quantity of interest, which is ``ridge'' in this paper, and model the signal. In Step \texttt{II}, analyze its properties and define it precisely. In Step \texttt{III}, develop extraction algorithms.
In Step \texttt{IV}, assess algorithm performance and statistical behavior.
To our knowledge, most existing ridge extraction methods skip Steps \texttt{I} and \texttt{II}, relying on heuristics with limited theoretical grounding.
Within this framework, a {\em signal} refers either to the observed time series, modeled as a random process $Y(t)=f(t)+\Phi(t)$, where $f$ is a deterministic oscillatory component and $\Phi$ is a mean-zero stochastic process representing noise, or to its TFR. The ridge is a structural property of the TFR. Therefore, the signal we analyze here is the TFR of $Y$, which inherits the statistical and structural characteristics of the process $Y$.
This paper focuses on Steps \texttt{I} and \texttt{II}, establishing the ridge's statistical and structural characteristics as a foundation for practical applications. While we defer algorithm development to future work, we include numerical results and extensively discuss their relationship to existing algorithms and to Steps \texttt{III} and \texttt{IV} in the Discussion section.

Following this framework, in Step \texttt{I}, we model the noisy signal as $Y=f+\Phi$, where $f$ fulfills the adaptive harmonic model (AHM) \cite{2011synchrosqueezed} with possibly multiple oscillatory components having well-separated IFs, and $\Phi$ is modeled as a stationary Gaussian random process with potential long- or short-range dependence. We analyze the TFR of $Y$ obtained using the continuous wavelet transform with an analytic mother wavelet, referred to as the analytic wavelet transform (AWT), which is a complex-valued function $W_{Y}:(t,s)\to \mathbb{C}$, where $t\in \mathbb{R}$ is time and $s>0$ is scale, whose squared modulus  $S_{Y}:(t,s)\to \mathbb{R}_+\cup\{0\}$,  the {\em scalogram}, describes the spectral energy distribution over time. We call the set $\{t\in \mathbb{R},s>0\}$ the {\em time-scale domain}. It is important to note that when the analysis centers on IF, researchers typically employ the AHM to facilitate the estimation of IFs by extracting energy-concentrated curves from the scalograms (or spectrograms) of noisy signals \cite{carmona1997characterization,carmona1999multiridge,iatsenko2016extraction,stiles2004wavelet,su2024ridge}.

In Step \texttt{II}, we examine the behavior of $W_Y$. When the IFs of $f$ are well-separated, the scalogram of $f$,
denoted by $S_f$, encodes IF information through its local maxima along the scale axis \cite{delprat1992asymptotic}.
The AWT of the Gaussian noise process $\Phi$, denoted by $W_\Phi$, interacts with the AWT of $f$, denoted by $W_f$,
via the relation $S_{Y}=|W_{f}+W_{\Phi}|^2$. How the random field $W_{\Phi}$ perturbs the local maxima in $S_Y$ is the main focus of this paper. By viewing $S_Y$ as an inhomogeneous random field, the study of its local maxima is naturally related to excursion set analysis \cite{adler2010geometry}, though the current theory does not fully account for the nonstationarity or the temporal coherence critical to ridge structure.
Additional questions, such as how the local maxima of $S_f$, and thus of $S_Y$, are affected by spectral inference when the IFs of $f$ are close \cite{delprat2002global} or what local maxima mean for the general $f\in L^2$ are important but deferred to future work.

As an illustrative example, consider Figure~\ref{fig:wavelet_potential_ridge_Surface plot}.
The top panel shows that the TFR of a clean signal $f$, modeled by the AHM (\ref{AMFM_signal}),
produces continuous and well-defined ridges.
In contrast, the bottom panel displays the TFR of the noisy signal $Y$.
The TFR of $Y$ exhibits multiple local maxima per time slice and more irregular structures, clearly differing from the clean case.
Notably,  noise can cause abrupt jumps in the ridge, as seen around the 6th, 8th, and 26th seconds in Figures~\ref{fig:wavelet_ridge} and \ref{fig:wavelet_potential_ridge}. These observations underscore the need to rigorously study the statistical behavior of $W_Y$ and to define ridges precisely under noisy conditions.

Based on empirical observations and existing theoretical results, we define the ridge in Step \texttt{II} for a clean signal $f$ with a single oscillatory component as the graph of
\begin{align*}
 s_{f}(t)= \arg \underset{s>0}{\max}\ S_{f}(t,s),\ t\in \mathbb{R},
\end{align*}
commonly known as the wavelet ridge \cite{mallat1999wavelet} or wavelet skeleton
\cite[pp. 14--18]{antoine2008two}. This definition extends naturally to the observed random process by
\begin{align*}
s_{Y}(t)= \arg \underset{s>0}{\max}\ S_{Y}(t,s),\ t\in \mathbb{R},
\end{align*}
with $s_{Y}$ referred to as the ridge of $Y$. Although widely used, the properties of $s_Y$ remain largely unexplored. A central challenge lies in the fact that as a nonlinear functional of the random field $W_Y$, $s_Y$ is a random process that may be discontinuous or multi-valued.
When $f$ consists of multiple oscillatory components, each corresponds to a local maximum in  $S_f$ \cite{delprat1992asymptotic}, and under sufficient IF separation, individual ridges can be defined by segmenting the time-scale domain. We explain this in detail later. We shall mention that if the IFs are not well separated or even crossover, ridge definitions become ambiguous, even without noise \cite{delprat2002global,meignen2021study}, and this is beyond the scope of this paper.

Our contributions in this paper are multifaceted.
First, we analyze the AWT, focusing on its regularity and a novel covariance structure of the scalogram of $\Phi$. These aspects are critical but underexplored aspects of statistical inference.
Second, we propose a ridge definition for noisy signals suited to practical needs. Under mild conditions, we show that the set-valued random process $t\mapsto s_Y(t)$ is almost surely single-valued on any countable and dense subset of $\mathbb{R}$, and upper hemicontinuous on $\mathbb{R}$. This means that for any time $t$,  as $t'$ approaches $t$, the maxima of $S_{Y}(t',\cdot)$ converge arbitrarily close to some of the maxima of $S_{Y}(t,\cdot)$. Based on these properties, we further discuss the local $C^{1}$-smoothness of $s_{Y}(t)$.
Third, we analyze how $s_Y(t)$ deviates from $s_f(t)$. Since the scalogram involves a nonlinear modulus operator, analyzing $|s_{Y}(t)-s_{f}(t)|$ is difficult, particularly from the perspectives of time-frequency analysis and probability. This problem has received limited attention in the literature, with the exception of \cite{sejdic2008quantitative}, where the authors focused on a frequency-modulated signal contaminated with additive Gaussian white noise and provided approximations for $\mathbb{E}[s_{Y}(t)-s_{f}(t)]$ and $\textup{Var}(s_{Y}(t)-s_{f}(t))$.
Our work extends this analysis to temporally correlated Gaussian noise models.
In Theorems~\ref{lemma:conditional} and \ref{mainresult:deviation}, we provide upper bounds on the probability that $|s_{Y}(t)-s_{f}(t)|$ exceeds a given threshold, where the bounds decay exponentially as the noise strength decreases, offering theoretical support for the numerical observations in Figure \ref{fig:effect_SNR}.
These bounds are derived from new maximal inequalities for $|W_\Phi|$ and its derivative (Lemma \ref{lemma:complex-Borell}).
These extend classical maximal inequalities, such as the Borell-TIS inequality and Dudley's theorem, to complex-valued Gaussian processes arising from the AWT of $\Phi$ and may be of independent interest because the classical Gaussian concentration inequalities \cite{kuhn2023maximal}, Dudley's theorem, and the Borell-TIS inequality \cite{talagrand2014upper,vershynin2018high} are formulated for real-valued Gaussian elements.
Furthermore, Theorem \ref{mainresult:deviation} reveals how the spectral energy distribution of the noise influences the probability of ridge deviations, yielding insights that go beyond what can be obtained through the approximations of $\mathbb{E}[s_{Y}(t)-s_{f}(t)]$ and $\textup{Var}(s_{Y}(t)-s_{f}(t))$.

The rest of the paper is organized as follows.
In Section \ref{sec:preliminary}, we summarize the necessary material for AWT and present the definitions of the adaptive harmonic model and the noise we consider.
Section \ref{sec:AWT_old_new} is dedicated to AWT for Gaussian random processes, where we present both old and new results.
In Section \ref{sec:def:ridge}, we introduce a definition of the ridge that takes into account practical needs when the data is noisy.
In Section \ref{sec:mainresult}, we state our main results about the probability of the occurrence of ridge deviation. Section \ref{sec:conclusion} concludes the paper with a discussion of the results and future directions.
The proofs of our main results and some technical lemmas are provided in Section \ref{sec:proof}.
Finally, a list of frequently used symbols is summarized in Table~\ref{List_symbols} in the appendix.

\section{Preliminary}\label{sec:preliminary}

Consider the following {\em adaptive harmonic model (AHM)}. Take a real-valued signal $f$:
\begin{equation}\label{AMFM_signal}
f(t) = \overset{M}{\underset{m=1}{\sum}} A_{m}(t)\cos\left(2\pi \phi_{m}(t)\right),\ t\in \mathbb{R},
\end{equation}
where $A_{m}(t)\in C^{1}(\mathbb{R})$ are positive functions representing the amplitude modulation (AM) and $\phi_{m}(t)\in C^{2}(\mathbb{R})$ are strictly monotonically increasing functions representing the phase so that $\phi'_{m}(t)>0$ represents the IF.
We impose two conditions on the AHM model. The first one is the {\em slowly varying} condition: $|A'_{m}(t)|\leq \epsilon \phi'_{m}(t)$ and $|\phi''_{m}(t)|\leq \epsilon \phi'_{m}(t)$, where $\epsilon\geq 0$ is a small constant.
The second one is the {\em spectral separation} condition when $M>1$:
\begin{align}\label{assumption:spectral_separation}
\phi'_{1}(t)<\phi'_{2}(t)<\cdots<\phi'_{M}(t)
\end{align}
for all $t\in \mathbb{R}$ so that the {\em IF separation condition} $\phi'_{m+1}(t)-\phi'_{m}(t)\geq \delta(\phi'_{m+1}(t)+\phi'_{m}(t))$ for some $\delta\in (0,1)$ is satisfies.
When $M=1$,
we denote $A(t)=A_{1}(t)$ and $\phi(t) = \phi_{1}(t)$. We call $A_m(t)\cos(2\pi\phi_m(t))$ the $m$-th {\em intrinsic mode type (IMT)} function. While it is possible to consider more complicated models, such as those involving crossover IF \cite{chen2023disentangling}, IFs with insufficient separation \cite{meignen2021study}, or non-sinusoidal oscillation \cite{su2024ridge}, in this paper we focus on the AHM model defined above.

In the presence of noise, the noise-affected AHM may be expressed as
\begin{align}\label{noisycos}
Y(t) = f(t)+\Phi(t),\ t\in \mathbb{R},
\end{align}
where $\Phi$ is a random process satisfying the following assumption.
\begin{Assumption}\label{assumption:Gaussian}
The noise $\Phi$ is a mean-square continuous and stationary Gaussian process.
Additionally, we assume that the covariance function of $\Phi$ has a spectral density function
$p: \mathbb{R}\rightarrow (0,\infty)$
such that
\begin{equation*}
\textup{Cov}(\Phi(t),\Phi(0)) = \int_{\mathbb{R}}e^{it\lambda}p(\lambda)d\lambda.
\end{equation*}
\end{Assumption}
Under Assumption \ref{assumption:Gaussian},
the process $\Phi$ can be represented by
\begin{equation}\label{spectral_rep:samplepath}
\Phi(t) = \mathbb{E}[\Phi(0)] +
\int_{\mathbb{R}}e^{it\lambda}\sqrt{p(\lambda)}Z(d\lambda),
\end{equation}
where $Z$ is an orthogonally scattered Gaussian random measure on $\mathbb{R}$
satisfying
\begin{align}\label{ortho}
Z(\Delta_{1})=\overline{Z(-\Delta_{1})},\ \
\mathbb{E}[Z(\Delta_{1})]=0,\
\textup{and}\ \ \mathbb{E}\left[Z(\Delta_{1})
\overline{Z(\Delta_{2})}\right]=\textup{Leb}(\Delta_{1}\cap\Delta_{2})
\end{align}
for any $\Delta_{1},\Delta_{2}\in\mathcal{B}(\mathbb{R})$, where Leb is the Lebesgue measure on $\mathbb{R}$ and $\mathcal{B}(\mathbb{R})$ is the Borel {$\sigma$-}algebra on $\mathbb{R}$ \cite{major1981lecture}.

The analytic wavelet transform (AWT) of an $L^2$ function $f$, denoted by $\{W_{f}(t,s)\mid t\in \mathbb{R}, s>0\}$, is defined as follows \cite{daubechies1992ten}:
\begin{align}\label{def:WT}
W_{f}(t,s)
=\int_{\mathbb{R}} f(\tau) \frac{1}{s}\overline{\psi\left(\frac{\tau-t}{s}\right)} d\tau,
\end{align}
where $t\in \mathbb{R}$ represents time, $s>0$ represents scale,
and $\psi\in L^{1}(\mathbb{R})\cap L^{2}(\mathbb{R})$ is the chosen  analytic mother wavelet.
Recall that the mother wavelet $\psi$ is called analytic if $\widehat{\psi}(\omega)=0$ for $\omega\leq0$,
where $\widehat{\psi}$ is the Fourier transform of $\psi$ defined by
$\widehat{\psi}(\omega) = \int_{\mathbb{R}}e^{-i t\omega}\psi(t)dt,\ \omega\in \mathbb{R}.$
The squared magnitude of the complex-valued function $W_{f}$,
\begin{align}\label{def:scalogram}
S_{f}(t,s):=|W_{f}(t,s)|^{2},
\end{align}
is called the {\em scalogram}. See Figure \ref{fig:wavelet_potential_ridge_Surface plot} for an example.

Similar to (\ref{def:WT}) and \eqref{def:scalogram},
we denote the AWTs of $\Phi$ and $Y$ by $W_{\Phi}$ and $W_{Y}$, respectively, and the scalograms of $\Phi$ and $Y$ by $S_{\Phi}$ and $S_{Y}$, respectively. By the linearity of AWT, we know $W_Y=W_f+W_\Phi$.
By (\ref{spectral_rep:samplepath}) and the stochastic Fubini Theorem \cite{pipiras2010regularization},
\begin{align}\label{spect:WPhi}
W_{\Phi}(t,s) = \int_{0}^{\infty}e^{it\lambda} \overline{\widehat{\psi}(s\lambda)}\sqrt{p(\lambda)}Z(d\lambda).
\end{align}
The AWT of both $\Phi$ and $Y$ are complex-valued and nonstationary Gaussian random processes on the scale axis for each $t\in\mathbb{R}$.

\section{AWT for Gaussian random processes--some old and new results}\label{sec:AWT_old_new}

The first result focuses on the covariance structure of $W_\Phi$, which serves as a foundational step in analyzing the ridges of the TFR determined by AWT.
\begin{Proposition}\label{lemma:exp}
Suppose that $\Phi$ is a stationary Gaussian process satisfying Assumption \ref{assumption:Gaussian}.

\noindent$(a)$ For any $t\in \mathbb{R}$, $d\in \mathbb{N}$, and $s_{1},s_{2},\ldots,s_{d}>0$,
the complex Gaussian random vector
$$\mathbf{W}:=\left[W_{\Phi}(t,s_{1})\ W_{\Phi}(t,s_{2})\ \cdots\ W_{\Phi}(t,s_{d})\right]^{\top}$$
is circularly symmetric, where $\top$ represents the transpose. That is, $e^{i\theta}\mathbf{W}$ has the same probability distribution as $\mathbf{W}$ for any $\theta\in \mathbb{R}$.

\noindent$(b)$ For any $t\in \mathbb{R}$, $s_{1},s_{2}>0$ with $s_{1}\neq s_{2}$,
$$\frac{|W_{\Phi}(t,s_{1})-W_{\Phi}(t,s_{2})|^{2}}{d^{2}_{W_{\Phi}}(s_{1},s_{2})}$$
follows an exponential distribution with a mean of one,
where
\begin{align}\label{def:canonical_distance}
d_{W_{\Phi}}(s_{1},s_{2})=\left\{\mathbb{E}\left[\left|W_{\Phi}(0,s_{1})-W_{\Phi}(0,s_{2})\right|^{2}\right]\right\}^{1/2}.
\end{align}

\noindent$(c)$
For any $t\in \mathbb{R}$ and $s_{1},s_{2}>0$,
\begin{align*}
\textup{Cov}\left(S_{\Phi}(t,s_{1}),S_{\Phi}(t,s_{2})\right)
= 4\left\{\mathbb{E}\left[W^{R}_{\Phi}(t,s_{1})W^{R}_{\Phi}(t,s_{2})\right]\right\}^{2}
+4\left\{\mathbb{E}\left[W^{R}_{\Phi}(t,s_{1})W^{I}_{\Phi}(t,s_{2})\right]\right\}^{2},
\end{align*}
where $W^{R}_{\Phi}$ and $W^{I}_{\Phi}$ are the real and imaginary parts of $W_{\Phi}$, respectively.
Especially, for any $s>0$,
\begin{align*}
\textup{Var}\left(S_{\Phi}(t,s)\right) =  4\left\{\mathbb{E}\left[|W^{R}_{\Phi}(t,s)|^{2}\right]\right\}^{2}
=\left\{\mathbb{E}\left[S_{\Phi}(t,s)\right]\right\}^{2}.
\end{align*}
\end{Proposition}

When the noise is modeled as a stationary Gaussian process, the scalogram $S_{\Phi}$ at a given time is a nonstationary random process along the scale axis.
Proposition \ref{lemma:exp}(c) shows that for any $t\in \mathbb{R}$,
the nonstationary process $S_{\Phi}(t,\cdot)$ is non-negatively correlated.
Because this property is not found in existing literature, we provide its proof in Section \ref{sec:proof:lemma:exp}.

The goal of this work is to analyze the location of the maximum of $|W_{Y}(t,\cdot)|^{2}$ along the scale axis
and to examine its deviation from the location of the maximum of $|W_{f}(t,\cdot)|^{2}$ along the same axis.
For any fixed $t\in \mathbb{R}$, finding the maximizers of $|W_{Y}(t,\cdot)|^{2}$ requires
computing the partial derivatives of both
$W_{f}$ and $W_{\Phi}$ with respect to the scale variable since $W_Y=W_f+W_\Phi$.
Furthermore, it is essential to ensure the finiteness of the expectation of the maximum of $|W_{\Phi}|$ and $|\partial W_{\Phi}/\partial s|$
(see Theorem \ref{lemma:conditional} and Lemma \ref{lemma:complex-Borell} below),
so we make the following assumption.
\begin{Assumption}\label{assumption:boundedness:psi}
The Fourier transform $\widehat{\psi}$ of the analytic mother wavelet $\psi$ is
three times differentiable on $(0,\infty)$, and the following conditions hold:
\begin{align*}
\textup{($\mathrm{D}^{0}_{1}$-$\mathrm{D}^{0}_2$):}\  \underset{\lambda>0}{\sup}|\lambda^{p}\widehat{\psi}(\lambda)|<\infty,\ p=1,2;\ \ \
&\textup{($\mathrm{D}^{1}_{0}$-$\mathrm{D}^{1}_2$):}\  \underset{\lambda>0}{\sup}|\lambda^{p}D\widehat{\psi}(\lambda)|<\infty,\ p=0,1,2;
\\
\textup{($\mathrm{D}^{2}_2$-$\mathrm{D}^{2}_3$):}\  \underset{\lambda>0}{\sup}|\lambda^{p}D^{2}\widehat{\psi}(\lambda)|<\infty,\ p=2,3;\ \ \
&\textup{($\mathrm{D}^{3}_3$):}\  \underset{\lambda>0}{\sup}|\lambda^{3}D^{3}\widehat{\psi}(\lambda)|<\infty,
\end{align*}
where $D^{k}$ represents the $k$th derivative operator.
Without loss of generality, we also assume that $|\widehat{\psi}|$ is a unimodal function.
\end{Assumption}
Assumption \ref{assumption:boundedness:psi}
imposes mild smoothness and polynomial decay conditions on the Fourier transform of the analytic mother wavelet.
Since the polynomial decay of the Fourier transform suppresses high-frequency components in the noise, these conditions ensure that the scalogram of the analytic wavelet transform of the stationary Gaussian noise process $\Phi$ admits sample paths that are twice continuously differentiable on the time-scale domain; see Proposition \ref{lemma:twice_differentiability}.
The regularity is essential for establishing the local uniqueness and hemicontinuity properties of ridge curves (Theorems \ref{prop:unique_argmax} and \ref{lemma:compact&hemiconti}).
In particular, the uniqueness of the ridge point for the noise-affected signal guarantees that the deviation of the ridge point from its noise-free counterpart is well defined.
On the other hand, under the conditions $(\mathrm{D}^{0}_{1})$,
$(\mathrm{D}^{1}_{0})$, and $(\mathrm{D}^{1}_{1})$,
the continuity of the AWT of the Gaussian process $\Phi$ can be extended to the boundary of the time-scale domain by defining $W_{\Phi}(t,0)=0$. This boundary continuity is a key ingredient in our extension of Dudley's theorem to $\mathbb{C}$-valued nonstationary Gaussian random processes indexed by scale and generated by the AWT (see the proof of
Lemma \ref{lemma:Dudley}).

\begin{Example}\label{example:wavelet}
The set of wavelets that satisfy Assumption \ref{assumption:boundedness:psi}
includes some of the generalized Morse wavelets discussed in \cite[(14)]{lilly2010analytic}
and the analytic wavelets derived in \cite[(6)]{holighaus2019characterization}.
For instance, the generalized Morse wavelets, whose Fourier transform takes the form
$\widehat{\psi}(\lambda)= a_{\beta_{1},\beta_{2}}\lambda^{\beta_{1}}e^{-\lambda^{\beta_{2}}}1_{[0,\infty)}(\lambda),
\ \beta_{1}\geq1,\ \beta_{2}>0$,
where $a_{\beta_{1},\beta_{2}}$ is a normalizing constant.
Another example is the Klauder wavelet, whose Fourier transform is given by
$\widehat{\psi}(\lambda) = \lambda^{\alpha}e^{-\gamma \lambda}e^{i\beta \log \lambda} 1_{[0,\infty)}(\lambda)$,
where $\alpha\geq1$, $\beta\in \mathbb{R}$, and $\gamma\in \mathbb{C}$ with $\textup{Re}(\gamma)>0$.

\end{Example}

\begin{Proposition}\label{lemma:twice_differentiability}
Under Assumptions \ref{assumption:Gaussian} and \ref{assumption:boundedness:psi},
the sample paths of $W_{\Phi}(\cdot,\cdot)$  are two times continuously differentiable on $\mathbb{R}\times (0,\infty)$ almost surely.
\end{Proposition}

The proof of Proposition \ref{lemma:twice_differentiability} is provided
in Section \ref{sec:proof:lemma:twice_differentiability}. It shows that
if the mother wavelet $\psi$ has higher regularity, the sample paths of
$W_Y$ will also exhibit higher-order smoothness.
A surface plot of $|W_{Y}|$ is shown in the bottom portion of Figure \ref{fig:wavelet_potential_ridge_Surface plot}.

\begin{Assumption}\label{assumption:spectral_density}
For the behavior of the spectral density function $p$ of the noise $\Phi$ at infinity, we assume that $(a)$ there exist constants $\gamma\in(0,2)$ and $C_{1}>0$ such that
\begin{equation*}
p(\lambda) \leq  C_{1}|\lambda|^{-(1+\gamma)}\ \textup{for any}\ \lambda\neq0.
\end{equation*}
For the behavior of $p$ near the origin,
we consider two scenarios.
\begin{itemize}
\item $\mathrm{(b1)}$ Long-range dependence scenario:
\begin{equation}\label{spectral:slowly}
p(\lambda) =
L(|\lambda|^{-1})|\lambda|^{H-1},
\end{equation}
where $H$ is a long-memory parameter belonging to $(0,1)$ and $L: (0,\infty) \rightarrow
(0,\infty)$ is a locally bounded function which is slowly varying at infinity,  meaning $L(cr)/L(r)\rightarrow 1$ as $r\rightarrow\infty$ for any $c>0$.
\item $\mathrm{(b2)}$  Short-range dependence scenario: $p$ is a bounded function continuous at the origin.
\end{itemize}
\end{Assumption}

\begin{Example}\label{example:covariance}
Let us consider the stationary random process $\Phi$ with the generalized Linnik
covariance function \cite{lim2010analytic}, defined as
\begin{align*}
C_{\Phi}(t)=(1+|t|^{\gamma})^{-H/\gamma},\ \gamma\in(0,2],\ H>0.
\end{align*}
According to \cite[Propositions 3.1 and 3.4]{lim2010analytic}, $p(\lambda)=\mathcal{O}(|\lambda|^{-(1+\gamma)})$  when $|\lambda|\rightarrow\infty$.
For the case $H>1$, because $C_{\Phi}\in L^{1}(\mathbb{R})$, the spectral density function $p$ is continuous everywhere.
For the case $0<H<1$, by \cite[Corollaries 3.3 and 3.10]{lim2010analytic},
\begin{align*}
p(\lambda)\sim |\lambda|^{H-1}\Gamma(\frac{1-H}{2})\left[2^{H}\pi^{\frac{1}{2}}\Gamma\left(\frac{H}{2}\right)\right]^{-1}
\end{align*}
as $|\lambda|\rightarrow0$.
\end{Example}

\begin{Corollary}\label{lemma:continuity_extension}
Suppose that Assumption \ref{assumption:Gaussian} and conditions $(\mathrm{D}^{0}_{1})$,
$(\mathrm{D}^{1}_{0})$, and $(\mathrm{D}^{1}_{1})$ in Assumption \ref{assumption:boundedness:psi} hold.

\noindent(a) If Assumption \ref{assumption:spectral_density}(a) is satisfied, then
\begin{equation}\label{proof:continuity_0}
\mathbb{P}\left(\underset{s\rightarrow 0+}{\lim} W_{\Phi}(t,s)=0\ \textup{for all}\ t\in \mathbb{R} \right)=1.
\end{equation}
\noindent(b) If either Assumption \ref{assumption:spectral_density}(b1) or Assumption \ref{assumption:spectral_density}(b2) is satisfied,
then
\begin{equation}\label{proof:continuity_00a}
\mathbb{P}\left(\underset{s\rightarrow \infty}{\lim} W_{\Phi}(t,s)=0\ \textup{for all}\ t\in \mathbb{R} \right)=1.
\end{equation}
\end{Corollary}
We note that the definition of the AWT does not include the scale $s=0$.
Corollary \ref{lemma:continuity_extension} shows that, under Assumption \ref{assumption:spectral_density},
the continuity of $W_{\Phi}$ can be extended to the boundary $s=0$ by defining $W_{\Phi}(\cdot,0)=0$.
The proof of Corollary \ref{lemma:continuity_extension} is provided in Section \ref{sec:proof:lemma:continuity_extension}.

\section{Definition of ridges}\label{sec:def:ridge}
We start the definition of ridges by recalling existing results when the input signal is noise-free.
With AHM and the slowly varying condition specified in (\ref{AMFM_signal}), we have the following approximation for the AWT of $f$ fulfilling the AHM \cite[Estimate 3.5]{2011synchrosqueezed}:
\begin{align}
\left|W_{f}(t,s)\right| 
\label{relation:scalogram_ridge0}
=&\left|\frac{1}{2}\overset{M}{\underset{m=1}{\sum}}A_{m}(t)e^{i2\pi \phi_{m}(t)}\widehat{\psi}\left(2\pi s\phi'_{m}(t)\right)+\epsilon E(t,s)\right|\,,
\end{align}
where $\epsilon$ is a sufficiently small constant,
$$0\leq E(t,s)\leq s\left[\mathcal{I}_1\sum_{m=1}^M|\phi'_m(t)|+\frac{1}{2}\mathcal{I}_2s\sum_{m=1}^M[M_m''+|A_m(t)||\phi'_m(t)|]+\frac{1}{6}\mathcal{I}_3s^2\sum_{m=1}^{M}M''_m|A_m(t)| \right],
$$
$\mathcal{I}_k:=\int_{\mathbb{R}} |u|^k|\psi(u)|du,\ k=1,2,3$,
and $M_m'':=\underset{t\in \mathbb{R}}{\sup}|\phi_m''(t)|<\infty$.
In other words, when $\epsilon>0$ is sufficiently small, we have
\begin{equation}\label{relation:scalogram_ridge}
|W_{f}(t,s)|\approx \frac{1}{2} \left|\sum_{m=1}^MA_m(t)e^{i2\pi \phi_m(t)}\widehat{\psi}(2\pi s\phi'_m(t))\right|\,.
\end{equation}

To analyze the ridges in the scalogram of $f$, we make the following observations.
First, if the amplitude of the $m$-th IMT function  in $f$, where $m\in \{1,\ldots,M\}$, is very small, its contribution to the scalogram may be buried by spectral leakage caused by the time-varying nature of the frequencies and amplitudes. This can make detecting the $m$-th IF, $\phi'_{m}$, from the scalogram challenging. Conversely, if one component has a very large amplitude, its spectral leakage may dominate the scalogram, potentially masking the contributions of other components.
Second, if the spectral energy of the mother wavelet $\psi$ is spread over a broad frequency range rather than being tightly concentrated, the components in (\ref{relation:scalogram_ridge}) may overlap significantly, making it harder to isolate individual components or accurately infer the IFs.
To build intuition about ridge definitions for noise-free AMH, we focus on the scenario where the impact of spectral leakage in the scalogram of the noise-free AHM is negligible. This holds, for example, when the support of $\hat{\psi}$ is concentrated around $[1-B,1+B]$, where $0<B<\delta/(1+\delta)$ and
$\delta$ is the constant mentioned in the spectral separation condition (\ref{assumption:spectral_separation}).
Under this assumption, the regions where the values of the functions $\{s\mapsto |\widehat{\psi}\left(2\pi s\phi'_{m}(t)\right)|\}_{m=1,2,\ldots,M}$ are non-negligible are nearly non-overlapping.
Consequently, (\ref{relation:scalogram_ridge}) can be further approximated as follows:
\begin{align}\label{relation:scalogram_ridge2}
\left|W_{f}(t,s)\right|
\approx
\frac{1}{2}\overset{M}{\underset{m=1}{\sum}}A_{m}(t)\left|\widehat{\psi}\left(2\pi s\phi'_{m}(t)\right)\right|.
\end{align}

Denote
\begin{align}\label{def:centralfrequency_psi}
\omega_{\psi}= \frac{1}{2\pi}\arg\underset{\omega>0}{\max} |\widehat{\psi}(\omega)|,
\end{align}
which is a singleton under Assumption \ref{assumption:boundedness:psi}.
By (\ref{relation:scalogram_ridge2}), for any $t\in \mathbb{R}$, the local maxima of $\left|W_{f}(t,\cdot)\right|$ approximately
occur at $\frac{\omega_{\psi}}{\phi'_{1}(t)}, \frac{\omega_{\psi}}{\phi'_{2}(t)},\ldots, \frac{\omega_{\psi}}{\phi'_{M}(t)}$, which are inversely related to the IFs $\{\phi'_{1}(t),\phi'_{2}(t),\ldots,\phi'_{M}(t)\}$.
The observation above leads to the following definition of ridges \cite[Definitions 2.1 and 2.2]{lilly2010analytic} (see also \cite[Section B]{meignen2023new}).

\begin{Definition}\label{definition:ridgepoint}
Suppose $f$ satisfies AHM (\ref{AMFM_signal}).  A ridge point of $S_{f}(t,s)$ is a time/scale pair $(t,s)$ satisfying the two conditions:
\begin{equation*}
(R1)\ \frac{\partial}{\partial s} S_{f}(t,s) = 0
\quad
\mbox{and}\quad
(R2)\ \frac{\partial^{2}}{\partial s^{2}} S_{f}(t,s) < 0\,;
\end{equation*}
that is $S_f(t,s)$ is a local maximum along the scale axis.
A {\em ridge curve} in $S_f(t,s)$ 
is the graph of a function of $t$
such that  the conditions (R1) and (R2) are satisfied along the curve, and
(R3) the function is continuously differentiable.
\end{Definition}

Note that for a signal fulfilling the AHM, there might be multiple ridge curves, and a ridge curve is a set of ridge points fulfilling the regularity condition (R3).
To utilize this definition, consider
\begin{align}\label{def:s_{f}(t)}
s_{f}(t)=\arg\underset{s>0}{\max}\ S_{f}(t,s)\,.
\end{align}
When $f$ satisfies the AHM (\ref{AMFM_signal}) with $M=1$, $s_{f}$ is a single-valued function.
By \eqref{relation:scalogram_ridge2},
\begin{align*}
s_{f}(t)\approx 
\arg\underset{s>0}{\max}\left|\widehat{\psi}\left(2\pi s\phi'(t)\right)\right|
=\frac{\omega_{\psi}}{\phi'(t)}\,;
\end{align*}
that is, the ridge captures the IF information.
See Figure \ref{fig:wavelet_potential_ridge_Surface plot} for an example of the graph of $s_{f}$
and its relationship with the IF $\phi'(t)$. Based on the proof technique using the implicit function theorem presented in \cite{colominas2023iterative}, the functions $s_{f}(t)$ is continuously differentiable and hence constitutes a ridge. Note that \cite{colominas2023iterative} analyzes the ridge of the spectrogram generated by the short-time Fourier transform. Since a similar proof can be carried out to the scalogram, we skip the details here.

When $f$ satisfies the AHM (\ref{AMFM_signal}) with $M\geq 2$, \eqref{def:s_{f}(t)} might be problematic.
First, due to the time-varying amplitude, $s_f(t)$ might alternately capture the IF of different IMT functions. For example, consider the case where $M=2$, $\phi_m'(t)=\xi_m>0$ are constants, and $\xi_2>2\xi_1$.
If $A_1(t)>2A_2(t)$ for $t<-1$ and $A_2(t)> 2A_1(t)$ for $t>1$, then $s_f(t)\approx \omega_{\psi}/\phi'_1(t)$ when $t<-1$ and $s_f(t)\approx \omega_{\psi}/\phi'_2(t)$ when $t>1$.
On the other hand, when the spectral leakage from other IMT functions of $f$ is smaller than $A_m(t)|\widehat{\psi}(2\pi s\phi'_m(t))|$,
the cross-section of the scalogram at time $t$ exhibits $M$ local maxima.
In this case, to utilize the definition of ridges, consider
\begin{align}\label{prior_information}
\frac{\omega_{\psi}}{\phi'_{m}(t)} \in B_{m}(t):= \left[\underline{b}_{m}(t), \overline{b}_{m}(t)\right]\ \textup{for}\ 1\leq m\leq M,
\end{align}
where $\underline{b}_{1},\overline{b}_{1},\underline{b}_{2},\overline{b}_{2},\ldots,\underline{b}_{M},\overline{b}_{M}$ are continuous on $\mathbb{R}$ and
satisfy
\begin{align*}
0<\underline{b}_{M}(t)<\overline{b}_{M}(t)\leq\cdots\leq\underline{b}_{2}(t)<\overline{b}_{2}(t)\leq\underline{b}_{1}(t)<\overline{b}_{1}(t)< \infty\,.
\end{align*}
In other words, we segment the time-scale domain into pieces using the IF.
Then, define
\begin{align*}
s_{f,m}(t) = \arg\underset{s\in B_{m}(t)}{\max} S_{f}(t,s).
\end{align*}

Similarly, by the same argument as in \cite{colominas2023iterative}, we can prove that $s_{f,m}(t)$ are continuously differentiable, and hence they are ridge curves by definition.
Hereafter, we use the graphs of $s_{f}$ and $\{s_{f,m}\}_{m=1,\ldots,M}$ to represent the ridge curves of the noise-free function satisfying the AHM,
and we refer to them collectively as the {\em ridge curves} for simplicity.


\begin{Remark}
In this work, we focus on the local properties of ridge curves of the noisy signal $Y$, which will be defined in Section \ref{sec:def:ridge:Y},
and on their deviation from those of the noise-free signal $f$.
Accordingly, we suppose that the scalogram $S_f(t,\cdot)$ of the noise-free signal exhibits
$M$ well-separated local maxima at each time $t$ with known locations.
This assumption ensures the existence of the interval functions
$\underline{b}_m$ and $\overline{b}_m$.

The choice of $\underline b_m$ and $\overline b_m$ is not unique.
When additional information about $f$ is available, the width of the ring-shaped regions
$\{(t,s)\mid \underline b_m(t)\leq s \leq \overline b_m(t)\}_{m=1}^{M}$ may be reduced.
In practice, the IFs of $f$, and hence the intervals $B_m(t)$, must be estimated.
An example of an algorithm addressing this estimation problem is given in \cite{Laurent_Meignen_2021}.
This algorithm, jointly with its theoretical analysis, falls in Steps \texttt{III} and \texttt{IV}, and lies beyond the scope of this paper.

Conditions ensuring the existence of well-separated local maxima in time-frequency representations have been studied
for specific signal classes, such as purely harmonic modes or parallel linear chirps, primarily in the context of the
spectrogram; see, for example, \cite{meignen2021study}.
Extending such results to the situation that two IMT functions have close IFs, or even crossover IFs, constitutes an important but challenging problem. Handling this challenge is related to analyzing the spectral interference behavior of CWT. It is also important to study the case when the IF or AM fail the slowly varying condition. Since these challenges are out of the scope of this paper, we leave them to our future work.
\end{Remark}

\subsection{Definition of ridges for noise-affected AHM}\label{sec:def:ridge:Y}

We proceed to the main focus of this paper: define and analyze ridges when the signal is noisy. Motivated by Definition \ref{definition:ridgepoint} and real data experience, we consider the following definition.

\begin{Definition}\label{definition:ridgepointNOISY}
Denote $S_{Y}(t,s)=|W_{Y}(t,s)|^{2}.$
For the noise-affected AHM (\ref{noisycos}) with $M\geq1$, and analogous to (\ref{def:s_{f}(t)}),  we define
\begin{align}\label{proxy:s_Y}
s_{Y}(t)= \arg \underset{s>0}{\max}\ S_{Y}(t,s)
\end{align}
for any $t\in \mathbb{R}$, and refer the graph of $s_{Y}$ as a ridge curve of $Y$,
noting that $s_{Y}$ may be a multivalued function at some time points.

For the noise-affected AHM (\ref{noisycos}) with $M\geq1$,
we also define
\begin{align}\label{proxy:s_Ym}
s_{Y,m}(t)= \arg \underset{s\in B_{m}(t)}{\max}\ S_{Y}(t,s)
\end{align}
for any $t\in \mathbb{R}$ and $m\in\{1,2,\ldots,M\}$, and refer to the graph of $s_{Y,m}$ as the $m$th ridge curve of $Y$.
\end{Definition}

\begin{figure}[hbt!]
  \centering
  \includegraphics[width=0.99\textwidth]{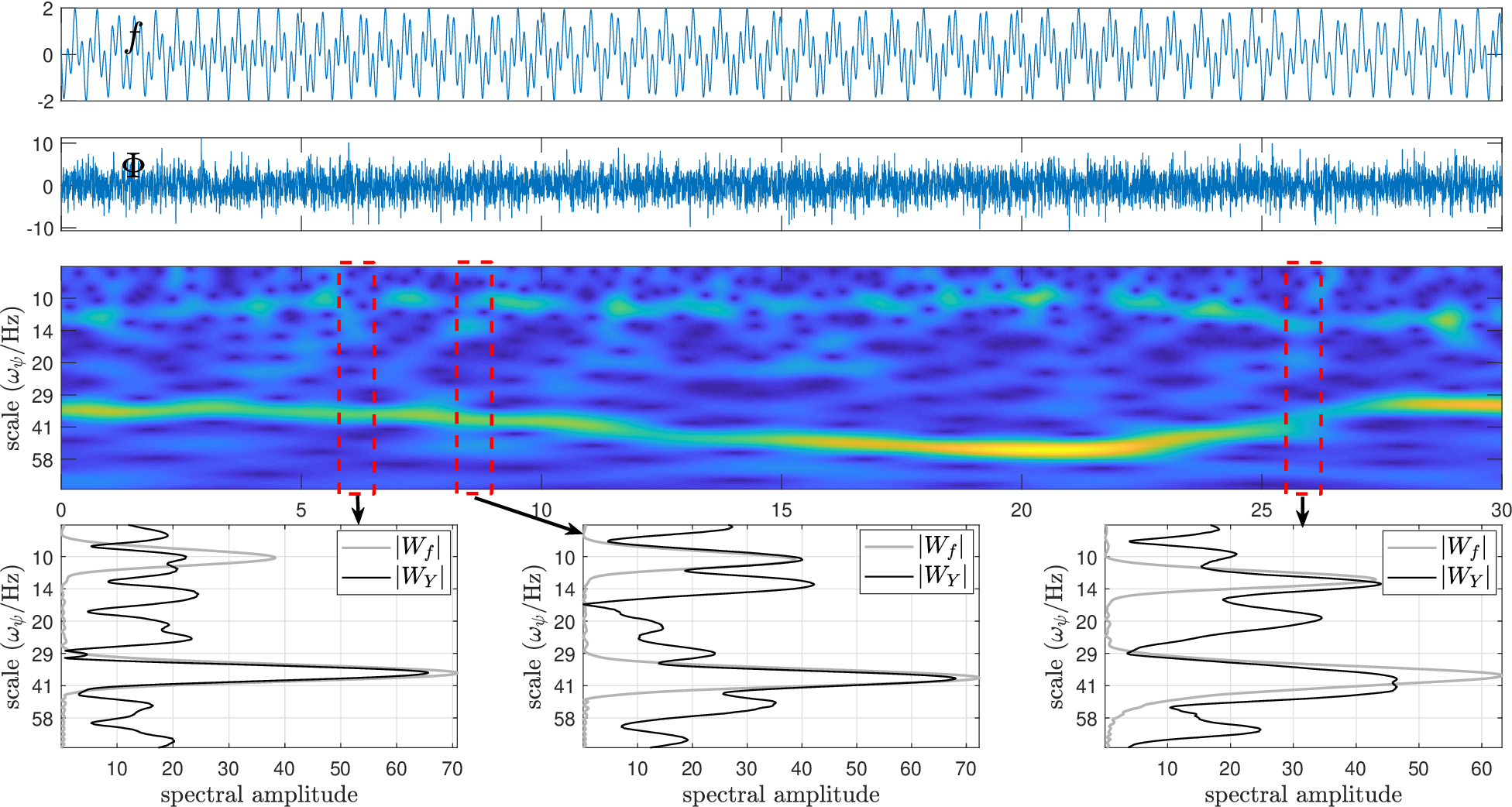}
 \caption{Example of the complex modulus of the AWT for a 30-second noise-contaminated signal composed of two frequency-modulated sinusoids.
The first row displays the frequency-modulated signal $f$, while the second row shows the noise, a sample from the stationary Gaussian process $\Phi$, identical to that in Figure \ref{fig:wavelet_potential_ridge_Surface plot}. The third row illustrates the complex modulus of the AWT, i.e., the spectral amplitude, of the noise-contaminated signal $Y$, where $Y=f+\Phi$. In the last three columns, the black curves represent the spectral amplitude of $Y$ at specific times, corresponding to cross-sections of the image in the third row. For comparison, the spectral amplitudes of $f$ are shown in gray.
For the single-component case, see Figure \ref{fig:wavelet_ridge:one} in the appendix.}
\label{fig:wavelet_ridge}
\end{figure}

Note that compared with the ridge curve defined in Definition \ref{definition:ridgepoint} when the signal is noise-free, in the above definition we do not specify the regularity of the ridge curve. Instead, we need to explore its regularity.
To this end, observe that unlike $|W_{f}(t,\cdot)|$ shown in Figure \ref{fig:wavelet_potential_ridge_Surface plot}, $|W_{Y}(t,\cdot)|$ has multiple local maxima even in the case where $M=1$, as illustrated in Figure \ref{fig:wavelet_ridge}. We further show the scatter plot of the local maximizers of $|W_{Y}(t,\cdot)|$ for any $t\in \mathbb{R}$ in Figure \ref{fig:wavelet_potential_ridge}. We find many paths formed by the local maximizers of $|W_{Y}(t,\cdot)|$ for each $t\in \mathbb{R}$, and they are piecewise continuous. Moreover, we find that near the ``true'' IF, the associated local maximum that we expect to be the ridge occasionally jumps; for example, around the 6th, 8th, and 26th seconds in the top panel of Figure \ref{fig:wavelet_potential_ridge}. This motivates us to analyze the graphs of $s_{Y}$ and $\{s_{Y,m}\}_{m=1}^{M}$, as well as
the difference between $s_{f}$ (and  $\{s_{f,m}\}_{m=1}^{M}$)
and $s_{Y}$ (and $\{s_{Y,m}\}_{m=1}^{M}$).

\section{Main results}\label{sec:mainresult}

In this section, we present our main results. The first part is about the singleton property of the set-valued functions $s_{Y}$ (or $s_{Y,1},\ldots,s_{Y,M}$ when $M>1$), as well as the hemicontinuity and local $C^{1}$-smoothness properties of their graphs. The second part is the deviation of $s_Y$ (or $s_{Y,1},\ldots,s_{Y,M}$ when $M>1$) from $s_f$ (or $s_{f,1},\ldots,s_{f,M}$ when $M>1$).

\subsection{Local properties of ridge curves}

\begin{Theorem}\label{prop:unique_argmax}
Suppose $f\in L^2$, $\Phi$ is a stationary Gaussian process satisfying Assumptions \ref{assumption:Gaussian} and \ref{assumption:spectral_density},
and $\psi$ is a mother wavelet satisfying Assumption \ref{assumption:boundedness:psi}.
Then, we have
\begin{align*}
\mathbb{P}\left(s_{Y}(t)\ \textup{is a singleton at any}\  t\in \mathbb{Q}\right)=1,
\end{align*}
where $\mathbb{Q}$ denotes the set of rational numbers.
Moreover, when $f$ satisfies the AHM defined in (\ref{AMFM_signal}) with $M\geq 1$,
we have
\begin{align*}
\mathbb{P}\left(s_{Y,m}(t)\ \textup{is a singleton}\ \textup{for any}\  t\in \mathbb{Q}\ \textup{and}\ m\in\{1,\ldots,M\}\right)=1.
\end{align*}
\end{Theorem}

We note that Theorem \ref{prop:unique_argmax} holds regardless of whether $f\equiv 0$ or not, and in general, the singleton property of $s_Y$ does not require the AHM assumption. The proof of Theorem \ref{prop:unique_argmax} is provided in Section \ref{sec:proof:prop:unique_argmax}.
In addition to investigating the singleton property of the sets $ s_{Y}(t)$, or $s_{Y,1}(t),\ldots,s_{Y,M}(t)$ for $f$ satisfying the AHM, at $\mathbb{Q}$, we are also interested in the behavior of these sets for the remaining $t\in \mathbb{R}\backslash \mathbb{Q}$.

\begin{Theorem}\label{lemma:compact&hemiconti}
Suppose $f\in L^2$,  $\Phi$ is a stationary Gaussian process satisfying Assumptions \ref{assumption:Gaussian} and \ref{assumption:spectral_density},
and $\psi$ is a mother wavelet satisfying Assumption \ref{assumption:boundedness:psi}.
With probability one,
the set-valued function $ s_{Y}$ is upper hemicontinuous on $\mathbb{R}$.
That is, almost surely, for any $t\in \mathbb{R}$ and any open set $U$ with $ s_{Y}(t)\subset U$,
there exists $\delta>0$ such that
$ s_{Y}(t')\subset U$ for any $t'\in (t-\delta,t+\delta)$.
Moreover, if $f$ is an adaptive harmonic signal as defined in (\ref{AMFM_signal}) with $M\geq 1$,
then almost surely,
the set-valued functions $s_{Y,1},\ldots,s_{Y,M}$, defined in (\ref{proxy:s_Ym}), are also upper hemicontinuous on $\mathbb{R}$.
\end{Theorem}

\begin{figure}[hbt!]
\centering
  \includegraphics[width=0.99\textwidth]{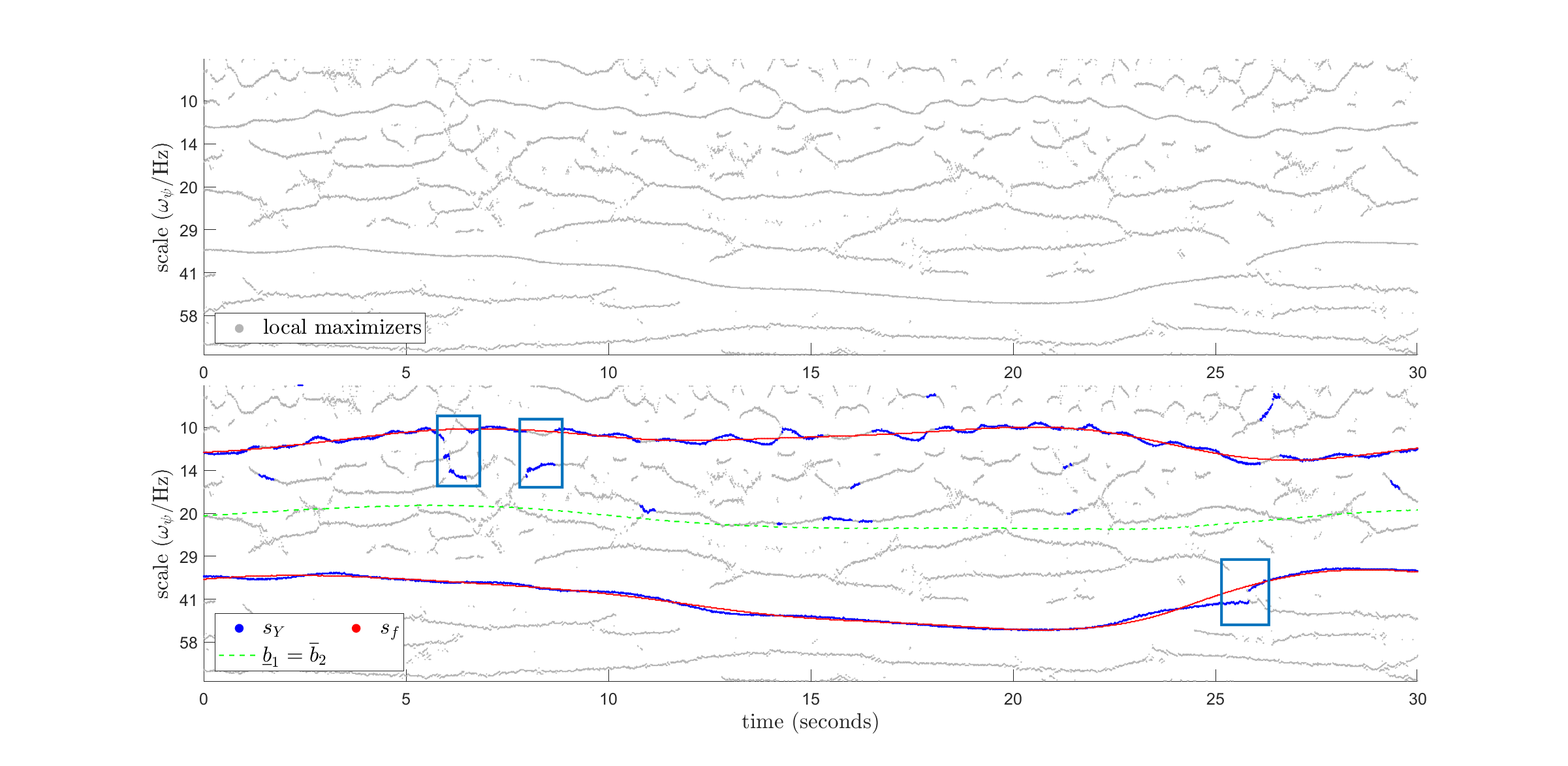}
\caption{Example of ridge points (local maxima) and a ridge curve $ s_{Y}$, where $Y = f+\Phi$, and $s_{f}$. The signal $f$ and the realization of $\Phi$ are identical to those in Figure \ref{fig:wavelet_potential_ridge_Surface plot}. The blue boxes in the second row highlight the discontinuities in the ridge curve $s_{Y}$ caused by the presence of noise.
For the single-component case, refer to Figure \ref{fig:wavelet_potential_ridge:one} in the appendix.
  \label{fig:wavelet_potential_ridge}}
\end{figure}

The upper hemicontinuity of $ s_{Y}$ implies that as $t'$ approaches
$t$, the maximizers of $S_{Y}(t',\cdot)$ get arbitrarily close to some of the maximizers of $S_{Y}(t,\cdot)$, as illustrated in Figure \ref{fig:wavelet_potential_ridge}.
The proof of Theorem \ref{lemma:compact&hemiconti} is provided in Section \ref{sec:proof:lemma:compact&hemiconti}, which relies on Berge's maximum theorem
\cite[Theorem 3.4]{hu1997handbook} (see also \cite[p. 570]{aliprantis2006infinite} and Lemma \ref{lemma:Berge} in the appendix).

Theorems \ref{prop:unique_argmax} and \ref{lemma:compact&hemiconti} imply that, with probability one,
the restriction of $ s_{Y}$ on the dense set $\mathbb{Q}$, denoted by $ s_{Y}|_{\mathbb{Q}}$, is continuous.
By \cite[Lemma (4.3.16)]{MR1039321} (see also Lemma \ref{lemma:continuous_dense} in the appendix),
$ s_{Y}|_{\mathbb{Q}}$
is extendable to a continuous function defined on a $G_\delta$-set containing
$\mathbb{Q}$. Here, a $G_{\delta}$-set refers to a subset of $\mathbb{R}$ that is a countable intersection of open sets.
Next, we provide more insights into this $G_{\delta}$-set.

\begin{Corollary}\label{prop:C1_argmax}
Assume that Assumptions \ref{assumption:Gaussian}-\ref{assumption:spectral_density} hold.
For any $t_{0}\in \mathbb{Q}$, if
\begin{align}\label{assumption:atomless}
\mathbb{P}\left(\frac{\partial^{2}S_{Y}}{\partial s^{2}}\left(t_{0}, s_{Y}(t_{0})\right)= 0\right)=0,
\end{align}
then, almost surely, there exists a unique continuously differentiable function $r$ defined on
$(t_{0}-\delta,t_{0}+\delta)$ for some $\delta>0$
such that
\begin{itemize}
\item $r(t) =  s_{Y}(t)$ for any $t\in  (t_{0}-\delta,t_{0}+\delta)\cap \mathbb{Q}$, and
\item
$r(t) \in  s_{Y}(t)$ for any $t\in  (t_{0}-\delta,t_{0}+\delta)\cap \mathbb{Q}^{c}$.
\end{itemize}
\end{Corollary}

\begin{figure}[htb!]
\centering
\subfigure[][$f(t)=\cos(2\pi t)$]
{\includegraphics[width=0.495\textwidth]{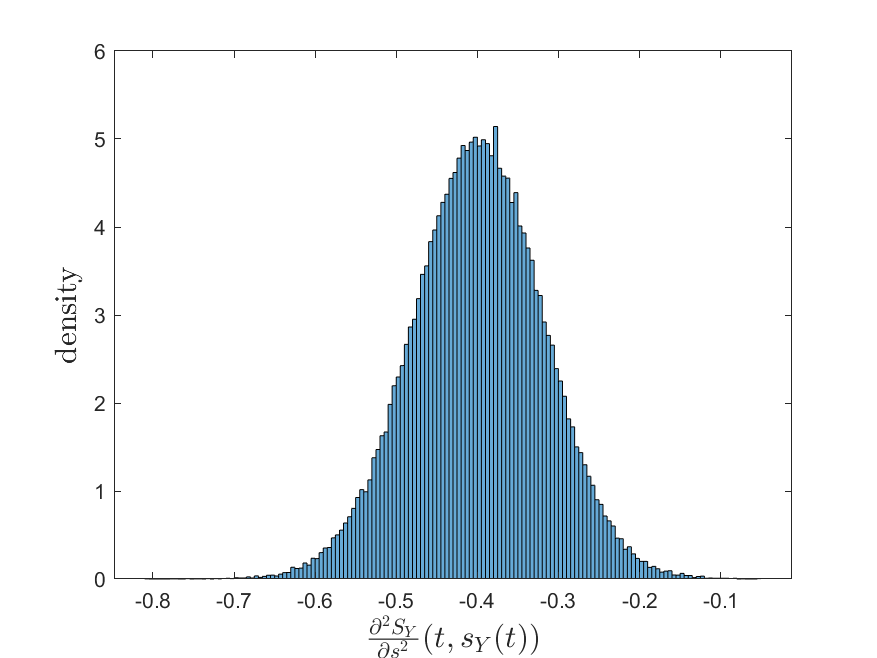}}
\subfigure[][$f(t)=\cos(2\pi t^2)$]
{\includegraphics[width=0.495\textwidth]{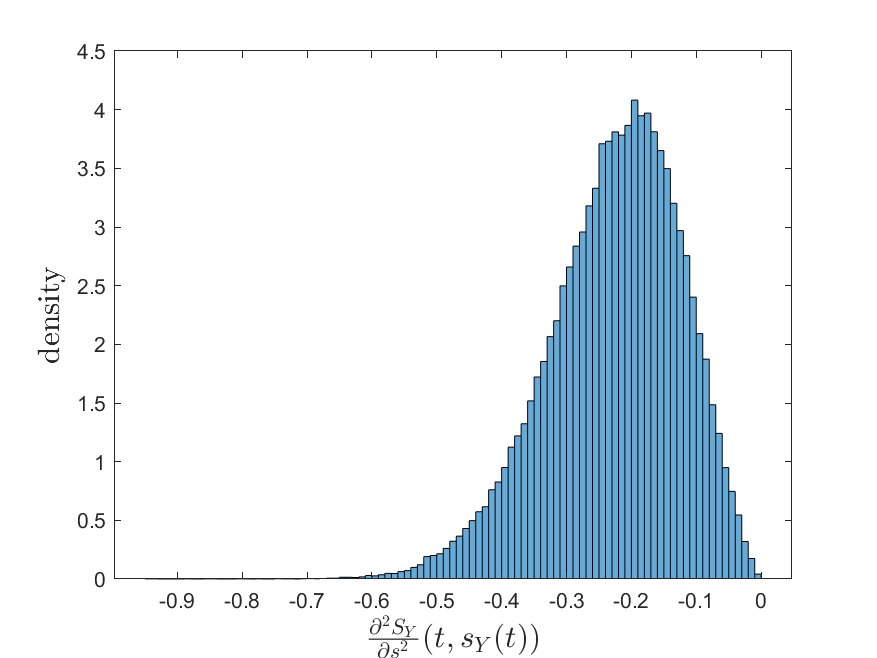}}
\caption{Histogram of the random variable $\frac{\partial^{2}S_{Y}}{\partial s^{2}}\left(t, s_{Y}(t)\right)$, where $Y(t) = f(t)+\Phi(t)$ and $\Phi$ is a stationary Gaussian process. The time variable is fixed throughout the experiment, and the histogram is based on $10^5$ realizations of $\Phi$.
 \label{fig:histogramD2}}
\end{figure}

The proof of Corollary \ref{prop:C1_argmax} is provided in
Section \ref{sec:proof:prop:C1_argmax}. Corollary \ref{prop:C1_argmax} provides additional insight into the nature of the $G_{\delta}$-set,
on which $s_{Y}$ is continuous, under the assumption (\ref{assumption:atomless}).
This assumption is supported by the histogram of $\frac{\partial^{2}S_{Y}}{\partial s^{2}}\left(t, s_{Y}(t)\right)$
shown in Figure \ref{fig:histogramD2}.
Corollary \ref{prop:C1_argmax} also explains why the ridge curves observed in the numerical experiments consist of piecewise continuous curves. See Figure \ref{fig:wavelet_potential_ridge} for an illustration. Note that while it is visually possible to identify if the ridge points form a continuous path, the local maxima plot further highlights and emphasizes any discontinuities.

Furthermore, if the signal $f$ fulfills the AHM, we have the following property of the paths $\left\{(t,s_{Y,m}(t))\mid t\in \mathbb{R}\right\}_{m=1,\ldots,M}$, which is parallel to Corollary \ref{prop:C1_argmax}.

\begin{Corollary}\label{prop:unique_argmax:m}
Assume that Assumptions \ref{assumption:Gaussian}-\ref{assumption:spectral_density} hold.
For any $t_{0}\in \mathbb{Q}$ and $m\in\{1,\ldots,M\}$, if
$s_{Y,m}(t_{0})$ lies in the interior of $B_{m}(t_{0})$ and
\begin{align*}
\mathbb{P}\left(\frac{\partial^{2}S_{Y}}{\partial s^{2}}\left(t_{0}, s_{Y,m}(t_{0})\right)= 0\right)=0,
\end{align*}
then, almost surely,
there exists a unique continuously differentiable function $r_{m}$ defined on
$(t_{0}-\delta,t_{0}+\delta)$ for some $\delta>0$ such that
\begin{itemize}
\item $r_{m}(t) \in B_{m}(t)$ for any $t\in  (t_{0}-\delta,t_{0}+\delta)$,
\item $r_{m}(t) =  s_{Y,m}(t)$ for any $t\in  (t_{0}-\delta,t_{0}+\delta)\cap \mathbb{Q}$, and
\item
$r_{m}(t) \in  s_{Y,m}(t)$ for any $t\in  (t_{0}-\delta,t_{0}+\delta)\cap \mathbb{Q}^{c}$.
\end{itemize}

\end{Corollary}

Based on our literature review, the analysis of the uniqueness of maxima for a random field determined by AWT along the scale axis, as well as the behavior of the locations of these maxima along another axis, has rarely been discussed.

\subsection{Deviation of ridge points caused by noise}

We continue to quantify the probability that $|s_{Y}(t)-s_{f}(t)|$ and $|s_{Y,m}(t)-s_{f,m}(t)|$
exceed a given threshold, where $m=1,\ldots,M$.

\begin{Theorem}\label{lemma:conditional}
Suppose that the stationary Gaussian process $\Phi$ satisfies Assumptions \ref{assumption:Gaussian} and \ref{assumption:spectral_density}, and the mother wavelet $\psi$ satisfies Assumption \ref{assumption:boundedness:psi}.
Fix $t\in \mathbb{R}$. For any interval $I$ containing
$s_{f}(t)$,
denote
\begin{align}
\triangle_{I} =  |W_{f}(t,s_{f}(t))|-\underset{s\in I^{c}}{\max}|W_{f}(t,s)|, \ \
\mu_{I}  =\mathbb{E}\left[|W_{\Phi}(t,s_{f}(t))|\right] +\mathbb{E}[\underset{s\in I^{c}}{\max}|W_{\Phi}(t,s)|],\label{def:mu_ab}
\end{align}
and
\begin{align}\label{def:sigma_ab}
\sigma_{I} = &
\sqrt{\mathbb{E}\left[|W_{\Phi}(t,s_{f}(t))|^{2}\right]}
+\underset{s\in I^{c}}{\max}\sqrt{\mathbb{E}\left[|W_{\Phi}(t,s)|^{2}\right]}.
\end{align}
If $\triangle_{I}>\mu_{I}$, then
\begin{align}\label{ineq:sYinI}
\mathbb{P}\left( s_{Y}(t)\in I\right)\geq
1-\textup{exp}\left(-\frac{(\triangle_{I}-\mu_{I})^{2}}{\sigma_{I}^{2}}\right),
 \end{align}
where $ s_{Y}$ is defined in (\ref{proxy:s_Y}).
\end{Theorem}

Theorem \ref{lemma:conditional} shows that if the peak of $|W_{f}(t,\cdot)|$ is distinct (i.e. $\triangle_{I}$ is ``large'', or $\epsilon E(t,\cdot)$ in \eqref{relation:scalogram_ridge0} is sufficiently small) and much higher than the peaks of $\{|W_{\Phi}(t,s)|\mid s\in I^{c}\}$ (i.e. $\triangle_{I}>\mu_{I}$), then the deviation of the ridge point $s_{f}(t)$ caused by noise is controlled with high probability.
Deviations in the ridge point of $f$ may also arise if the spectral component of $f$ at the scale $s_{f}(t)$ is dominated by the spectral component of the noise at the same scale. Hence, the moments of the wavelet coefficients of the noise $\Phi$ at the scale $s_{f}(t)$, including $\mathbb{E}\left[|W_{\Phi}(t,s_{f}(t))|\right]$ and $\mathbb{E}\left[|W_{\Phi}(t,s_{f}(t))|^{2}\right]$,
are contained in $\mu_{I}$ and $\sigma_{I}$.
Most components in (\ref{def:mu_ab}) and  (\ref{def:sigma_ab}) can be calculated
as follows:
\begin{align*}
\mathbb{E}\left[|W_{\Phi}(t,s_{f}(t))|\right] = \frac{1}{2}\pi^{\frac{1}{2}}\left[\int_{0}^{\infty}|\widehat{\psi}\left(s_{f}(t)\lambda\right)|^{2}p(\lambda)d\lambda\right]^{\frac{1}{2}}
\end{align*}
and
\begin{align*}
\sigma_{I} =
\sqrt{\int_{0}^{\infty}|\widehat{\psi}\left(s_{f}(t)\lambda\right)|^{2}p(\lambda)d\lambda}
+\underset{s\in I^{c}}{\max}\sqrt{\int_{0}^{\infty}|\widehat{\psi}\left(s\lambda\right)|^{2}p(\lambda)d\lambda}.
\end{align*}

The proof of Theorem \ref{lemma:conditional} is
presented in Section \ref{sec:proof:lemma:conditional}.
To prove Theorem \ref{lemma:conditional} and derive an upper bound for the term
$\mathbb{E}[\underset{s\in I^{c}}{\max}|W_{\Phi}(t,s)|]$
in (\ref{def:mu_ab}),
we need to extend the classical Gaussian concentration inequality for Lipschitz functions of real-valued Gaussian elements \cite{kuhn2023maximal}, Dudley's theorem, and the Borell-TIS inequality \cite{talagrand2014upper,vershynin2018high},
to account for complex-valued Gaussian elements arising from the AWT of $\Phi$.
To keep the focus on the discussion of ridge deviation,
these extensions are provided in Section \ref{sec:Dudley}.

 \begin{Remark}
We remark that the lower bound in Theorem~\ref{lemma:conditional} and the upper bound in Theorem~\ref{mainresult:deviation} are presented from a theoretical perspective, aiming to explain how key factors--such as the peak location $s_{f}(t)$ on the scalogram of the clean signal $f$, the prominence $\triangle_{I}$ of the wavelet ridge of $f$, and the terms $\mu_{I}$ and $\sigma_{I}$ related to the spectral energy distribution of the noise--influence the probability that the ridge deviation exceeds a given threshold. Although these quantities are not accessible from noisy observations, their presence in the bounds highlights the underlying mechanisms governing ridge stability.
\end{Remark}

\begin{Corollary}\label{corollary:sYm_in_I}
Suppose that $f$ satisfies the AHM, the stationary Gaussian process $\Phi$ satisfies Assumptions \ref{assumption:Gaussian} and \ref{assumption:spectral_density}, and the mother wavelet $\psi$ satisfies Assumption \ref{assumption:boundedness:psi}.
Fix $t\in \mathbb{R}$.
For any $m\in\{1,\ldots,M\}$ and any interval $I_{m}\subset B_{m}(t)$ containing
$s_{f,m}(t)$,
denote
\begin{align}\label{def:triangle_mu_sigma_I}
\begin{aligned}
\triangle_{I_{m}} =&  |W_{f}(t,s_{f,m}(t))|-\underset{s\in B_{m}(t)\setminus I_{m}}{\max}|W_{f}(t,s)|,
\\
\mu_{I_{m}}=&\mathbb{E}\left[|W_{\Phi}(t,s_{f,m}(t))|\right] +\mathbb{E}\left[\underset{s\in B_{m}(t)\setminus I_{m}}{\max}|W_{\Phi}(t,s)|\right],
\end{aligned}
\end{align}
and
$$\sigma_{I_{m}} =
\sqrt{\mathbb{E}\left[|W_{\Phi}(t,s_{f,m}(t))|^{2}\right]}
+\underset{s\in B_{m}(t)\setminus I_{m}}{\max}\sqrt{\mathbb{E}\left[|W_{\Phi}(t,s)|^{2}\right]}.$$

If $\triangle_{I_{m}}>\mu_{I_{m}}$, then
\begin{align*}
\mathbb{P}\left( s_{Y,m}(t)\in I_{m}\right)\geq
1-\textup{exp}\left(-\frac{(\triangle_{I_{m}}-\mu_{I_{m}})^{2}}{\sigma_{I_{m}}^{2}}\right)
 \end{align*}
where $ s_{Y,m}$ is defined in (\ref{proxy:s_Ym}).
\end{Corollary}

\begin{figure}
  \centering
  \includegraphics[width=0.99\textwidth]{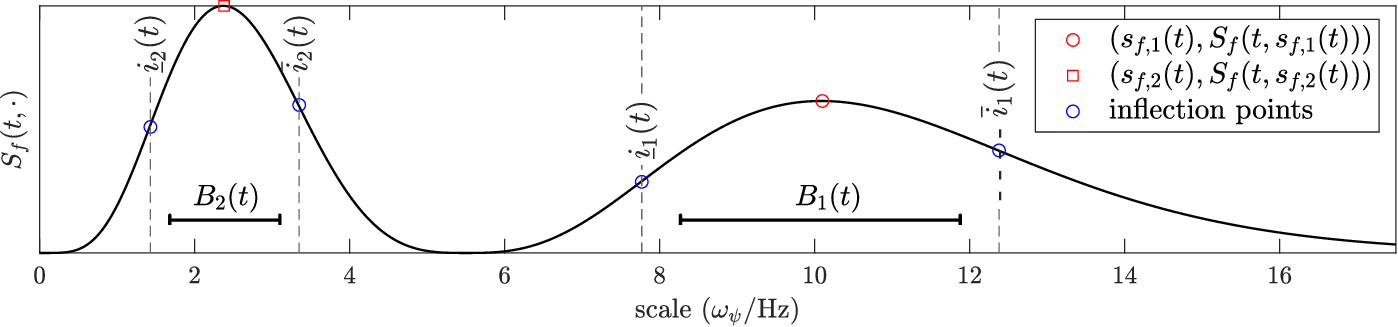}
  \caption{The relationship between the scale region $B_{m}(t)$ and the interval $[\underline{\mathfrak{i}}_{m}(t),\overline{\mathfrak{i}}_{m}(t)]$ considered in Theorem \ref{mainresult:deviation}.}\label{fig:inflection_interval}
\end{figure}

For every $m\in\{1,2,\ldots,M\}$ and $t\in \mathbb{R}$, let $\underline{\mathfrak{i}}_{m}(t)$ and $\overline{\mathfrak{i}}_{m}(t)$ be the inflection points of $S_{f}(t,\cdot)$ closest to $s_{f,m}(t)$
with $\underline{\mathfrak{i}}_{m}(t)<s_{f,m}(t)<\overline{\mathfrak{i}}_{m}(t)$, as depicted in Figure \ref{fig:inflection_interval}.
As shown in Figures \ref{fig:wavelet_ridge} and \ref{fig:wavelet_potential_ridge},
local maxima of $S_{Y}(t,\cdot)$ are likely to occur within
neighborhoods of the local maxima of $S_{f}(t,\cdot)$.
Hence, for the function $s_{Y,m}$ defined in (\ref{proxy:s_Ym})
with
$B_{m}(t)\subset (\underline{\mathfrak{i}}_{m}(t),\overline{\mathfrak{i}}_{m}(t))$, we are also interested in the probability that $| s_{Y,m}(t)-s_{f,m}(t)|$ exceeds a threshold
conditioned on the event $ s_{Y,m}(t)\in B^{\circ}_{m}(t)$,
where $B^{\circ}_{m}(t)$ is the interior of $B_{m}(t)$.

Before presenting our final main result, we introduce the following notation:
For every $m\in\{1,2,\ldots,M\}$ and $t\in \mathbb{R}$, let
\begin{align}\notag
  &\mu_{1,m} = \mathbb{E}\left[\underset{s\in B_{m}(t)}{\max}\left|\frac{\partial W_{f}}{\partial s}(t,s)W_{\Phi}(t,s)\right|\right],
  & \quad &
\sigma_{1,m}^{2} =
\underset{s\in  B_{m}(t)}{\max} \left|\frac{\partial W_{f}}{\partial s}(t,s)\right|^{2}\mathbb{E}\left[\left|W_{\Phi}(t,s)\right|^{2}\right],
\\\notag
&\mu_{2,m} = \mathbb{E}\left[\underset{s \in B_{m}(t)}{\max}\left|W_{f}(t,s)\frac{\partial W_{\Phi}}{\partial s}(t,s)\right|\right],
& \quad &
\sigma_{2,m}^{2} = \underset{s\in B_{m}(t)}{\max} \left|W_{f}(t,s)\right|^{2}\mathbb{E}\left[\left|\frac{\partial W_{\Phi}}{\partial s}(t,s)\right|^{2}\right],
\\\label{def:mu1_4}
&\mu_{3,m}=
\mathbb{E}\left[\underset{s\in  B_{m}(t)}{\max}\left|W_{\Phi}(t,s)\right|\right],
& \quad &
\sigma_{3,m}^{2} = \underset{s\in  B_{m}(t)}{\max} \mathbb{E}\left[\left|W_{\Phi}(t,s)\right|^{2}\right],
\\\notag
&\mu_{4,m}=\mathbb{E}\left[\underset{s\in  B_{m}(t)}{\max}\left|\frac{\partial W_{\Phi}}{\partial s}(t,s)\right|\right],
& \quad &
\sigma_{4,m}^{2} = \underset{s\in  B_{m}(t)}{\max} \mathbb{E}\left[\left|\frac{\partial W_{\Phi}}{\partial s}(t,s)\right|^{2}\right],
\end{align}
where
\begin{align*}
\mathbb{E}\left[\left|W_{\Phi}(t,s)\right|^{2}\right]= \int_{0}^{\infty}|\widehat{\psi}(s\lambda)|^{2}p(\lambda)d\lambda,
\ \ \ \mathbb{E}\left[\left|\frac{\partial W_{\Phi}}{\partial s}(t,s)\right|^{2}\right]= \int_{0}^{\infty}
|\lambda [D\widehat{\psi}](s\lambda) |^{2}p(\lambda)d\lambda,
\end{align*}
and $D\widehat{\psi}$ represents the derivative of $\widehat{\psi}$.

\begin{Theorem}\label{mainresult:deviation}
Suppose that $f$ satisfies the AHM, the stationary Gaussian process $\Phi$ satisfies Assumptions \ref{assumption:Gaussian} and \ref{assumption:spectral_density}, and the mother wavelet $\psi$ satisfies Assumption \ref{assumption:boundedness:psi}.
For every $m\in\{1,2,\ldots,M\}$, $t\in \mathbb{R}$, and
$B_{m}(t)\subset (\underline{\mathfrak{i}}_{m}(t),\overline{\mathfrak{i}}_{m}(t)),$
where the interval $B_{m}(t)$ is defined in (\ref{prior_information}),
let
\begin{align}\label{def:lower_bound_L}
L_{m} = \underset{s\in B_{m}(t)}{\min}
\left|\frac{\partial^{2}S_{f}}{\partial s^{2}}(t,s)\right|.
\end{align}
For any $\varepsilon>0$ with
\begin{equation}\label{available_region}
\underset{s\in B_{m}(t)}{\min}|s-s_{f,m}(t)| >\varepsilon>\frac{6}{L_{m}}\max\left\{\mu_{1,m},
\mu_{2,m},\mu_{3,m}^{2},\mu_{4,m}^{2}\right\},
\end{equation}
we have
\begin{align}\notag
&\mathbb{P}\left(|s_{Y,m}(t)-s_{f,m}(t)|>\varepsilon \mid   s_{Y,m}(t)\in B^{\circ}_{m}(t)\right)
\leq
\left[1-\textup{exp}\left(-\frac{(\triangle_{B^{\circ}_{m}(t)}-\mu_{B^{\circ}_{m}(t)})^{2}}
{\sigma_{B^{\circ}_{m}(t)}^{2}}\right)\right]^{-1}
\\\label{main_4_exponent}&\qquad\qquad\times\left\{\overset{2}{\underset{k=1}{\sum}}\textup{exp}\left[-\frac{L_{m}^{2}}{\sigma_{k,m}^{2}}\left(\frac{\varepsilon }{6}-\frac{\mu_{k,m}}{L_{m}}\right)^{2}\right]
+
\overset{4}{\underset{k=3}{\sum}}
\textup{exp}\left[-\frac{L_{m}}{\sigma_{k,m}^{2}}\left(\sqrt{\frac{\varepsilon }{6}}-\frac{\mu_{k,m}}{\sqrt{L_{m}}}\right)^{2}\right]
\right\},
\end{align}
where $\triangle_{B^{\circ}_{m}(t)}$,
$\mu_{B^{\circ}_{m}(t)}$, and $\sigma_{B^{\circ}_{m}(t)}$
are defined in (\ref{def:triangle_mu_sigma_I}) with $I =B^{\circ}_{m}(t)$ and $B_{m}(t)\setminus I = \partial B_{m}(t)$, i.e., the endpoints of the interval $B_{m}(t)$.
\end{Theorem}

The proof of Theorem \ref{mainresult:deviation} is presented in Section
\ref{sec:proof:mainresult:deviation}.

\begin{Remark}[Interpretation of Theorem \ref{mainresult:deviation}]
The bound in Theorem \ref{mainresult:deviation} is expressed in terms of the threshold $\varepsilon$ and the following ratios:
\begin{align*}
\frac{\mu_{1,m}}{L_{m}}, \frac{\mu_{2,m}}{L_{m}}, \frac{\mu_{3,m}^{2}}{L_{m}},  \frac{\mu_{4,m}^{2}}{L_{m}},\frac{\sigma_{m,1}^{2}}{L^{2}_{m}}, \frac{\sigma_{m,2}^{2}}{L^{2}_{m}},
\frac{\sigma_{m,3}^{2}}{L_{m}}, \frac{\sigma_{m,4}^{2}}{L_{m}}.
\end{align*}
To highlight that these ratios become small when the noise strength is low relative to the clean signal, that is, when the signal-to-noise ratio (SNR) is high, we consider the single-component case $M=1$ for brevity.
For notational simplicity, below we write
$A$, $L$, $\mu_{1},\ldots,\mu_{4}$, and $\sigma_{1},\ldots,\sigma_{4}$
in place of
$A_{1}$, $L_{1}$, $\mu_{1,1},\ldots,\mu_{1,4}$, and $\sigma_{1,1},\ldots,\sigma_{1,4}$.

\begin{itemize}
\item $L$: A larger value of the second derivative of $|W_{f}(t,\cdot)|^{2}$ with respect to the scale variable
at $s_{f}(t)$
indicates that $|W_{f}(t,\cdot)|$ exhibits a sharper peak at $s_{f}(t)$.
This can be associated with a rapid change in the spectral content of the signal $f$.
Hence, when
$L$ is larger, the ridge point of $f$ is less likely to be disturbed.
Due to the slowly varying of the magnitude function $A$ and hence the approximation
\begin{align}\label{normalization_approx}
S_{f}(t,s) \approx A^{2}(t)S_{f/A}(t,s),
\end{align}
$L$ is approximately proportional to $A^{2}(t)$.
In this case, a larger
$A$ reduces the possibility of ridge point disturbance.

\item $\mu_{1}$, $\mu_{2}$, $\mu_{3}$, $\mu_{4}$: Under the approximation (\ref{normalization_approx}),
\begin{align*}
\mu_{1} \approx
A(t)\sqrt{\textup{Var}(\Phi)}\mathbb{E}\left[\underset{s\in  I}{\max}\left|\frac{\partial W_{f/A}}{\partial s}(t,s)W_{\widetilde{\Phi}}(t,s)\right|\right],
\end{align*}
where $\widetilde{\Phi}=\Phi/\sqrt{\textup{Var}(\Phi)}$.
Similarly,
$\mu_{2}$ is approximately proportional to
$A(t)\sqrt{\textup{Var}(\Phi)}$.
Both $\mu_{3}$ and $\mu_{4}$ are approximately proportional to $\sqrt{\textup{Var}(\Phi)}$.
This implies that
\begin{align}\label{interpretation:mu/L}
\frac{6}{L}\max\left\{\mu_{1},
\mu_{2},\mu_{3}^{2},\mu_{4}^{2}\right\}
\approx \max\left\{\mathcal{O}\left(\frac{\sqrt{\textup{Var}(\Phi)}}{A}\right),\mathcal{O}\left(\frac{\textup{Var}(\Phi)}{A^{2}}\right)\right\}.
\end{align}
Therefore, the available region
for $\varepsilon$ in (\ref{available_region}) becomes larger when the ratio $A(t)/\sqrt{\textup{Var}(\Phi)}$ increases.

\item For any fixed $\varepsilon$ satisfying (\ref{available_region}), (\ref{interpretation:mu/L}) shows that the exponents
$$
\frac{\varepsilon }{6}-\frac{\mu_{1}}{L},\ \frac{\varepsilon }{6}-\frac{\mu_{2}}{L},\ \sqrt{\frac{\varepsilon }{6}}-\frac{\mu_{3}}{\sqrt{L}},\ \sqrt{\frac{\varepsilon }{6}}-\frac{\mu_{4}}{\sqrt{L}}
$$
become larger when the ratio $A/\sqrt{\textup{Var}(\Phi)}$ increases.
Thus,  (\ref{main_4_exponent})  implies that the probability of $ s_{Y}(t)$ deviating from $s_{f}(t)$ by more than $\varepsilon$
decreases exponentially fast when the ratio $A/\sqrt{\textup{Var}(\Phi)}$ increases.

\item Exponents $L^{2}/\sigma_{1}^{2}$, $L^{2}/\sigma_{2}^{2}$, $L/\sigma_{3}^{2}$, and $L/\sigma_{4}^{2}$:
Since both $\sigma_{1}^{2}$ and
$\sigma_{2}^{2}$ are approximately proportional to $A^{2}\textup{Var}(\Phi)$,
and both $\sigma_{3}^{2}$ and
$\sigma_{4}^{2}$ are approximately proportional to $\textup{Var}(\Phi)$,
we have
\begin{align*}
L^{2}/\sigma_{k}^{2}=\mathcal{O}\left(\frac{A^{2}}{\textup{Var}(\Phi)}\right)\ \textup{for}\ k=1,2,
\end{align*}
and
\begin{align*}
L/\sigma_{k}^{2}=\mathcal{O}\left(\frac{A^{2}}{\textup{Var}(\Phi)}\right)\ \textup{for}\ k=3,4.
\end{align*}
This also implies that, for any $\varepsilon$ satisfying (\ref{available_region}),
the probability that $ s_{Y}(t)$ deviates from $s_{f}(t)$ by more than $\varepsilon$
decreases exponentially fast when the ratio $A/\sqrt{\textup{Var}(\Phi)}$ increases.

\end{itemize}

\end{Remark}

\begin{Remark}
The exponential bounds in Theorems~\ref{lemma:conditional} and~\ref{mainresult:deviation} are derived using a generalized Borell-TIS inequality (Lemma~\ref{lemma:complex-Borell}) for the modulus of Gaussian processes.
Since this generalized Borell-TIS inequality neither requires stationarity nor imposes specific structural assumptions on the underlying random process, it is applicable to the nonstationary random processes arising from the modulus of the AWT of an adaptive harmonic signal corrupted by Gaussian noise.
On the other hand, the Borell-TIS inequality does not account for the dependence structure induced by the AWT.
As a consequence of this generality, the constants in Theorems~\ref{lemma:conditional} and~\ref{mainresult:deviation}
are conservative.

\end{Remark}

\subsection{Numerical simulations}
The dependence of the bound \eqref{main_4_exponent} on SNR is intuitive. To illustrate this intuition, we present a numerical experiment.
Throughout all figures in this paper, we use a mother wavelet with center frequency $\omega_{\psi}= 80$ Hz, as defined in (\ref{def:centralfrequency_psi}).
Take a 60-second frequency-modulated signal $f$ sampled at $f_s:=100$ Hz and contaminated by a Gaussian noise $\Phi$ of different intensities. We generate $10^4$ realizations of the noise and carry out the following analysis.
For each realization of the noise, the SNR (in dB) of the noise-contaminated signal $Y=f+\Phi$ is defined as $10  \log_{10}(\|f\|^{2}_{\ell^2} / \|\Phi\|_{\ell^{2}}^{2})$.
The ridge deviation is quantified in two forms:
 \begin{align}\label{def:triangle_ridgedeviation}
 \triangle = \frac{1}{60f_s}\overset{60f_s}{\underset{k=1}{\sum}} |s_{f}(k/f_s)-s_{Y}(k/f_s)|
 \ (\textup{units}: \omega_{\psi}/\textup{Hz})
\end{align}
 and
 \begin{align*}
 \widetilde{\triangle} = \frac{1}{60f_s}\overset{60f_s}{\underset{k=1}{\sum}} |\omega_{\psi}s_{f}^{-1}(k/f_s)-\omega_{\psi}s_{Y}^{-1}(k/f_s)|\ (\textup{units: Hz}),
 \end{align*}
where $s_{f}$ and $s_{Y}$ denote the ridge curves of the clean and noisy signals, respectively. Note that $\triangle$ measures how accurate the ridge is when viewing $s_f$ as the ground truth, and $\widetilde{\triangle}$ expresses the deviation in the more familiar unit of Hertz.

\begin{figure}[hbt!]
\centering
\includegraphics[width=0.454\textwidth]{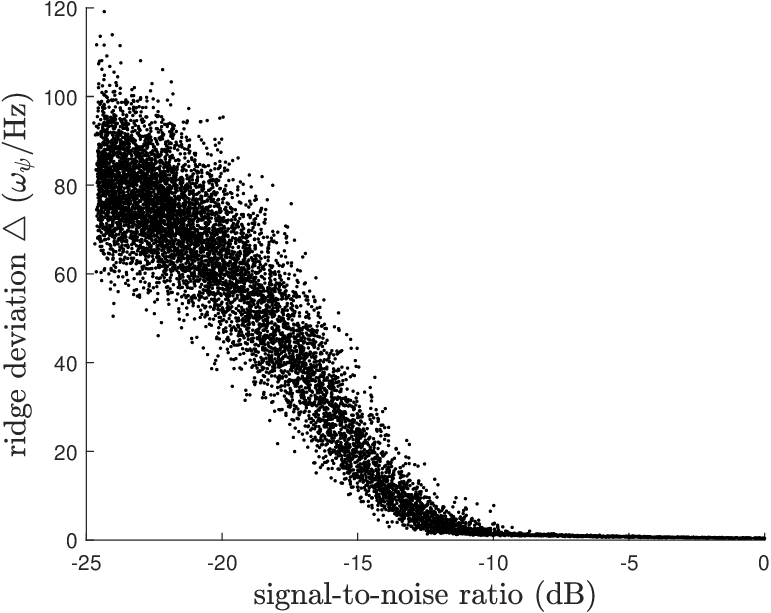}
\includegraphics[width=0.454\textwidth]{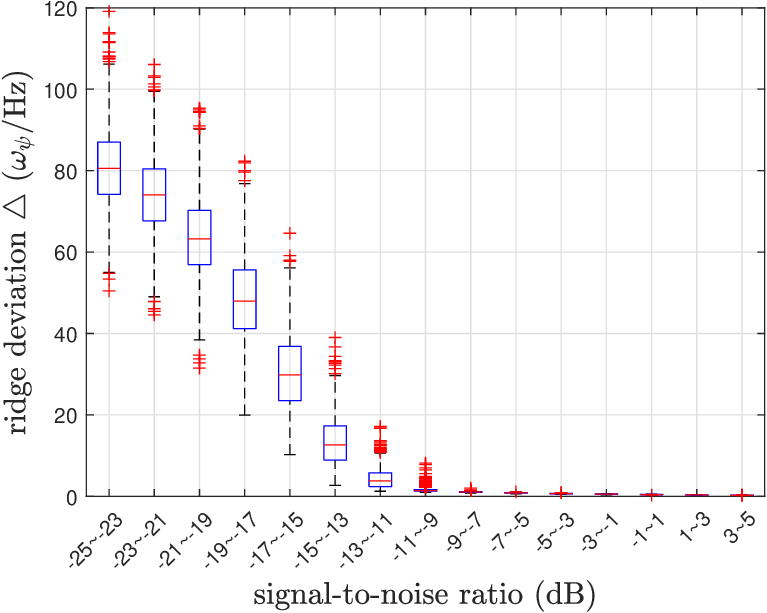}
\includegraphics[width=0.454\textwidth]{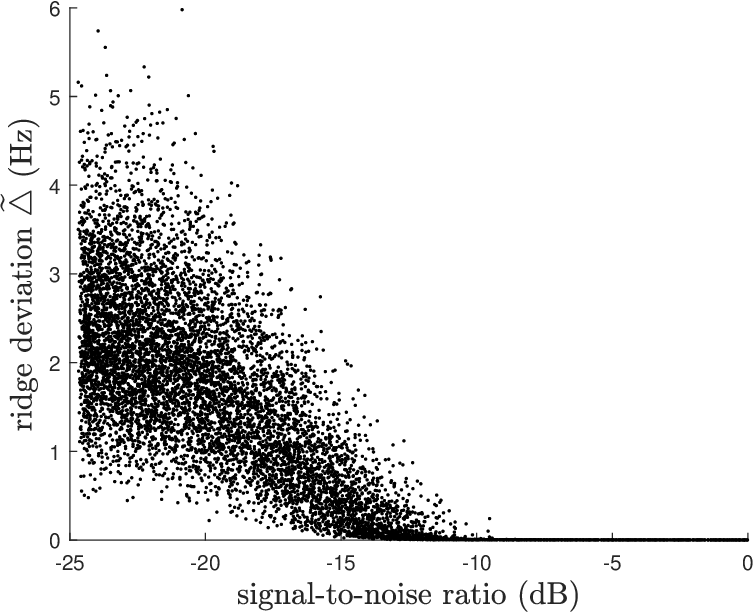}
\includegraphics[width=0.454\textwidth]{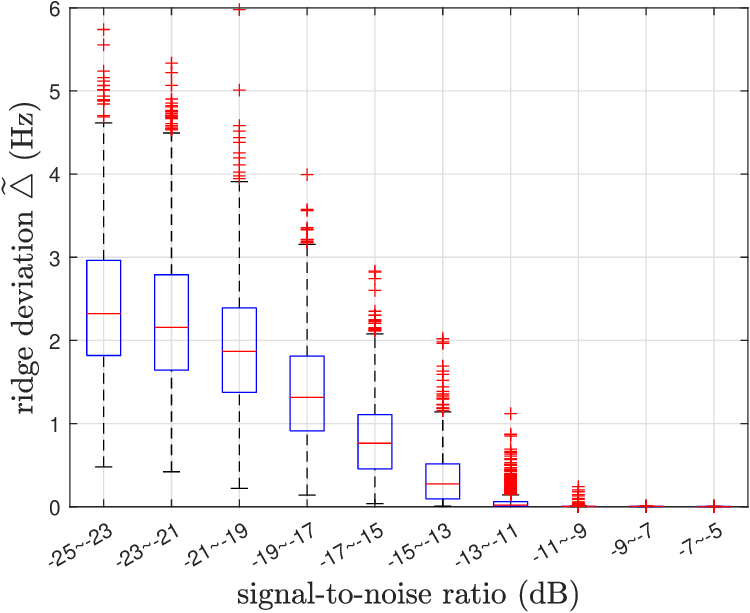}
 \caption{Effect of signal-to-noise ratio on wavelet ridge deviation.
  \label{fig:effect_SNR}}
\end{figure}

The result is shown in Figure \ref{fig:effect_SNR}. The left panel displays a scatter plot of all $10^4$ trials, showing the relationship between the SNR (in dB) and $\triangle$ or $\widetilde{\triangle}$. The right panel presents a box plot summarizing the ridge deviations within different SNR intervals.
The result shows that when the SNR is relatively high (resp. low), the ridge point $s_{Y,m}(t)$ of the noisy signal $Y$ is unlikely (resp. likely) to deviate too far from the ridge point $s_{f,m}(t)$ of the clean signal $f$.
This numerical result echoes results shown in Theorem \ref{mainresult:deviation}, providing a bound on the probability that the deviation $|s_{Y,m}(t)-s_{f,m}(t)|$ exceeds a given threshold.

\subsection{Dudley's theorem and the Borell-TIS inequality for AWT}\label{sec:Dudley}

To derive Theorems \ref{lemma:conditional} and \ref{mainresult:deviation},
a key challenge is estimating the tail probabilities of the following random variables arising from
the AWT of the Gaussian process $\Phi$:

\begin{itemize}
\item The sum of  $|W_{\Phi}(t, s_{f}(t))|$ and the supremum of $|W_{\Phi}(t, s)|$ over $s \in I^c$,
where $I$ is an interval containing $s_{f}(t)$ and $I^{c}$ is the complement of $I$, i.e.,
\begin{align}\label{type1_max}
|W_{\Phi}(t,s_{f}(t))| +\underset{s\in I^{c}}{\max}|W_{\Phi}(t,s)|.
\end{align}

\item The supremum of $\left|W_{\Phi}(t,s)\right|$ and $\left|\frac{\partial W_{\Phi}}{\partial s}(t,s)\right|$  over $s \in I$, i.e.,
\begin{align}\label{type2_max}
\underset{s\in I}{\max}\left|W_{\Phi}(t,s)\right|\ \textup{and}\ \underset{s\in I}{\max}\left|\frac{\partial W_{\Phi}}{\partial s}(t,s)\right|.
\end{align}

\item The maximum of the products
\begin{align}\label{type3_max}
\left|\frac{\partial W_{f}}{\partial s}(t,s)\right|\left|W_{\Phi}(t,s)\right|\ \textup{and}\ \left|W_{f}(t,s)\right|\left|\frac{\partial W_{\Phi}}{\partial s}(t,s)\right|
\end{align}
over $s \in I.$

\end{itemize}
The tail probabilities of the random variables above, with the replacement of the interval $I$ and its complement $I^{c}$
by a subinterval $I_{m}$ of $B_{m}(t)$ and the relative complement $B_{m}(t)\setminus I_{m}$, respectively,
also play an important role in the proof of Theorem \ref{mainresult:deviation}.

Our estimation method is inspired by the Borell-TIS inequality
\cite{talagrand2014upper,vershynin2018high},
which provides a tail probability estimate for the supremum of a real-valued non-stationary Gaussian process over an index set. In the proof of the Borell-TIS inequality,
the Gaussian concentration inequality for Lipschitz functions of real-valued Gaussian random vectors \cite{kuhn2023maximal} plays an important role.
Lemma \ref{lemma:Lip} below generalizes the concentration inequality in \cite{kuhn2023maximal}
to accommodate complex-valued and circularly symmetric Gaussian random vectors,
whose proof is provided in Section \ref{sec:proof:lemma:Lip}.

\begin{Lemma}\label{lemma:Lip}
Let $\mathbf{J}$ be a circularly symmetric Gaussian random vector in $\mathbb{C}^{d\times 1}$
with $\mathbb{E}\left[\mathbf{J}\right]=\mathbf{0}$, where $d\in \mathbb{N}$.
For any Lipschitz continuous function $g: \mathbb{C}^{d\times 1} \rightarrow \mathbb{R}$
and any $u>\mathbb{E}[g(\mathbf{J})]$,
\begin{align*}
\mathbb{P}\left(g(\mathbf{J})>u\right)\leq \textup{exp}\left(-\frac{(u-\mathbb{E}[g(\mathbf{J})])^{2}}{2\|g(\mathbf{A}\bullet)\|_{\textup{Lip}}^{2}}\right),
\end{align*}
where  $\mathbf{A}$ is a matrix in $\mathbb{C}^{d\times d}$ satisfying $\mathbf{A}\mathbf{A}^{*}=2^{-1}\mathbb{E}\left[\mathbf{J}\mathbf{J}^{*}\right]$.
Here, the symbol * represents the conjugate transpose.
The Lipschitz norm is defined by
\begin{align*}
\|g(\mathbf{A}\bullet)\|_{\textup{Lip}}= \underset{\begin{subarray}{c}\mathbf{x},\mathbf{y}\in \mathbb{C}^{d}\\ \mathbf{x}\neq\mathbf{y} \end{subarray}}
{\sup} \frac{|g(\mathbf{A}\mathbf{x})-g(\mathbf{A}\mathbf{y})|}{|\mathbf{x}-\mathbf{y}|}.
\end{align*}
\end{Lemma}

The Lipschitz continuous function $g$ in Lemma \ref{lemma:Lip} can be chosen as follows:
$$g(\mathbf{x}) = |x_{1}|+\underset{2\leq \ell \leq d}{\max}|x_{\ell}|$$
for case (\ref{type1_max}), and
$$g(\mathbf{x}) =\|\mathbf{x}\|_{\ell^{\infty}}:=\underset{1\leq \ell \leq d}{\max}|x_{\ell}|$$
for case (\ref{type2_max}), where  $\mathbf{x}=(x_{1},x_{2},\ldots,x_{d})\in \mathbb{C}^{d}$.

Note that in cases (\ref{type1_max}) and (\ref{type2_max}), the supremum is taken over the continuous index sets $I$ and $I^{c}$, respectively.
In order to apply the dominated convergence theorem to extend the result in Lemma \ref{lemma:Lip} to the continuous index set case, we need to ensure the convergence of the expected value of the supremum.
Details of this extension are provided in the proof of Theorem \ref{lemma:conditional}.
In \cite{dudley1967sizes,talagrand2014upper},
Dudley's theorem provides an upper bound
for the expected maximum of a centered real-valued Gaussian field over a compact set.
Because the AWT of the Gaussian process $\Phi$  is a complex-valued random process and the index set $I^{c}$ for case (\ref{type1_max}) is noncompact,
we present a complex-valued version of Dudley's theorem for the expectation of the supremum of the process $|W_{\Phi}(t,\cdot)|$ over the noncompact set $[0,\infty)$  in Lemma \ref{lemma:Dudley}.
The statement and proof of Lemma \ref{lemma:Dudley} rely on estimates of $\mathbb{E}[|W_{\Phi}(0,s)|^{2}]$ for large scale $s$ and the distance function $d_{W_{\Phi}}$ as follows.

\begin{Lemma}\label{lemma:d_s1_s2}
Suppose that $\psi$ is an analytic wavelet belonging to $L^{1}\cap L^{2}$, and its Fourier transform $\widehat{\psi}$
is Lipschitz continuous on $[0,\infty)$.
For the random process $\Phi$, if Assumption \ref{assumption:Gaussian} and Assumptions \ref{assumption:spectral_density}(a) and \ref{assumption:spectral_density}(b1) hold,
then for any constant $H^{-}\in(0,H)$, where $H$ is the long-memory parameter of $\Phi$ defined in (\ref{spectral:slowly}),  there exists a threshold $T>0$
such that for any $s\geq T$,
\begin{align}\label{decay_var_X1}
\mathbb{E}[|W_{\Phi}(0,s)|^{2}]\leq s^{-H^{-}}.
\end{align}
Under the same assumptions,
for any $s_{1},s_{2}>0$ satisfying $|s_{1}-s_{2}|\leq 1$, the distance $d_{W_{\Phi}}(s_{1},s_{2})$, defined in (\ref{def:canonical_distance}),
is $\gamma/2$-H$\ddot{o}$lder continuous:
\begin{align}\label{d_s1_s2_in_lemma}
d_{W_{\Phi}}(s_{1},s_{2})\leq C_{2}|s_{1}-s_{2}|^{\gamma/2},
\end{align}
where $\gamma$ is given in Assumption \ref{assumption:spectral_density}(a),
\begin{align}\label{def:C2}
C_{2} = 1 \vee \left[C_{1}\left(\frac{1}{\gamma}+\frac{2}{2-\gamma}\right)(2\|\widehat{\psi}\|_{\infty})^{2-\gamma}\|\widehat{\psi}\|_{\textup{Lip}}^{\gamma}
+2\|\widehat{\psi}\|_{\infty}^{2}\textup{Var}\left(\Phi(0)\right)\right]^{1/2},
\end{align}
and $\|\widehat{\psi}\|_{\textup{Lip}}$ represents the Lipschitz constant of $\widehat{\psi}$.

If the spectral density function of $\Phi$ satisfies Assumption \ref{assumption:spectral_density}(b2) instead of Assumption \ref{assumption:spectral_density}(b1), the results in (\ref{decay_var_X1}) and (\ref{d_s1_s2_in_lemma})
still hold with the substitution $H^{-}=1$.
\end{Lemma}

\begin{Lemma}[Dudley's theorem for AWT]\label{lemma:Dudley}
Suppose that $\psi$ is an analytic wavelet belonging to $L^{1}\cap L^{2}$, and its Fourier transform $\widehat{\psi}$
is Lipschitz continuous on $[0,\infty)$.
For the random process $\Phi$, if Assumptions \ref{assumption:Gaussian} and \ref{assumption:spectral_density} hold,
then
\begin{align}\label{lemma:ineq:supX}
\mathbb{E}\left[\underset{s> 0}{\sup}\left|W_{\Phi}(0,s)\right|\right]
\leq 4\left[\sqrt{\ln\mathcal{D}}\frac{\mathcal{D}^{3}}{(\mathcal{D}-1)^{2}}\sqrt{\frac{1}{\gamma}+\frac{1}{H^{-}}}
+\left(\sqrt{\frac{1}{\gamma}\ln C_{2}}+\sqrt{\ln T}\right)\frac{\mathcal{D}^{2}}{\mathcal{D}-1}\right],
\end{align}
where $H^{-}$, $C_{2}$, and $T$ are constants mentioned in Lemma \ref{lemma:d_s1_s2} and
\begin{align}\label{def:D}
\mathcal{D} = 2\vee \left\{\underset{s> 0}{\sup}\ \mathbb{E}\left[\left|W_{\Phi}(0,s)\right|^{2}\right]\right\}^{\frac{1}{2}}.
\end{align}

\end{Lemma}

In \cite{dudley1967sizes,talagrand2014upper}, the expectation
of the maximum of a centered real-valued Gaussian field over a compact set is bounded above by the Dudley entropy integral.
Lemma \ref{lemma:Dudley} extends this bound, providing an exact upper bound
for the expectation of the maximum of the modulus of a circularly symmetric and nonstationary complex Gaussian process over $(0,\infty)$.
The upper bound given in (\ref{lemma:ineq:supX})
increases as the long-memory parameter $H$ and the decay rate parameter $\gamma$ of the spectral density function of $\Phi$ decrease to zero.
This relationship echoes Slepian's inequality and
the Sudakov-Fernique inequality \cite[Chapter 7]{vershynin2018high}, which state that, for two real-valued Gaussian processes with the same mean, the process with higher correlation exhibits a greater expected maximum.
The proof of Lemma \ref{lemma:Dudley} is provided in Section \ref{sec:proof:lemma:Dudley}.
Finally, to estimate the tail probabilities of the random variables in (\ref{type3_max}), we introduce the complex-valued version of the Borell-TIS inequality as follows.

\begin{Lemma}[Borell-TIS inequality for AWT]\label{lemma:complex-Borell}
Suppose that the stationary Gaussian process $\Phi$ satisfies Assumptions \ref{assumption:Gaussian} and \ref{assumption:spectral_density}, and the mother wavelet $\psi$ satisfies Assumption \ref{assumption:boundedness:psi}.
For any $a,b \in  (0,\infty)$ and any bounded and continuous function $w:[a,b]\rightarrow [0,\infty)$,
define either
\begin{align*}
M_{[a,b]} := \underset{s\in [a,b]}{\max}\ w(s)\left|W_{\Phi}(0,s)\right|,
\ \ \
\sigma_{[a,b]} :=
\underset{s\in [a,b]}{\max}\ w(s)\left(\int_{0}^{\infty}|\widehat{\psi}\left(s\lambda\right)|^{2}p(\lambda)d\lambda\right)^{\frac{1}{2}},
\end{align*}
or
\begin{align*}
M_{[a,b]} := \underset{s\in [a,b]}{\max}\ w(s)\left|\frac{\partial W_{\Phi}}{\partial s}(0,s)\right|,
\ \ \
\sigma_{[a,b]} :=
\underset{s\in [a,b]}{\max}\ w(s)\left(\int_{0}^{\infty}|\lambda [D\widehat{\psi}]\left(s\lambda\right)|^{2}p(\lambda)d\lambda\right)^{\frac{1}{2}}.
\end{align*}
In both cases, we have $\mathbb{E}\left[M_{[a,b]}\right]<\infty$, and
for any $u\geq \mathbb{E}\left[M_{[a,b]}\right]$,
\begin{align*}
\mathbb{P}\left(M_{[a,b]}\geq u\right)
\leq
\textup{exp}\left(-\frac{(u-\mathbb{E}\left[M_{[a,b]}\right])^{2}}{\sigma_{[a,b]}^{2}}\right).
\end{align*}
\end{Lemma}
The proof of Lemma \ref{lemma:complex-Borell} is provided in Section \ref{sec:proof:lemma:complex-Borell}.
Last but not least, our work bridges and extends classical maximum theorems and maximum inequalities to address a key question in signal processing: how noise perturbs the ridge curves of adaptive harmonic signals that encode the IF information.

\section{Conclusion}\label{sec:conclusion}

In this work, we focus on providing a systematic study of Steps \texttt{I} and \texttt{II} as outlined in the Introduction. In addition to proposing a mathematical framework and definition for ridge analysis, we establish key properties of the ridge when viewed as a random process. Although we do not address Steps \texttt{III} and \texttt{IV} in detail, we briefly mention that most existing methods aim to fit a curve that captures the maximal energy on the TFR, using various weighting and fitting strategies \cite{carmona1997characterization,carmona1999multiridge,Ozkurt_Savaci_2005,iatsenko2016extraction,meignen2017demodulation,Zhu_Zhang_Gao_Li_2019,Legros_Fourer_2021}.
When multiple oscillatory components are present, the IF separation condition has been employed to support algorithm design \cite{Laurent_Meignen_2021}. In cases involving non-sinusoidal oscillations, geometric structures in the TFR can be leveraged to improve extraction accuracy \cite{su2024ridge}, and chirp characteristics can also be incorporated to enhance performance \cite{Colominas_Meignen_Pham_2020}. Moreover, specialized algorithms have been proposed to handle scenarios with significant spectral interference \cite{meignen2023new}.

\begin{figure}[hbpt!]
\centering
\subfigure[][]{\includegraphics[width=1\textwidth]{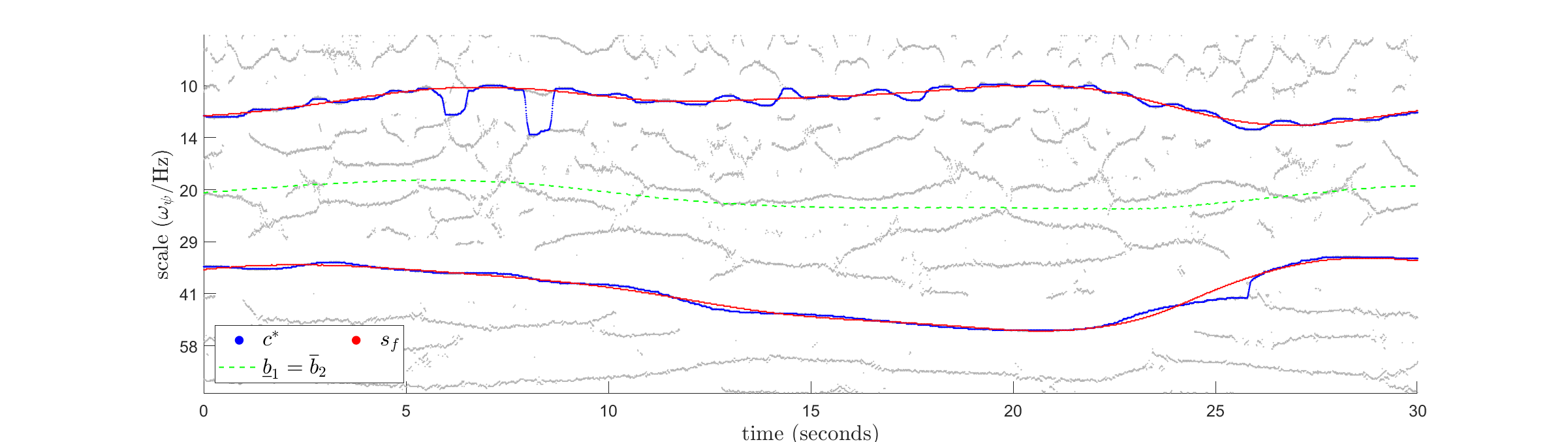}}
\subfigure[][]{
\includegraphics[width=0.454\textwidth]{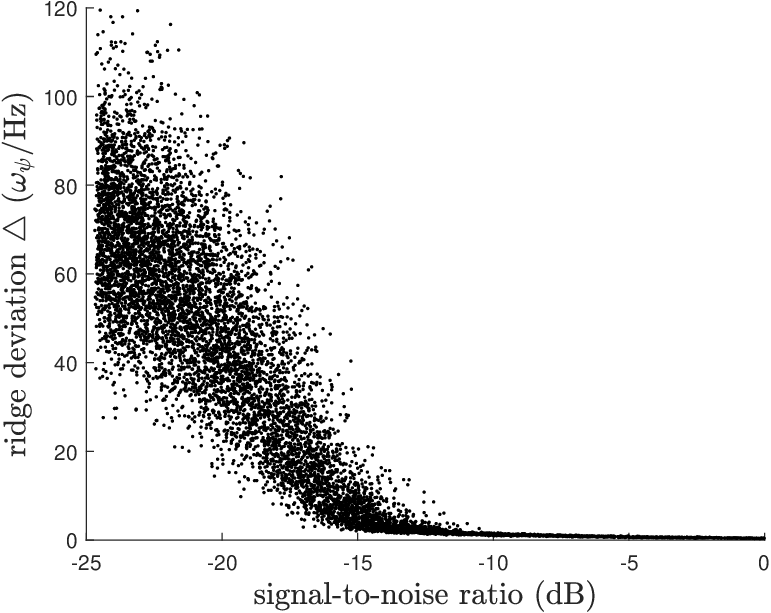}
\includegraphics[width=0.454\textwidth]{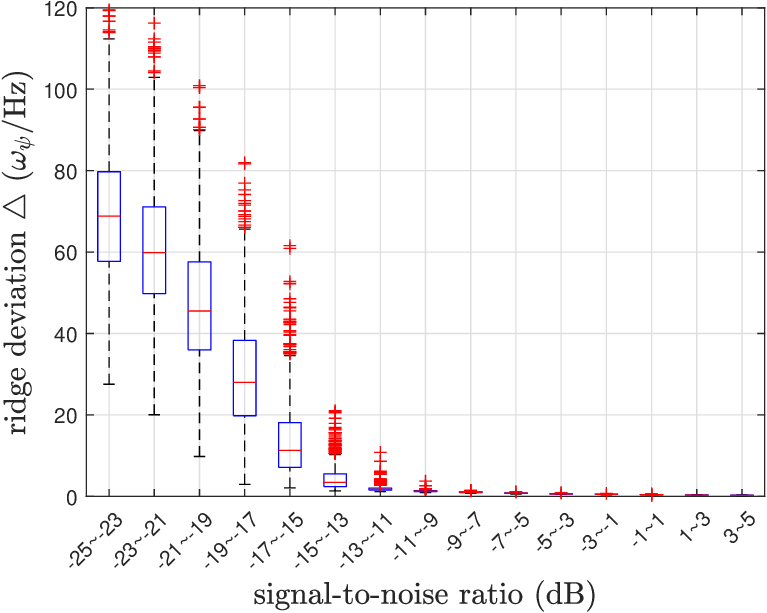}}
\caption{The top panel displays the ridge curve $c^{*}$, obtained by applying the ridge extraction algorithm
\eqref{eq ridge ex} to the scalogram of the noisy signal $Y$. The input signals $f$, $\Phi$, and $Y$ are the same as those in Figure \ref{fig:wavelet_potential_ridge}.
The bottom panel shows scatter and box plots of the deviation between $c^{*}$  and the ridge curve $s_{f}$ of the clean signal $f$. Here, we generate $10^4$ realizations of Gaussian noise, identical to those used in Figure \ref{fig:effect_SNR}.
For each realization, we compute the deviation value $\triangle$,
which is defined in (\ref{def:triangle_ridgedeviation}) with $s_{Y}$ replaced by $c^{*}$.
A Wilcoxon signed-rank test is performed on the $10^4$ paired deviation values $\triangle$, comparing the deviation of $s_{Y}$ from $s_{f}$ with that of $c^{*}$ from $s_{f}$.
The null hypothesis is that the median of the deviations $\triangle$ computed from $c^{*}$ is greater than or equal to that computed from $s_{Y}$.
The resulting $p$-value is less than $10^{-3}$,  indicating that the ridge $c^{*}$ is closer to $s_{f}$ than $s_{Y}$.
}
\label{fig:wavelet_potential_ridge2}
\end{figure}

We shall discuss the relationship among the established theory, the numerical simulation shown in Figure~\ref{fig:wavelet_potential_ridge}, and the results obtained by practitioners using a commonly applied ridge extraction algorithm. We also illustrate how this algorithm connects to the theoretical insights developed in this paper.
For the scalogram $S_Y$ of the noisy signal $Y$, denote its discretization as $R\in \mathbb{R}^{p\times n}$, where $p$ is the number of discretized scales and $n$ is the number of time samples. The ridge extraction algorithm fits a curve to the scalogram by solving the following optimization problem:
\begin{equation}\label{eq ridge ex}
c^*:=\arg\max_{c\in \{1,\ldots,p\}^n}\sum_{j=1}^n |\tilde R(c(j),j)|-\lambda\sum_{j=1}^{n-1}|c(j+1)-c(j)|^2\,,
\end{equation}
where $\tilde{R}\in \mathbb{R}^{p\times n}$, $$\tilde{R}(i,j)=\log\left(\frac{R(i,j)}{\sum_{i'=1}^p\sum_{j'=1}^n|R(i',j')|}\right),$$
and $\lambda\geq0$ is chosen by the user to control the smoothness of the extracted ridge.
The curve $c^{*}$, obtained by solving (\ref{eq ridge ex}) with $\lambda=0.1$, is illustrated in Figure \ref{fig:wavelet_potential_ridge2}(a) for the case where the clean signal $f$ consists of two frequency-modulated sinusoids, and in Figure \ref{fig:wavelet_potential_ridge:one}(b) for the case where
$f$ consists of a single frequency-modulated sinusoid.

Since $R$ is a discretized random field, technically $c^*$ is a random sequence as well. When $\lambda=0$, $c^*$ is the numerical realization of $s_Y$ defined in \eqref{proxy:s_Y}.
When $\lambda>0$, a smoothing effect is clearly observed in Figures \ref{fig:wavelet_potential_ridge2}(a) and \ref{fig:wavelet_potential_ridge:one}(b).
The penalty term in (\ref{eq ridge ex}) trades off local maxima for regularity, thereby encouraging the desired continuity of the ridge.
As a result, the ridge $c^*$ differs from the discretized ridge $s_Y$.
A comparison of the box plots in Figures \ref{fig:effect_SNR} and \ref{fig:wavelet_potential_ridge2}
further shows that
$c^{*}$ is closer to the discretized ridge of the clean signal than $s_{Y}$ is.
The parameter $\lambda$ can thus be interpreted as controlling the strength of ``denoising'' applied to $s_Y$.
Analyzing the statistical behavior of $c^*$ thus provides insight into the ridge extraction algorithm, which is the goal of Step \texttt{IV}.

It is worth noting that the ridge extraction algorithm \eqref{eq ridge ex}, like most existing methods, relies on an implicit assumption of ridge regularity that was adopted before the ridge's behavior was well understood. Our analysis suggests that this implicit assumption about ridge should be updated, potentially leading to improvements in the algorithm.
In particular, a natural and interesting question arises: how often does a ridge exhibit discontinuities or ``jumps''?
Numerically, this appears to depend on the SNR, with lower SNRs tending to produce more frequent jumps.
Quantifying this relationship involves studying the excursion properties of the ridge as a random process.
This is the joint objective of Steps \texttt{III} and \texttt{IV}.

Finally, we comment on the scope of applicability and limitations of the proposed framework.
The second half of this work focuses on the deviation of ridge points of a noisy signal from those of the corresponding noise-free signal.
In other words, our analysis concerns the difference between the local maximizers of the scalograms of the noise-affected and noise-free signals, rather than the difference between the ridge points of the noisy signal and the IFs of the noise-free signal.
As a result, the findings in this work can be extended to more general signal models, provided that the scalogram of the noise-free signal $f$ exhibits well-separated local maxima at each time. We shall mention that our current result is not sufficient to fully describe the general model like the adaptive non-harmonic model \cite{wu2013instantaneous}, defined as
\begin{align*}
f_w(t) = \sum_{m=1}^{M} A_m(t) s_m\left( \phi_m(t) \right),
\end{align*}
where $M$, $A_{m}(t)$, $\phi_m(t)$, and $\phi_{m}'(t)$ have the same interpretation as those of AHM, and
$s_{m}$ is the wave-shape function of the $m$-th component, which is
a $1$-periodic smooth function describing the non-sinusoidal oscillatory shape. Due to the non-sinusoidal oscillatory nature of $s_m$, $f_w$ can be written as
\begin{align*}
f_w(t) = \sum_{m=1}^{M} \sum_{k=1}^{D_m} A_m(t) a_{m,k}\cos( 2\pi k\phi_m(t) +\beta_{m,k}),
\end{align*}
where $D_m\geq 1$, $\beta_{m,k}\in [0,2\pi)$, and $\alpha_{m,k}\geq 0$ are the Fourier coefficients of $s_m$. With this expression, $f_w$ can be viewed as a generalized AHM with many IMT functions, while different IMT functions might have close or even crossover IFs. Since our results rely on the well-separated IFs of all IMT components, the analysis applies only to the special situation where $D_m$ is sufficiently small and the multiples of distinct $\phi_m'(t)$ are well separated. Extending the theory to settings with spectral interference; that is, when IFs are close or intersect, remains an open problem for future ridge analysis.

Regarding noise models, our deviation bounds rely on Gaussian concentration inequalities for complex-valued Gaussian random variables arising from the analytic wavelet transform. These tools inherently require a Gaussian assumption. Extending the analysis to non-Gaussian noise would necessitate new maximal inequalities for complex-valued transforms of non-Gaussian random processes. For non-stationary noise, the covariance structure of the analytic wavelet transform becomes more involved, and corresponding theoretical extensions require further technical development. These issues are directly relevant in practice, where signals are discretized, and noise is typically non-Gaussian. Under mild non-stationarity assumptions (including controlled dependence and moments) and uniform sampling, it has been shown in \cite{liu2025random} that the distribution of the discretized AWT of non-Gaussian, causal, non-stationary time series can be approximated by that of a suitable Gaussian sequence, so our results remain applicable. However, the case of non-uniform sampling remains open.
Moreover, in real-world applications, oscillatory patterns may change abruptly, making the assumed adaptive harmonic model unsuitable. In such situations, ridges may become discontinuous, and ridge detection must incorporate change-point structure. Extending the current framework to address these challenges is both difficult and important, and we leave these directions for future work.

\section{Proofs}\label{sec:proof}

\subsection{Proof of Proposition \ref{lemma:exp}}\label{sec:proof:lemma:exp}

(a) To prove that $\mathbf{W}$ is circularly symmetric, it suffices to show that any element in the
pseudo-covariance matrix of
$\mathbf{W}$ is zero.
Denote $W_{\Phi}^{R}=\textup{Re}\left(W_{\Phi}\right)$ and $W_{\Phi}^{I}=\textup{Im}\left(W_{\Phi}\right)$.
For any $k,k'\in \{1,\ldots,d\}$,
\begin{align}\notag
\mathbb{E}\left[W_{\Phi}(t,s_{k})W_{\Phi}(t,s_{k'})\right]
=& \mathbb{E}\left[W_{\Phi}^{R}(t,s_{k})W_{\Phi}^{R}(t,s_{k'})\right]
- \mathbb{E}\left[W_{\Phi}^{I}(t,s_{k})W_{\Phi}^{I}(t,s_{k'})\right]
\\\label{proof:circular}&+i \mathbb{E}\left[W_{\Phi}^{R}(t,s_{k})W_{\Phi}^{I}(t,s_{k'})\right]
+i \mathbb{E}\left[W_{\Phi}^{I}(t,s_{k})W_{\Phi}^{R}(t,s_{k'})\right].
\end{align}
Denote $\psi_{R}=\textup{Re}(\psi)$  and $\psi_{I}=\textup{Im}(\psi)$.
Analogous to the representation in (\ref{spect:WPhi}),
\begin{align}\label{spectral:realWPhi}
W^{R}_{\Phi}(t,s) =  \int_{\mathbb{R}}
 e^{it\lambda}  \overline{\widehat{\psi_{R}}(s\lambda)} \sqrt{p(\lambda)}Z(d\lambda)
\end{align}
and
\begin{align}\label{spectral:imagWPhi}
W^{I}_{\Phi}(t,s) =  \int_{\mathbb{R}}
e^{it\lambda}  \overline{\widehat{\psi_{I}}(s\lambda)} \sqrt{p(\lambda)}Z(d\lambda).
\end{align}
Because  $\psi$ is an analytic wavelet,
the Fourier transforms of $\psi_{R}$ and $\psi_{I}$ are related by the following equation:
\begin{align}\label{Fourier_relation_psihat}
\widehat{\psi_{I}}(\zeta) = -i\textup{sgn}(\zeta)\widehat{\psi_{R}}(\zeta),\ \zeta\in \mathbb{R},
\end{align}
where $\textup{sgn}$ is the signum function.
By the orthogonal property (\ref{ortho}) of the Gaussian random measure $Z$ and (\ref{Fourier_relation_psihat}),
\begin{align}\notag
\mathbb{E}\left[W_{\Phi}^{I}(t,s_{k})W_{\Phi}^{I}(t,s_{k'})\right]
=&\int_{-\infty}^{\infty}\overline{\widehat{\psi_{I}}(s_{k}\lambda)}
\widehat{\psi_{I}}(s_{k'}\lambda)p(\lambda)d\lambda
\\\notag=&\int_{-\infty}^{\infty}\overline{\widehat{\psi_{R}}(s_{k}\lambda)}
\widehat{\psi_{R}}(s_{k'}\lambda)p(\lambda)d\lambda
\\\label{proof:circular_1}=&\mathbb{E}\left[W_{\Phi}^{R}(t,s_{k})W_{\Phi}^{R}(t,s_{k'})\right]
\end{align}
and
\begin{align}\notag
\mathbb{E}\left[W_{\Phi}^{I}(t,s_{k})W_{\Phi}^{R}(t,s_{k'})\right]
=&\int_{-\infty}^{\infty}\overline{\widehat{\psi_{I}}(s_{k}\lambda)}\widehat{\psi_{R}}(s_{k'}\lambda)p(\lambda)
d\lambda
\\\notag=&-\int_{-\infty}^{\infty}\overline{\widehat{\psi_{R}}(s_{k}\lambda)}\widehat{\psi_{I}}(s_{k'}\lambda)p(\lambda)
p(\lambda)d\lambda
\\\label{proof:circular_2}=&-\mathbb{E}\left[W_{\Phi}^{R}(t,s_{k})W_{\Phi}^{I}(t,s_{k'})\right].
\end{align}
By substituting (\ref{proof:circular_1}) and (\ref{proof:circular_2}) into (\ref{proof:circular}),
we obtain  $\mathbb{E}\left[W_{\Phi}(t,s_{k})W_{\Phi}(t,s_{k'})\right] =0$.

(b) For any $t\in \mathbb{R}$, $s_{1},s_{2}>0$ with $s_{1}\neq s_{2}$,
by (\ref{spectral:realWPhi}) and (\ref{spectral:imagWPhi}),
\begin{align*}
W_{\Phi}^{R}(t,s_{1})-W_{\Phi}^{R}(t,s_{2})
= \int_{\mathbb{R}}e^{it\lambda}\left[\ \overline{\widehat{\psi_{R}}(s_{1}\lambda)}-\overline{\widehat{\psi_{R}}(s_{2}\lambda)}\ \right] \sqrt{p(\lambda)}Z(d\lambda)
\end{align*}
and
\begin{align*}
W_{\Phi}^{I}(t,s_{1})-W_{\Phi}^{I}(t,s_{2})
= \int_{\mathbb{R}}e^{it\lambda}\left[\ \overline{\widehat{\psi_{I}}(s_{1}\lambda)}-\overline{\widehat{\psi_{I}}(s_{2}\lambda)}\ \right] \sqrt{p(\lambda)}Z(d\lambda).
\end{align*}
By the orthogonal property (\ref{ortho}) of the Gaussian random measure $Z$,
\begin{align}\notag
&\mathbb{E}\left[\left(W_{\Phi}^{R}(t,s_{1})-W_{\Phi}^{R}(t,s_{2})\right)
\left(W_{\Phi}^{I}(t,s_{1})-W_{\Phi}^{I}(t,s_{2})\right)\right]
\\\label{increment_real_imag}=&
\int_{\mathbb{R}}\left[\widehat{\psi_{R}}(s_{1}\lambda)-\widehat{\psi_{R}}(s_{2}\lambda)\right]
\left[\ \overline{\widehat{\psi_{I}}(s_{1}\lambda)}-\overline{\widehat{\psi_{I}}(s_{2}\lambda)}\ \right]
p(\lambda)d\lambda.
\end{align}
Applying (\ref{Fourier_relation_psihat}) again,  and noting that the spectral density function
 $p$ is even, we can rewrite (\ref{increment_real_imag}) as follows
\begin{align}\notag
&\mathbb{E}\left[\left(W_{\Phi}^{R}(t,s_{1})-W_{\Phi}^{R}(t,s_{2})\right)
\left(W_{\Phi}^{I}(t,s_{1})-W_{\Phi}^{I}(t,s_{2})\right)\right]
\\\notag=&
i\int_{\mathbb{R}}\left|\widehat{\psi_{R}}(s_{1}\lambda)-\widehat{\psi_{R}}(s_{2}\lambda)\right|^{2}
\textup{sgn}(\lambda)p(\lambda)d\lambda=0,
\end{align}
which shows that
the real and imaginary parts of $W_{\Phi}(t,s_{1})-W_{\Phi}(t,s_{2})$ are independent.
Following the same procedure, we can obtain that
\begin{align}\notag
&\mathbb{E}\left[\left|W_{\Phi}^{R}(t,s_{1})-W_{\Phi}^{R}(t,s_{2})\right|^{2}\right]
\\\label{Ereal2=Eimag2}=&
\mathbb{E}\left[\left|W_{\Phi}^{I}(t,s_{1})-W_{\Phi}^{I}(t,s_{2})\right|^{2}\right]
=\frac{1}{2}d^{2}_{W_{\Phi}}(s_{1},s_{2}).
\end{align}
The equality (\ref{Ereal2=Eimag2}) implies that
\begin{align}\notag
&\frac{|W_{\Phi}(t,s_{1})-W_{\Phi}(t,s_{2})|^{2}}{d^{2}_{W_{\Phi}}(s_{1},s_{2})}
\\\notag=&\frac{1}{2}\left(\frac{|W_{\Phi}^{R}(t,s_{1})-W_{\Phi}^{R}(t,s_{2})|^{2}}
{\mathbb{E}\left[|W_{\Phi}^{R}(t,s_{1})-W_{\Phi}^{R}(t,s_{2})|^{2}\right]}
+\frac{|W_{\Phi}^{I}(t,s_{1})-W_{\Phi}^{I}(t,s_{2})|^{2}}
{\mathbb{E}\left[|W_{\Phi}^{I}(t,s_{1})-W_{\Phi}^{I}(t,s_{2})|^{2}\right]}\right).
\end{align}
Proposition \ref{lemma:exp}(b) follows from the fact that the distribution of half the sum of two independent chi-square random variables, each with one degree of freedom, is an exponential distribution with a mean of one.

(c) By (\ref{spectral:realWPhi}), (\ref{spectral:imagWPhi}), and It$\hat{\textup{o}}$'s formula,
\begin{align*}
\left(W^{R}_{\Phi}(t,s)\right)^{2} = \mathbb{E}\left[\left(W^{R}_{\Phi}(t,s)\right)^{2}\right]
+ \int^{'}_{\mathbb{R}^{2}}
e^{it(\lambda_{1}+\lambda_{2})}  \overline{\widehat{\psi_{R}}(s\lambda_{1})\widehat{\psi_{R}}(s\lambda_{2})} \sqrt{p(\lambda_{1})p(\lambda_{2})}Z(d\lambda_{1})Z(d\lambda_{2}),
\end{align*}
where $\int^{'}$ means that the integral excludes the diagonal hyperplanes
$\lambda_{1}=\mp \lambda_{2}$ \cite{major1981lecture}.
Similarly,
\begin{align*}
\left(W^{I}_{\Phi}(t,s)\right)^{2} = \mathbb{E}\left[\left(W^{I}_{\Phi}(t,s)\right)^{2}\right]
+ \int^{'}_{\mathbb{R}^{2}}
e^{it(\lambda_{1}+\lambda_{2})}  \overline{\widehat{\psi_{I}}(s\lambda_{1})\widehat{\psi_{I}}(s\lambda_{2})} \sqrt{p(\lambda_{1})p(\lambda_{2})}Z(d\lambda_{1})Z(d\lambda_{2}).
\end{align*}
Hence, for any $s>0$,
\begin{align}\notag
&S_{\Phi}(t,s)-\mathbb{E}\left[S_{\Phi}(t,s)\right]
\\\notag=&\int^{'}_{\mathbb{R}^{2}}
e^{it(\lambda_{1}+\lambda_{2})} \left[
 \ \overline{\widehat{\psi_{R}}(s\lambda_{1})\widehat{\psi_{R}}(s\lambda_{2})}
 +
 \overline{\widehat{\psi_{I}}(s\lambda_{1})\widehat{\psi_{I}}(s\lambda_{2})}
\  \right]\sqrt{p(\lambda_{1})p(\lambda_{2})}Z(d\lambda_{1})Z(d\lambda_{2}).
\end{align}
By the product formula for the second Wiener chaos \cite{major1981lecture,nourdin2012normal},
\begin{align}\notag
\textup{Cov}\left(S_{\Phi}(t,s_{1}),S_{\Phi}(t,s_{2})\right)
=&2\int_{\mathbb{R}^{2}}
\left[
 \ \overline{\widehat{\psi_{R}}(s_{1}\lambda_{1})\widehat{\psi_{R}}(s_{1}\lambda_{2})}
+
 \overline{\widehat{\psi_{I}}(s_{1}\lambda_{1})\widehat{\psi_{I}}(s_{1}\lambda_{2})}
\  \right]
\\\label{cov:sumvarphi}&\times\left[\widehat{\psi_{R}}(s_{2}\lambda_{1})\widehat{\psi_{R}}(s_{2}\lambda_{2})
 +
\widehat{\psi_{I}}(s_{2}\lambda_{1})\widehat{\psi_{I}}(s_{2}\lambda_{2})
\right]p(\lambda_{1})p(\lambda_{2})d\lambda_{1}d\lambda_{2}
\\\notag =& \overset{2}{\underset{i=1}{\sum}}\overset{2}{\underset{j=1}{\sum}}\varphi_{i,j}(s_{1},s_{2}),
\end{align}
where
\begin{align}\label{def:varph_iij}
\varphi_{i,j}(s_{1},s_{2})=2\left(\int_{\mathbb{R}}\overline{\widehat{\psi_{\alpha(i)}}(s_{1}\lambda)}\widehat{\psi_{\alpha(j)}}(s_{2}\lambda)
p(\lambda)d\lambda\right)^{2},
\end{align}
$\alpha(1)=R$, and $\alpha(2)=I$.
Because $\psi$ is an analytic wavelet, the relation (\ref{Fourier_relation_psihat}) implies that $\varphi_{1,1}(s_{1},s_{2})=\varphi_{2,2}(s_{1},s_{2})$ and $\varphi_{1,2}(s_{1},s_{2})=\varphi_{2,1}(s_{1},s_{2})$.

On the other hand, from (\ref{spectral:realWPhi}) and (\ref{spectral:imagWPhi}), we have
\begin{align}\label{cov:WRWR}
\mathbb{E}\left[W^{R}_{\Phi}(t,s_{1})W^{R}_{\Phi}(t,s_{2})\right]=
\int_{\mathbb{R}}\overline{\widehat{\psi_{R}}(s_{1}\lambda)}\widehat{\psi_{R}}(s_{2}\lambda)
p(\lambda)d\lambda
\end{align}
and
\begin{align}\label{cov:WRWI}
\mathbb{E}\left[W^{R}_{\Phi}(t,s_{1})W^{I}_{\Phi}(t,s_{2})\right]=
\int_{\mathbb{R}}\overline{\widehat{\psi_{R}}(s_{1}\lambda)}\widehat{\psi_{I}}(s_{2}\lambda)
p(\lambda)d\lambda.
\end{align}
By (\ref{def:varph_iij}), (\ref{cov:WRWI}), and (\ref{cov:WRWR}),
\begin{equation}\label{relation:varphi_ij}
\varphi_{i,j}(s_{1},s_{2}) =\left\{\begin{array}{ll}
2\left(\mathbb{E}\left[W^{R}_{\Phi}(t,s_{1})W^{R}_{\Phi}(t,s_{2})\right]\right)^{2}\ &\ \textup{if}\ i=j,
\\
2\left(\mathbb{E}\left[W^{R}_{\Phi}(t,s_{1})W^{I}_{\Phi}(t,s_{2})\right]\right)^{2}\ &\ \textup{if}\ i\neq j.
\end{array}\right.
\end{equation}
The proof of the first equality in Proposition \ref{lemma:exp}(c) is concluded by substituting (\ref{relation:varphi_ij}) into (\ref{def:varph_iij}).

When $s_{1}=s_{2}=s>0$, (\ref{proof:circular_2}) shows that
\begin{align}\label{proof:circular_2_s1=s2}
\mathbb{E}\left[W^{R}_{\Phi}(t,s)W^{I}_{\Phi}(t,s)\right] = 0.
\end{align}
The second equality in Proposition \ref{lemma:exp}(c) implies that
\begin{align*}
\textup{Var}\left(S_{\Phi}(t,s)\right)
=4\left(\mathbb{E}\left[|W^{R}_{\Phi}(t,s)|^{2}\right]\right)^{2}.
\end{align*}
The proof of the second equality in Proposition \ref{lemma:exp}(c) is concluded by noting that
\begin{align}\notag
\mathbb{E}[S_{\Phi}(t,s)] = \mathbb{E}\left[|W^{R}_{\Phi}(t,s)|^{2}\right]+\mathbb{E}\left[|W^{I}_{\Phi}(t,s)|^{2}\right]
\overset{(\ref{proof:circular_1})}{=}2\mathbb{E}\left[|W^{R}_{\Phi}(t,s)|^{2}\right].
\end{align}
\qed

\subsection{Proof of Proposition \ref{lemma:twice_differentiability}}\label{sec:proof:lemma:twice_differentiability}
The proof consists of three steps. First, we show that the sample paths of $W_{\Phi}$
are almost surely continuous on $\mathbb{R}\times(0,\infty)$.
Next, we prove that the sample paths of $W_{\Phi}$ are continuously differentiable on $\mathbb{R}\times(0,\infty)$
almost surely. Finally, the second-order differentiability of the sample paths of $W_{\Phi}$
follows by the same arguments.

{\it Step 1.} To prove that the sample paths of the Gaussian field $W_{\Phi}$ are almost surely continuous on $\mathbb{R}\times(0,\infty)$,
according to \cite[p.73 and p.93]{dudley1973sample}, it suffices to show the existence
of a constant $q>0$ such that
for any $\delta>0$, there exists a constant $C_{\delta}$ satisfying
\begin{align}\label{claim:continuity_realW_new}
\left\{\mathbb{E}\left[\left|W_{\Phi}(t',s')-
W_{\Phi}(t,s)\right|^{2}\right]\right\}^{1/2}
\leq C_{\delta}\left(|t-t'|^{2}+|s-s'|^{2}\right)^{q/2}
\end{align}
for any $t,t'\in \mathbb{R}$ and $s,s'>\delta$.
We prove (\ref{claim:continuity_realW_new}) as follows.
By the triangle inequality and (\ref{spect:WPhi}),
\begin{align}\notag
&\left\{\mathbb{E}\left[\left|W_{\Phi}(t',s')-
W_{\Phi}(t,s)\right|^{2}\right]\right\}^{1/2}
\\\notag\leq&
\left\{\mathbb{E}\left[\left|W_{\Phi}(t',s')-
W_{\Phi}(t',s)\right|^{2}\right]\right\}^{1/2}
+\left\{\mathbb{E}\left[\left|
W_{\Phi}(t',s)-
W_{\Phi}(t,s)\right|^{2}\right]\right\}^{1/2}
\\\label{proof:continuity:e1+e2}=&
\sqrt{\mathcal{E}_{1}}+\sqrt{\mathcal{E}_{2}},
\end{align}
where
\begin{align*}
\mathcal{E}_{1} = \int_{\mathbb{R}}
|\widehat{\psi}(s'\lambda)- \widehat{\psi}(s\lambda)|^{2}p(\lambda)d\lambda
\end{align*}
and
\begin{align*}
\mathcal{E}_{2} = \int_{\mathbb{R}}
|e^{i(t'-t)\lambda}-1|^{2} | \widehat{\psi}\left(s\lambda\right)|^{2}p(\lambda)d\lambda.
\end{align*}

{\it Estimate of $\mathcal{E}_{1}$:} By the mean value theorem, there exists a constant $c_{s,s',\lambda}$ lying between $s\lambda$ and $s'\lambda$ such that
\begin{align}\notag
\mathcal{E}_{1}
=\int_{0}^{\infty}
|D\widehat{\psi}(c_{s,s',\lambda})(s'\lambda-s\lambda)|^{2}p(\lambda)d\lambda.
\end{align}
Because $s,s'>\delta$, we have $c_{s,s',\lambda}>\delta\lambda$ and
\begin{align}\notag
\mathcal{E}_{1}
\leq &\delta^{-2}|s'-s|^{2}\int_{0}^{\infty}
|D\widehat{\psi}(c_{s,s',\lambda})c_{s,s',\lambda}|^{2}p(\lambda)d\lambda
\\\label{s_near_field_behavior1_new}\leq& \delta^{-2}|s'-s|^{2}\left(\underset{z>0}{\sup}|zD\widehat{\psi}(z)|^{2}\right)\mathbb{E}\left[\Phi^{2}(0)\right],
\end{align}
where the supremum is finite under Assumptions \ref{assumption:boundedness:psi}($\mathrm{D}^{1}_{1}$).

{\it Estimate of $\mathcal{E}_{2}$:}
Given that $|e^{ix}-1|\leq |x|$ for any $x\in \mathbb{R}$, and $s>\delta$,
\begin{align}\label{t_near_field_behavior1new}
\mathcal{E}_{2} \leq |t'-t|^{2}\int_{\mathbb{R}}
|\widehat{\psi}\left(s\lambda\right)|^{2} \lambda^{2} p(\lambda)d\lambda
\leq \delta^{-2}|t'-t|^{2}\left(\underset{z>0}{\sup}|z\widehat{\psi}(z)|^{2}\right)\mathbb{E}\left[|\Phi(0)|^{2}\right],
\end{align}
where the supremum is finite under Assumption \ref{assumption:boundedness:psi}\textup{($\mathrm{D}^{0}_{1}$)}.

Substituting (\ref{s_near_field_behavior1_new}) and (\ref{t_near_field_behavior1new}) into (\ref{proof:continuity:e1+e2}) shows that
(\ref{claim:continuity_realW_new}) holds with $q=1$ and
\begin{align*}
C_{\delta}= \delta^{-1}\left\{\left(\underset{z>0}{\sup}|z\widehat{\psi}(z)|^{2}\right)
+
\left(\underset{z>0}{\sup}|zD\widehat{\psi}(z)|^{2}\right)
\right\}^{1/2}
\left\{\mathbb{E}\left[\Phi^{2}(0)\right]\right\}^{1/2}.
\end{align*}

{\it Step 2.}
Motivated from (\ref{spect:WPhi}) , we consider two random fields
\begin{align}\label{predef:DWPhi_t}
\partial_{t} W_{\Phi}(t,s) =  i\int_{0}^{\infty}
\lambda e^{it\lambda}  \overline{\widehat{\psi}(s\lambda)} \sqrt{p(\lambda)}Z(d\lambda),
\end{align}
\begin{align}\label{predef:DWPhi}
\partial_{s}W_{\Phi}(t,s) =  \int_{0}^{\infty}\lambda
e^{it\lambda} \overline{[D\widehat{\psi}](s\lambda)} \sqrt{p(\lambda)}Z(d\lambda),
\end{align}
where $D\widehat{\psi}$ represents the derivative of $\widehat{\psi}$.
Under Assumption \ref{assumption:boundedness:psi},
conditions \textup{($\mathrm{D}^{0}_1$)} and \textup{($\mathrm{D}^{1}_1$)} ensure that
$$\mathbb{E}|\partial_{t}W_{\Phi}(t,s)|^{2}<\infty\ \ \textup{and}\  \ \mathbb{E}|\partial_{s}W_{\Phi}(t,s)|^{2}<\infty.$$
To prove that sample paths of $W_{\Phi}$ are continuously differentiable almost surely, according to \cite{potthoff2010sample}, it suffices to show that
(a) the difference quotient of $W_{\Phi}$ along the time direction converges to
$\partial_{t} W_{\Phi}$ in quadratic mean,
(b) the difference quotient of $W_{\Phi}$ along the scale direction converges to
$\partial_{s} W_{\Phi}$ in quadratic mean, and
(c) both $\partial_{t} W_{\Phi}$ and $\partial_{s}W_{\Phi}$ are continuous in quadratic mean
and almost surely have continuous realizations.
We verify these conditions as follows.

{\it Proof of (a):}
By (\ref{spect:WPhi}) and (\ref{predef:DWPhi_t}), for any $h\in \mathbb{R}\setminus\{0\}$,
\begin{align}\notag
&\frac{W_{\Phi}(t+h,s)-W_{\Phi}(t,s)}{h}-\partial_{t} W_{\Phi}(t,s)
\\\notag=&  \int_{0}^{\infty}
e^{it\lambda} \left(\frac{e^{ih\lambda}-1}{h}-i\lambda\right)
\overline{\widehat{\psi}(s\lambda)}\sqrt{p(\lambda)}Z(d\lambda).
\end{align}
By the isometry property of the Wiener integral,
\begin{align}\notag
&\mathbb{E}\left[\left|\frac{W_{\Phi}(t+h,s)-W_{\Phi}(t,s)}{h}-\partial_{t} W_{\Phi}(t,s)\right|^{2}\right]
\\\label{quotientdiff_time}=& \int_{0}^{\infty}
\left|\frac{e^{ih\lambda}-1}{h}-i\lambda\right|^{2}|\widehat{\psi}(s\lambda)|^{2}p(\lambda)d\lambda.
\end{align}
Note that for any $h\in \mathbb{R}\setminus\{0\}$ and $\lambda\in \mathbb{R}$,
\begin{align}\label{quotient:exp}
\left|\frac{e^{ih\lambda}-1}{h}-i\lambda\right|\leq 2|\lambda|.
\end{align}
Assumption \ref{assumption:boundedness:psi}\textup{($\mathrm{D}^{0}_1$)}, together with (\ref{quotient:exp}), enables us to apply the dominated convergence theorem to (\ref{quotientdiff_time}) and obtain
\begin{align*}
\underset{h\rightarrow 0}{\lim}\ \mathbb{E}\left[\left|\frac{W_{\Phi}(t+h,s)-W_{\Phi}(t,s)}{h}-\partial_{t} W_{\Phi}(t,s)\right|^{2}\right]
=  0,
\end{align*}
which shows that the difference quotient of $W_{\Phi}$ along the time direction converges to
$\partial_{t}W_{\Phi}$ in quadratic mean.

{\it Proof of (b):}
By (\ref{spect:WPhi}) and (\ref{predef:DWPhi}), for any $h\in (-s/2,s/2)\setminus\{0\}$,
\begin{align}\notag
&\frac{W_{\Phi}(t,s+h)-W_{\Phi}(t,s)}{h}-\partial_{s}W_{\Phi}(t,s)
\\\notag=&  \int_{0}^{\infty}
e^{it\lambda} \left[
\frac{\overline{\widehat{\psi}((s+h)\lambda)}-\overline{\widehat{\psi}(s\lambda)}}{h}
-\lambda \overline{[D\widehat{\psi}](s\lambda)}
\right]\sqrt{p(\lambda)}Z(d\lambda).
\end{align}
By the isometry property of the Wiener integral,
\begin{align}\label{quotientdiff}
\mathbb{E}\left[\left|\frac{W_{\Phi}(t,s+h)-W_{\Phi}(t,s)}{h}-\partial_{s}W_{\Phi}(t,s)\right|^{2}\right]
=  \int_{0}^{\infty}
\left|K_{h}(s,\lambda)\right|^{2}p(\lambda)d\lambda,
\end{align}
where
\begin{align*}
K_{h}(s,\lambda) =
\frac{\widehat{\psi}((s+h)\lambda)-\widehat{\psi}(s\lambda)}{h}
-\lambda [D\widehat{\psi}](s\lambda).
\end{align*}
For any $\lambda>0$ and $h\in (-s/2,s/2)\setminus\{0\}$,
\begin{align}\notag
\left|\frac{\widehat{\psi}((s+h)\lambda)-\widehat{\psi}(s\lambda)}{h}\right|
=& \left|\lambda\int_{0}^{1} [D\widehat{\psi}]\left((s+h)\lambda u+s\lambda(1-u)\right)du\right|
\\\notag\leq &\left(\underset{z>0}{\sup}\left|zD\widehat{\psi}(z)\right|\right)\lambda \int_{0}^{1} \left|\frac{1}{(s+h)\lambda u+s\lambda(1-u)}\right| du
\\\notag=&\left(\underset{z>0}{\sup}\left|zD\widehat{\psi}(z)\right|\right)  \frac{1}{|h|}\ln\left|1+\frac{h}{s}\right|
\\\notag\leq&  \left(\underset{z>0}{\sup}\left|zD\widehat{\psi}(z)\right|\right)\left(\frac{2\ln2}{s}\right),
\end{align}
where the supremum is finite under Assumption \ref{assumption:boundedness:psi}\textup{($\mathrm{D}^{1}_1$)}.
Hence, for any $s>0$,
\begin{align}\label{bound:Kh}
\underset{h\in (-s/2,s/2)\setminus\{0\}}{\sup}\  \underset{\lambda>0}{\sup}\ |K_{h}(s,\lambda)| <\left(\underset{z>0}{\sup}\left|zD\widehat{\psi}(z)\right|\right)\left(\frac{1+2\ln2}{s}\right)<\infty.
\end{align}
The differentiability of $\widehat{\psi}$ implies that $K_{h}(s,\lambda)\rightarrow 0$ as $h\rightarrow0$.
The uniform bound (\ref{bound:Kh}) enables us to apply the dominated convergence theorem to (\ref{quotientdiff}),
yielding
\begin{align*}
\underset{h\rightarrow 0}{\lim}\
\mathbb{E}\left[\left|\frac{W_{\Phi}(t,s+h)-W_{\Phi}(t,s)}{h}-\partial_{s}W_{\Phi}(t,s)\right|^{2}\right]
=  0.
\end{align*}
This result shows that the difference quotient of $W_{\Phi}$ along the scale direction converges to
$\partial_{s}W_{\Phi}$ in quadratic mean.

{\it Proof of (c):} To prove that $\partial_{t} W_{\Phi}$ and $\partial_{s} W_{\Phi}$
not only are continuous in the quadratic mean but also almost surely have modification that is sample continuous,
according to \cite[p.73 and p.93]{dudley1973sample} (see also \cite[Corollary 4.6]{potthoff2009sample}),
it suffices to show that
for any $\delta>0$, there exists a constant $C_{\delta}$ such that
\begin{align}\notag
\left\{\mathbb{E}\left[\left|\partial_{t} W_{\Phi}(t',s')-
\partial_{t} W_{\Phi}(t,s)\right|^{2}\right]\right\}^{1/2}
\leq C_{\delta}\left(|t-t'|^{2}+|s-s'|^{2}\right)^{1/2}
\end{align}
and
\begin{align}\notag
\left\{\mathbb{E}\left[\left|\partial_{s} W_{\Phi}(t',s')-
\partial_{s} W_{\Phi}(t,s)\right|^{2}\right]\right\}^{1/2}
\leq C_{\delta}\left(|t-t'|^{2}+|s-s'|^{2}\right)^{1/2}
\end{align}
for any $t,t'\in \mathbb{R}$ and $s,s'>\delta$.
This step requires additional conditions, including conditions \textup{($\mathrm{D}^{0}_2$)}, \textup{($\mathrm{D}^{1}_2$)} and \textup{($\mathrm{D}^{2}_2$)} in Assumption \ref{assumption:boundedness:psi}.
Since the proof is similar to the argument used to prove the continuity of the sample paths of $W_{\Phi}$
on $\mathbb{R}\times (0,\infty)$, it is omitted here.
Based on the verification process (a)-(c), we conclude that the sample paths of $W_{\Phi}$ are continuously differentiable almost surely.

{\it Step 3.} The twice continuous differentiability of $W_{\Phi}$ follows by applying the same arguments, i.e.,
the verification process (a)-(c),
to the derivatives $\partial_{t}W_{\Phi}$ and $\partial_{s}W_{\Phi}$.
This part of proof relies on conditions \textup{($\mathrm{D}^{2}_3$)} and
\textup{($\mathrm{D}^{3}_3$)} in Assumption \ref{assumption:boundedness:psi}.
\qed

\subsection{Proof of Corollary \ref{lemma:continuity_extension}}\label{sec:proof:lemma:continuity_extension}

{\it  Proof of (\ref{proof:continuity_0}):}
Based on the continuity of the sample paths of $W_{\Phi}$ on $\mathbb{R}\times (0,\infty)$,
it suffices to show that the sample paths of the Gaussian field $W_{\Phi}$
are almost surely continuous at the boundary $\mathbb{R}\times\{0\}$ if we define
$W_{\Phi}(t,0)=0$ for any $t\in \mathbb{R}$.
According to \cite[p.73 and p.93]{dudley1973sample}, it suffices to show the existence
of a constant $q>0$ such that there exists a constant $C$ satisfying
\begin{align}\label{proof:continuity0a}
&\mathbb{E}\left[\left|W_{\Phi}(t,s)\right|^{2}\right] = \int_{0}^{\infty}
|\widehat{\psi}(s\lambda)|^{2}p(\lambda)d\lambda\leq C s^{q}
\end{align}
for any $t\in \mathbb{R}$ and sufficiently small $s>0$.
Under Assumption \ref{assumption:boundedness:psi}($\mathrm{D}^{1}_{0}$), we have
$|\widehat{\psi}(s\lambda)|\leq \|D\widehat{\psi}\|_{\infty}|s\lambda|$.
Hence,
\begin{align}\label{proof:continuity0b}
\int_{0}^{s^{-2/(2+\gamma)}}
|\widehat{\psi}(s\lambda)|^{2}p(\lambda)d\lambda
\leq
\|D\widehat{\psi}\|_{\infty}^{2}\mathbb{E}[\Phi^{2}(0)] s^{\frac{2\gamma}{2+\gamma}}.
\end{align}
On the other hand, by Assumption \ref{assumption:spectral_density}(a),
\begin{align}\label{proof:continuity0c}
\int_{s^{-2/(2+\gamma)}}^{\infty}
|\widehat{\psi}(s\lambda)|^{2}p(\lambda)d\lambda
\leq \|\widehat{\psi}\|_{\infty}^{2}C_{1}\gamma^{-1}s^{\frac{2\gamma}{2+\gamma}}.
\end{align}
By combining (\ref{proof:continuity0b}) and (\ref{proof:continuity0c}), we obtain that
\begin{align}\notag
&\mathbb{E}\left[\left|W_{\Phi}(t,s)\right|^{2}\right]\leq   \left[\|D\widehat{\psi}\|_{\infty}^{2}\mathbb{E}[\Phi^{2}(0)]+\|\widehat{\psi}\|_{\infty}^{2}C_{1}\gamma^{-1}\right]s^{\frac{2\gamma}{2+\gamma}}.
\end{align}
Hence, (\ref{proof:continuity0a}) holds with $q=2\gamma/(2+\gamma)$.

{\it Proof of  (\ref{proof:continuity_00a}):} We first introduce a Gaussian random field $\varphi:\mathbb{R}\times[0,\infty)\mapsto \mathbb{C}$
as follows
\begin{equation*}
\varphi(t,\tau) =\left\{
\begin{array}{ll}
W_{\Phi}(t,1/\tau) &\ \textup{for}\ \tau>0,
\\
0 &\ \textup{for}\ \tau=0.
\end{array}
\right.
\end{equation*}
Similarly to the proof of (\ref{proof:continuity_0}), proving (\ref{proof:continuity_00a}) is equivalent to showing that
the function $\varphi$ is continuous at the boundary $\mathbb{R}\times \{0\}$.
For any $t\in \mathbb{R}$ and $\tau>0$,
\begin{align}\label{proof:continuity_00b}
\mathbb{E}\left[\left|\varphi(t,\tau)\right|^{2}\right]
=\int_{0}^{\infty}
|\widehat{\psi}(\frac{\lambda}{\tau})|^{2}p(\lambda)d\lambda
=\tau\int_{0}^{\infty}
|\widehat{\psi}(u)|^{2}p(\tau u)du.
\end{align}
Under Assumption \ref{assumption:spectral_density}(b1),
(\ref{proof:continuity_00b}) can be rewritten as
\begin{align*}
\mathbb{E}\left[\left|\varphi(t,\tau)\right|^{2}\right]=
\tau\int_{0}^{\infty}
|\widehat{\psi}(u)|^{2}\frac{L(|\tau u|^{-1})}{|\tau u|^{1-H}}du.
\end{align*}
By \cite[Theorem 1.5.3]{bingham1989regular},
\begin{align*}
\underset{\tau\rightarrow 0+}{\lim}\left[L(\tau^{-1})\tau^{H}\right]^{-1}\mathbb{E}\left[\left|\varphi(t,\tau)\right|^{2}\right]
=
\underset{\tau\rightarrow 0+}{\lim}\int_{0}^{\infty}
\frac{|\widehat{\psi}(u)|^{2}}{|u|^{1-H}}\frac{L(|\tau u|^{-1})}{L(\tau^{-1})}du
=\int_{0}^{\infty}
\frac{|\widehat{\psi}(u)|^{2}}{|u|^{1-H}}du.
\end{align*}
Hence, $\mathbb{E}\left[\left|\varphi(t,\tau)\right|^{2}\right]\lesssim L(\tau^{-1})\tau^{H}$ as $\tau$ is sufficiently small.
On the other hand, under Assumption \ref{assumption:spectral_density}(b2),
it can similarly be shown that $\mathbb{E}\left[\left|\varphi(t,\tau)\right|^{2}\right]\lesssim \tau$
as $\tau$ is sufficiently small.
Therefore, the sample paths of $\varphi$ are continuous at the boundary $\mathbb{R}\times\{0\}$ almost surely.
\qed

\subsection{Proof of Theorem \ref{prop:unique_argmax}}\label{sec:proof:prop:unique_argmax}
(a) Because the time variable is fixed in the statement of Theorem \ref{prop:unique_argmax}, we omit it in the following derivation. That is, we denote $W_{Y}(s):=W_{Y}(t,s)$ and $S_{Y}(s):=S_{Y}(t,s)$ for a given time point $t$.
Next, we define the following matrices
\begin{align*}
\mathbf{W}_{f}(s) = \left[W^{R}_{f}(s)\ W^{I}_{f}(s)\right]^{\top},
\mathbf{W}_{\Phi}(s) = \left[W^{R}_{\Phi}(s)\ W^{I}_{\Phi}(s)\right]^{\top},
\mathbf{W}_{Y}(s) = \left[W^{R}_{Y}(s)\ W^{I}_{Y}(s)\right]^{\top},
\end{align*}
where $W^{R}_{f}(s) = \textup{Re}(W_{f}(s))$, $W^{I}_{f}(s) = \textup{Im}(W_{f}(s))$,
and similarly for other terms.
Clearly, $\mathbf{W}_{Y}(s) =\mathbf{W}_{f}(s) +\mathbf{W}_{\Phi}(s)$ and
$S_{Y}(s) = \mathbf{W}_{Y}^{\top}(s)\mathbf{W}_{Y}(s)$. To prove the uniqueness of the global maximizer of $S_{Y}$, it suffices to show that for any different $s_{1},s_{2}>0$,
there exist neighborhoods $N(s_{1})$ and $N(s_{2})$ of $s_{1}$ and $s_{2}$ such that
\begin{align}\label{claim:sup_vs_sup}
\mathbb{P}\left(\underset{s\in N(s_{1})}{\sup}S_{Y}(s) =
\underset{s\in N(s_{2})}{\sup}S_{Y}(s)\right)=0
\end{align}
because
\begin{align}\notag
&\mathbb{P}\left(S_{Y}(s_{1}) = \underset{s>0}{\sup}\ S_{Y}(s) = S_{Y}(s_{2})\ \textup{for some different}\ s_{1},s_{2}>0\right)
\\\notag\leq&  \mathbb{P}\left(\underset{s\in N(s_{1})}{\sup}S_{Y}(s) =
\underset{s\in N(s_{2})}{\sup}S_{Y}(s)\ \textup{for some different}\ s_{1},s_{2}\in \mathbb{Q}_{>0}\right)
\\\label{union_bound:Q}\leq&  \underset{\begin{subarray}{c} s_{1},s_{2}\in \mathbb{Q}_{>0}\\ s_{1}\neq s_{2}\end{subarray}}{\sum} \mathbb{P}\left(\underset{s\in N(s_{1})}{\sup}S_{Y}(s) =
\underset{s\in N(s_{2})}{\sup}S_{Y}(s)\right),
\end{align}
where $\mathbb{Q}_{>0}$ represents the set of positive rational numbers, and the first inequality follows from the property that, almost surely, the sample path of
$S_{Y}$ is continuous, as demonstrated in Proposition \ref{lemma:twice_differentiability}.

Denote the covariance matrix function of the process $\mathbf{W}_{Y}$ by
$\Gamma$.
Because $\mathbb{E}\left[\mathbf{W}_{Y}(s)\right] =\mathbf{W}_{f}(s)$ for any $s>0$,
\begin{align*}
\Gamma(s_{1},s_{2}) =\mathbb{E}\left[\mathbf{W}_{\Phi}(s_{1})\mathbf{W}_{\Phi}(s_{2})^{\top}\right]= \left[\begin{array}{cc}\mathbb{E}\left[W^{R}_{\Phi}(s_{1})W^{R}_{\Phi}(s_{2})\right]
& \mathbb{E}\left[W^{R}_{\Phi}(s_{1})W^{I}_{\Phi}(s_{2})\right] \\
\mathbb{E}\left[W^{I}_{\Phi}(s_{1})W^{R}_{\Phi}(s_{2})\right] &
\mathbb{E}\left[W^{I}_{\Phi}(s_{1})W^{I}_{\Phi}(s_{2})\right]\end{array}\right].
\end{align*}
By (\ref{proof:circular_1}) and (\ref{proof:circular_2_s1=s2}), for $k\in \{1,2\}$,
$$\Gamma(s_{k},s_{k}) =\mathbb{E}\left[|W^{R}_{\Phi}(s_{k})|^{2}\right]\mathbf{I}_{2\times 2}=2^{-1}\sqrt{\textup{Var}\left(S_{\Phi}(s_{k})\right)}\mathbf{I}_{2\times 2},$$
where $\mathbf{I}_{2\times 2}$ is the $2\times2$ identity matrix.
By the condition that the spectral density $p$ is positive almost everywhere in Assumption \ref{assumption:Gaussian}, we have $\textup{Var}\left(S_{\Phi}(s_{1})\right)>0$  and $\textup{Var}\left(S_{\Phi}(s_{2})\right)>0$.
Without loss of generality, we assume that
\begin{align}\label{var1>var2}
\textup{Var}\left(S_{\Phi}(s_{1})\right)> \textup{Var}\left(S_{\Phi}(s_{2})\right)>0.
\end{align}
Denote
\begin{align*}
\frac{1}{\Gamma(s_{1},s_{1})} = \left\{\mathbb{E}\left[|W^{R}_{\Phi}(s_{1})|^{2}\right]\right\}^{-1}
\end{align*}
and consider the following decomposition. For any $s>0$,
\begin{align*}
\mathbf{W}_{Y}(s) =  \mathbf{W}_{\perp}(s)+\frac{\Gamma(s,s_{1})}{\Gamma(s_{1},s_{1})}
\mathbf{W}_{\Phi}(s_{1}),
\end{align*}
where
\begin{align*}
\mathbf{W}_{\perp}(s) =\mathbf{W}_{f}(s)+\mathbf{W}_{\Phi}(s)-\frac{\Gamma(s,s_{1})}{\Gamma(s_{1},s_{1})}
\mathbf{W}_{\Phi}(s_{1}).
\end{align*}
Let $\delta_{1},\delta_{2}>0$ that will be determined later.
For $k\in\{1,2\}$ and any $\mathbf{z}=[z_{1}\ z_{2}]^{\top}\in \mathbb{R}^{2}$,
we define
\begin{align*}
M_{k}(\mathbf{z}) = \underset{s\in (s_{k}-\delta_{k},s_{k}+\delta_{k})}{\sup}
Q(\mathbf{z};s),
\end{align*}
where $Q(\cdot;s)$ is a quadratic function of $\mathbf{z}$ defined as follows
\begin{align*}
Q(\mathbf{z};s) = \left(\mathbf{W}_{\perp}(s)+\frac{\Gamma(s,s_{1})}{\Gamma(s_{1},s_{1})}
\mathbf{z}\right)^{\top}\left(\mathbf{W}_{\perp}(s)+\frac{\Gamma(s,s_{1})}{\Gamma(s_{1},s_{1})}
\mathbf{z}\right).
\end{align*}
By noting that
\begin{align*}
M_{k}\left(\mathbf{W}_{\Phi}(s_{1})\right) = \underset{s\in (s_{k}-\delta_{k},s_{k}+\delta_{k})}{\sup}
Q(\mathbf{W}_{\Phi}(s_{1});s)= \underset{s\in (s_{k}-\delta_{k},s_{k}+\delta_{k})}{\sup}
S_{Y}(s),
\end{align*}
the left-hand side of (\ref{claim:sup_vs_sup}) with $N(s_{1})= (s_{1}-\delta_{1},s_{1}+\delta_{1})$
and $N(s_{2})= (s_{2}-\delta_{2},s_{2}+\delta_{2})$
can be rewritten as
\begin{align}\notag
\mathbb{P}\left(\underset{s\in (s_{1}-\delta_{1},s_{1}+\delta_{1})}{\sup}\hspace{-0.5cm}S_{Y}(s) =
\underset{s\in (s_{2}-\delta_{2},s_{2}+\delta_{2})}{\sup}\hspace{-0.5cm}S_{Y}(s)\right)
=
\mathbb{E}\left[\mathbb{P}\left(M_{1}(\mathbf{W}_{\Phi}(s_{1})) =
M_{2}(\mathbf{W}_{\Phi}(s_{1}))\mid \mathbf{W}_{\perp}\right) \right].
\end{align}
Thus, to prove (\ref{claim:sup_vs_sup}), it suffices to show that conditioned on $\mathbf{W}_{\perp}$,
\begin{equation}\label{proof:unique:equiv}
\mathbb{P}\left(M_{1}(\mathbf{W}_{\Phi}(s_{1})) =
M_{2}(\mathbf{W}_{\Phi}(s_{1}))\mid \mathbf{W}_{\perp}\right)=0.
\end{equation}
Our strategy is to show that the Lebesgue measure of the set
$$\mathcal{I}:=\{\mathbf{z}\in\mathbb{R}^{2}\mid  M_{1}(\mathbf{z}) = M_{2}(\mathbf{z})\}$$
is zero.

By a direct expansion,
\begin{align}\notag
Q(\mathbf{z};s) =
\mathbf{W}_{\perp}^{\top}(s)\mathbf{W}_{\perp}(s)
+2\mathbf{W}^{\top}_{\perp}(s)\frac{\Gamma(s,s_{1})}{\Gamma(s_{1},s_{1})}\mathbf{z}
+\mathbf{z}^{\top}\frac{\Gamma(s,s_{1})^{\top}}{\Gamma(s_{1},s_{1})}\frac{\Gamma(s,s_{1})}{\Gamma(s_{1},s_{1})}
\mathbf{z}.
\end{align}
By (\ref{proof:circular_1}) and (\ref{proof:circular_2}),
\begin{align}\notag
\frac{\Gamma(s,s_{1})^{\top}}{\Gamma(s_{1},s_{1})}\frac{\Gamma(s,s_{1})}{\Gamma(s_{1},s_{1})}
=&\frac{\left(\mathbb{E}\left[W^{R}_{\Phi}(s)W^{R}_{\Phi}(s_{1})
\right]\right)^{2}+ \left(\mathbb{E}\left[W^{I}_{\Phi}(s)W^{R}_{\Phi}(s_{1})
\right]\right)^{2}}{\left(\mathbb{E}\left[W^{R}_{\Phi}(s_{1})W^{R}_{\Phi}(s_{1})
\right]\right)^{2}}\mathbf{I}_{2\times 2}
\\\notag=&
\frac{\textup{Cov}(S_{\Phi}(s),S_{\Phi}(s_{1}))}{\textup{Var}\left(S_{\Phi}(s_{1})\right)}\mathbf{I}_{2\times 2},
\end{align}
where the second equality stems from Proposition \ref{lemma:exp}(c).
Therefore, for any $s>0$,
\begin{align}\notag
Q(\mathbf{z};s) =
\mathbf{W}_{\perp}^{\top}(s)\mathbf{W}_{\perp}(s)
+2\mathbf{W}^{\top}_{\perp}(s)\frac{\Gamma(s,s_{1})}{\Gamma(s_{1},s_{1})}\mathbf{z}
+
\frac{\textup{Cov}(S_{\Phi}(s),S_{\Phi}(s_{1}))}{\textup{Var}\left(S_{\Phi}(s_{1})\right)}
\|\mathbf{z}\|^{2},
\end{align}
where $\|\mathbf{z}\|^{2}=\mathbf{z}^{\top}\mathbf{z}$.
By the result $\textup{Cov}\left(S_{\Phi}(s_{1}),S_{\Phi}(s_{2})\right)\geq 0$ given in Proposition \ref{lemma:exp}(c),
the Cauchy-Schwarz inequality, the condition in Assumption \ref{assumption:Gaussian} that the spectral density $p$ is positive almost everywhere, and (\ref{var1>var2}), we have
\begin{align}\label{cov11>cov21}
0\leq \textup{Cov}(S_{\Phi}(s_{2}),S_{\Phi}(s_{1}))<\textup{Var}\left(S_{\Phi}(s_{1})\right)=\textup{Cov}(S_{\Phi}(s_{1}),S_{\Phi}(s_{1})).
\end{align}
The continuity of the function $s\mapsto \textup{Cov}(S_{\Phi}(s),S_{\Phi}(s_{1}))$ and (\ref{cov11>cov21})
ensure that
there exist constants $c_{1}>c_{2}>0$ and sufficiently small $\delta_{1},\delta_{2}>0$ such that
the quadratic coefficients of $\{Q(\cdot;s)\mid s\in (s_{1}-\delta_{1},s_{1}+\delta_{1})\}$ and $\{Q(\cdot;s)\mid s\in (s_{2}-\delta_{2},s_{2}+\delta_{2})\}$ can be separated as follows
\begin{align}\label{separate_hessian}
\underset{s\in (s_{1}-\delta_{1},s_{1}+\delta_{1})}{\sup}\frac{\textup{Cov}(S_{\Phi}(s,s_{1}))}{\textup{Var}\left(S_{\Phi}(s_{1})\right)}
>c_{1}>c_{2}>\underset{s\in (s_{2}-\delta_{2},s_{2}+\delta_{2})}{\sup}\frac{\textup{Cov}(S_{\Phi}(s,s_{1}))}{\textup{Var}\left(S_{\Phi}(s_{1})\right)}.
\end{align}
Furthermore, for any fixed $s>0$, the graph of $Q(\cdot,s)$ can be categorized as follows:
\begin{itemize}
\item For $s$ satisfying $\textup{Cov}(S_{\Phi}(s),S_{\Phi}(s_{1}))>0$,
 the surface $\{( \mathbf{z},Q(\mathbf{z},s))\mid \mathbf{z}\in \mathbb{R}^{2}\}$ is an upward-opening paraboloid;
\item For $s$ satisfying $\textup{Cov}(S_{\Phi}(s),S_{\Phi}(s_{1}))=0$, by Proposition \ref{lemma:exp}(c),
we have
$\mathbb{E}\left[W^{R}_{\Phi}(s)W^{R}_{\Phi}(s_{1})\right]=0$
and
$\mathbb{E}\left[W^{R}_{\Phi}(s)W^{I}_{\Phi}(s_{1})\right]=0,$
which imply that $\Gamma(s,s_{1})=\mathbf{0}_{2\times 2}$ and $Q(\mathbf{z};s)=\mathbf{W}_{\perp}^{\top}(s)\mathbf{W}_{\perp}(s)$ for any $\mathbf{z}\in \mathbb{R}^{2}$.
This indicates that $\{( \mathbf{z},Q(\mathbf{z},s))\mid \mathbf{z}\in \mathbb{R}^{2}\}$ is a horizontal plane.
\end{itemize}

By the fact that the supremum of a family of convex functions is still convex,
both $M_{1}$ and $M_{2}$ are convex functions defined on $\mathbb{R}^{2}$.
According to Alexandrov's theorem
\cite[Theorem 2.1]{rockafellar1999second}, $M_{1}$ and $M_{2}$ have the second derivatives almost everywhere.
More precisely, there exists a set $D\subset \mathbb{R}^{2}$, whose complement has Lebesgue measure zero, such that
for $k\in\{1,2\}$ and every point $\mathbf{z}\in D$, there is a
quadratic expansion in the form
\begin{align}\label{second_derivative_convex}
M_{k}(\mathbf{z}+\mathbf{w}) = M_{k}(\mathbf{z})+\mathbf{v}_{k}(\mathbf{z})^{\top}\mathbf{w}
+\frac{1}{2}\mathbf{w}^{\top}\mathbf{H}_{k}(\mathbf{z})\mathbf{w}+o(\|\mathbf{w}\|^{2}),\ \mathbf{w}\in \mathbb{R}^{2},
\end{align}
where $\mathbf{v}_{k}(\mathbf{z})$ is the gradient, and $\mathbf{H}_{k}(\mathbf{z})$ is the Hessian matrix of  $M_{k}$ at $\mathbf{z}$.
The function $o(\|\mathbf{w}\|^{2})$ satisfies $o(\|\mathbf{w}\|^{2})\|\mathbf{w}\|^{-2}\rightarrow 0$ as $\|\mathbf{w}\|\rightarrow 0$.
By (\ref{separate_hessian}) (see also \cite[Example 13.10]{rockafellar2009variational}),
\begin{align}\label{second_derivative_convex0}
\frac{1}{2}\mathbf{w}^{\top}\mathbf{H}_{1}(\mathbf{z})\mathbf{w}\geq c_{1}\|\mathbf{w}\|^{2}
>c_{2}\|\mathbf{w}\|^{2}\geq \frac{1}{2}\mathbf{w}^{\top}\mathbf{H}_{2}(\mathbf{z})\mathbf{w},
\end{align}
which implies that
\begin{align}\label{lebesgue_zero_intersection}
\textup{Leb}\left(\mathcal{I}\right)=0,
\end{align}
where $\textup{Leb}$ represents the Lebesgue measure on $\mathbb{R}^{2}$. For a further explanation of (\ref{lebesgue_zero_intersection}), see Remark \ref{remark:Leb0} at the end of this section.
From the covariance computation, the normal random vector $\mathbf{W}_{\Phi}(s_{1})$ is independent of the Gaussian random process $\mathbf{W}_{\perp}$. Thus, conditioning on $\mathbf{W}_{\perp}$, the probability density of
$\mathbf{W}_{\Phi}(s_{1})$ remains absolutely continuous with respect to the Lebesgue measure on $\mathbb{R}^{2}$.
Therefore, the observation (\ref{lebesgue_zero_intersection}) implies that
\begin{align}\notag
\mathbb{P}\left(M_{1}\left(\mathbf{W}_{\Phi}(s_{1})\right)=M_{2}\left(\mathbf{W}_{\Phi}(s_{1})\right)\mid \mathbf{W}_{\perp}\right)
=
\mathbb{P}\left(\mathbf{W}_{\Phi}(s_{1})\in \mathcal{I}\mid \mathbf{W}_{\perp}\right)=0,
\end{align}
which establishes (\ref{proof:unique:equiv}).

(b) Regarding the second part of Theorem \ref{prop:unique_argmax}, similar to (\ref{union_bound:Q}),
for any $t\in \mathbb{R}$ and $m\in \{1,\ldots,M\}$, we have
\begin{align}\notag
&\mathbb{P}\left(s_{Y,m}(t)\ \textup{contains at least two scale points}\right)
\\\notag=&\mathbb{P}\left(S_{Y}(t,s_{1}) = \underset{s\in B_{m}(t)}{\sup}\ S_{Y}(t,s) = S_{Y}(t,s_{2})\ \textup{for some different}\ s_{1},s_{2}\in B_{m}(t)\right)
\\\notag\leq&  \mathbb{P}\left(\underset{s\in N(s_{1})\cap B_{m}(t)}{\sup}S_{Y}(t,s) =
\underset{s\in N(s_{2})\cap B_{m}(t)}{\sup}S_{Y}(t,s)\ \textup{for some different}\ s_{1},s_{2}\in
\mathbb{Q}\cap B_{m}(t)\right)
\\\notag\leq&  \underset{\begin{subarray}{c} s_{1},s_{2}\in \mathbb{Q}\cap B_{m}(t)\\ s_{1}\neq s_{2}\end{subarray}}{\sum} \mathbb{P}\left(\underset{s\in N(s_{1})\cap B_{m}(t)}{\sup}S_{Y}(t,s) =
\underset{s\in N(s_{2})\cap B_{m}(t)}{\sup}S_{Y}(t,s)\right).
\end{align}
Thus, to prove the second part of Theorem \ref{prop:unique_argmax},
it suffices to show that,
for any different $s_{1},s_{2}\in \mathbb{Q}\cap B_{m}(t)$, there exist neighborhoods $N(s_{1})$ and $N(s_{2})$ of $s_{1}$ and $s_{2}$, respectively, such that
\begin{align}\label{claim:sup_vs_sup:m}
\mathbb{P}\left(\underset{s\in N(s_{1})\cap B_{m}(t)}{\sup}S_{Y}(t,s) =
\underset{s\in N(s_{2})\cap B_{m}(t)}{\sup}S_{Y}(t,s)\right)=0.
\end{align}
The proof of (\ref{claim:sup_vs_sup:m}) is similar to that of (\ref{claim:sup_vs_sup}).
Therefore, we omit the details here.
\qed

\begin{Remark}\label{remark:Leb0}
Suppose that (\ref{lebesgue_zero_intersection}) does not hold. Then, $\textup{Leb}\left(\mathcal{I}\cap D\right)>0$.
By the Lebesgue density theorem, there exists $\mathbf{z}\in \mathcal{I}\cap D$ such that for any
$\varepsilon>0$,
\begin{align*}
\textup{Leb}\left(\mathcal{I}\cap D \cap B_{\varepsilon}(\mathbf{z})\right)>0,
\end{align*}
where $B_{\varepsilon}(\mathbf{z})$ is the disk centered at $\mathbf{z}$ and with radius $\varepsilon$.
For any $\varepsilon>0$ and $\mathbf{w}\in \mathbb{R}^{2}$
satisfying
$\mathbf{z}+\mathbf{w}\in \mathcal{I}\cap B_{\varepsilon}(\mathbf{z})$,
because $M_{1}(\mathbf{z}+\mathbf{w})=M_{2}(\mathbf{z}+\mathbf{w})$,
(\ref{second_derivative_convex}) implies that
\begin{align}\label{second_derivative_convex2}
\left(\mathbf{v}_{1}(\mathbf{z})-\mathbf{v}_{2}(\mathbf{z})\right)^{\top}\mathbf{w}
+\frac{1}{2}\mathbf{w}^{\top}\left(\mathbf{H}_{1}(\mathbf{z})-\mathbf{H}_{2}(\mathbf{z})\right) \mathbf{w}=o(\|\mathbf{w}\|^{2}).
\end{align}
For the left-hand side of (\ref{second_derivative_convex2}),  from (\ref{second_derivative_convex0}),
\begin{align}\notag
&\left(\mathbf{v}_{1}(\mathbf{z})-\mathbf{v}_{2}(\mathbf{z})\right)^{\top}\mathbf{w}
+\frac{1}{2}\mathbf{w}^{\top}\left(\mathbf{H}_{1}(\mathbf{z})-\mathbf{H}_{2}(\mathbf{z})\right) \mathbf{w}
\\\notag\geq&\left(\mathbf{v}_{1}(\mathbf{z})-\mathbf{v}_{2}(\mathbf{z})\right)^{\top}\mathbf{w}
+(c_{1}-c_{2})\|\mathbf{w}\|^{2},
\end{align}
which contradicts the right-hand side of (\ref{second_derivative_convex2}).
\end{Remark}

\subsection{Proof of Theorem \ref{lemma:compact&hemiconti}}\label{sec:proof:lemma:compact&hemiconti}
Let $\Omega$ be the space of random elements.
For any $\omega\in \Omega$, we denote the associated sample path of the random process $\Phi$ by $\{\Phi(t;\omega)\mid t\in \mathbb{R}\}$.
Besides, to emphasize the dependence of the scalogram of $Y$ on the sample path of $\Phi$,
we denote the scalogram by $S_{Y}(\cdot,\cdot;\omega)$ for the noisy signal $Y(t;\omega) = f(t)+\Phi(t;\omega)$.
Proposition \ref{lemma:twice_differentiability}
shows that there exists a subset $\widetilde{\Omega}$ of $\Omega$ with
$\mathbb{P}(\widetilde{\Omega})=1$ such that for any $\omega\in \widetilde{\Omega}$,
$S_{Y}(\cdot,\cdot;\omega)$ is two times continuously differentiable on $\mathbb{R}\times (0,\infty)$.

In the following, we focus on the case where $\omega\in \widetilde{\Omega}$.
In order to show that the set-valued function $t\in \mathbb{R} \mapsto s_{Y}(t;\omega) \subseteq [0,\infty)$
is upper hemicontinuous, we classify time points into two groups:
\begin{align*}
\mathcal{M}:= \left\{t\in \mathbb{R}\mid S_{Y}(t,s;\omega)>0\ \textup{for a certain}\ s>0\right\}
\end{align*}
and
\begin{align*}
\mathcal{N}:= \left\{t\in \mathbb{R}\mid S_{Y}(t,s;\omega)=0\ \textup{for all}\ s>0\right\}.
\end{align*}

For any $t_{0}\in \mathcal{N}$, it is clear that $s_{Y}(t_{0};\omega)=[0,\infty)$.
Thus, $s_{Y}(t';\omega)\subseteq s_{Y}(t_{0};\omega)$ for any time point $t'$.
This directly implies that $s_{Y}(\cdot;\omega)$ is upper hemicontinuous at every time point in $\mathcal{N}$.
It remains to prove that $s_{Y}(\cdot;\omega)$ is also upper hemicontinuous at every $t_{0}\in \mathcal{M}$.
The proof is structured as follows.

{\it Step 1.}
We begin by showing that there exist an interval $[a,b]$, whose interior contains the time point $t_{0}$, and a threshold $\mathcal{S}(\omega)<\infty$ such that
\begin{align}\label{claim:compact_valued}
 s_{Y}(t;\omega) = \arg\underset{s>0}{\max}\ S_{Y}(t,s;\omega)\subseteq [0,\mathcal{S}(\omega)]
\end{align}
for all $t\in [a,b].$
Because $S_{Y}(t_{0},s_{0};\omega)>0$ for some $s_{0}>0$ and $S_{Y}(\cdot,\cdot;\omega)$ is continuous
on $\mathbb{R}\times (0,\infty)$, there exist $\delta>0$ and an interval $[a,b]$ with $a<t_{0}<b$
such that
\begin{equation}\label{mountain}
S_{Y}(\cdot,\cdot;\omega)>\frac{S_{Y}(t_{0},s_{0};\omega)}{2}\ \textup{on}\ [a,b]\times [s_{0}-\delta,s_{0}+\delta].
\end{equation}
We now fix the interval $[a,b]$ and proceed to prove the finiteness of the threshold $\mathcal{S}(\omega)$.
Assume, for the sake of contradiction, that no such threshold exists.
This would imply that for any $n\in \mathbb{N}$, there exist  $t_{n}\in [a,b]$ and
$s_{n}\geq n$ such that  $s_{n}\in  s_{Y}(t_{n};\omega)$. That is,
\begin{align}\label{proxy:s_n}
S_{Y}(t_{n},s_{n};\omega) = \underset{s>0}{\sup}\ S_{Y}(t_{n},s;\omega).
\end{align}
To establish the impossibility of this scenario, we define a function $\mathcal{R}$ on  $[a,b]\times [0,1]$ as follows
\begin{align}\label{reciprocal:R}
\mathcal{R}(t,u;\omega) =\left\{\begin{array}{ll}S_{Y}(t,1/u;\omega)\ &\ \textup{if}\ u\in (0,1],
\\
0\ &\ \textup{if}\ u=0.
   \end{array}\right.
\end{align}
Because $S_{Y}$ is continuous on $\mathbb{R}\times (0,\infty)$ and $\mathbb{P}(\underset{s\rightarrow \infty}{\lim} S_{Y}(t,s)=0)=1$ for any $t\in \mathbb{R}$, as shown in Corollary \ref{lemma:continuity_extension},
the function $\mathcal{R}$ is continuous on $[a,b]\times [0,1]$.
Because $\{t_{n}\}_{n\in \mathbb{N}}\subset [a,b]$,
there exists a convergent subsequence of  $\{t_{n}\}_{n\in \mathbb{N}}$ which converges to some time point $t^{*}\in [a,b]$.
For simplicity of notation, assume that  $t_{n}\rightarrow t^{*}$ as $n\rightarrow \infty$.
Because the pointwise supremum of any collection of continuous functions is lower semi-continuous, we have
\begin{align}\notag
\underset{s>0}{\sup}\ S_{Y}(t^{*},s;\omega) \leq& \underset{n\rightarrow\infty}{\liminf}\ \underset{s>0}{\sup}\ S_{Y}(t_{n},s;\omega)
\\\notag\overset{(\ref{proxy:s_n})}{=}&\underset{n\rightarrow\infty}{\liminf}\ S_{Y}(t_{n},s_{n};\omega)
\\\notag\overset{(\ref{reciprocal:R})}{=}&\underset{n\rightarrow\infty}{\liminf}\ \mathcal{R}(t_{n},1/s_{n};\omega)
\\\notag=&\mathcal{R}(t^{*},0;\omega)=0.
\end{align}
The derivation above  implies that $S_{Y}(t^{*},s;\omega)=0$ for any $s>0$, which contradicts
 (\ref{mountain}). Hence, (\ref{claim:compact_valued}) holds.

{\it Step 2.}
By  (\ref{claim:compact_valued}), for any $t\in[a,b]$,
\begin{align*}
 s_{Y}(t;\omega) = \arg\underset{s>0}{\max}\ S_{Y}(t,s;\omega) = \arg\underset{s\in [0, \mathcal{S}(\omega)]}{\max}\ S_{Y}(t,s; \omega).
\end{align*}
The continuity of $S_{Y}(t,s;\omega)$ and the compactness of  $[0, \mathcal{S}(\omega)]$ allow us to apply
Berge's maximum theorem \cite[Theorem 3.4]{hu1997handbook} (see also \cite[p. 570]{aliprantis2006infinite} or
Lemma \ref{lemma:Berge} in the appendix), implying that the set-valued function $ s_{Y}(\cdot;\omega)$
is upper hemicontinuous at every time point $t_{0}\in \mathcal{M}$.
\qed

\subsection{Proof of Corollary \ref{prop:C1_argmax}}\label{sec:proof:prop:C1_argmax}
Using the notation introduced at the beginning of Section \ref{sec:proof:lemma:compact&hemiconti},
Theorems \ref{prop:unique_argmax} and \ref{lemma:compact&hemiconti}
show that there exists a subset $\widetilde{\Omega}$ of $\Omega$ with
$\mathbb{P}(\widetilde{\Omega})=1$ such that for any $\omega\in \widetilde{\Omega}$,
\begin{itemize}
\item $S_{Y}(\cdot,\cdot;\omega)$ is two times continuously differentiable on $\mathbb{R}\times (0,\infty)$,
\item $ s_{Y}(t;\omega)$ is a singleton for every $t\in \mathbb{Q}$, and
\item $t\mapsto  s_{Y}(t;\omega)$ is upper hemicontinuous on $\mathbb{R}$.
\end{itemize}
In the following, we focus on the case where $\omega\in \widetilde{\Omega}$
and consider an arbitrary point $t_{0}\in \mathbb{Q}$.
The proof of Corollary \ref{prop:C1_argmax} is organized as follows.

{\it Step 1.}
According to the definition of upper hemicontinuity in Theorem \ref{lemma:compact&hemiconti},
for any $\varepsilon>0$, there exists $\delta>0$ such that
\begin{align}\label{upper:inclusion_relation}
 s_{Y}(t';\omega)\subset ( s_{Y}(t_{0};\omega)-\varepsilon, s_{Y}(t_{0};\omega)+\varepsilon)
\end{align}
for all $t'\in (t_{0}-\delta,t_{0}+\delta)$.
The inclusion relation (\ref{upper:inclusion_relation}) implies the continuity
of the restriction of $ s_{Y}$ on $\mathbb{Q}$.

{\it Step 2:} We aim to demonstrate the existence of a continuously differentiable function $r$ defined on $(t_{0}-\delta,t_{0}+\delta)$
for some $\delta>0$ such that $r(t_{0}) =  s_{Y}(t_{0};\omega)$ and
\begin{align}\label{IFT1s}
\frac{\partial  S_{Y}}{\partial s}(t,r(t))=0
\end{align}
for any $t\in (t_{0}-\delta,t_{0}+\delta)$.
Because the function $s\mapsto S_{Y}(t,s;\omega)$ reaches its global maximum at only
$ s_{Y}(t_{0};\omega)$,
\begin{align}\label{IFT1}
\frac{\partial  S_{Y}}{\partial s}(t_{0}, s_{Y}(t_{0});\omega) = 0.
\end{align}
On the other hand, if the assumption (\ref{assumption:atomless}) holds, then
\begin{align}\label{IFT_add_assumption}
\frac{\partial^{2}  S_{Y}}{\partial s^{2}}(t_{0}, s_{Y}(t_{0});\omega)\neq 0
\end{align}
for almost every $\omega\in \widetilde{\Omega}$.
By (\ref{IFT1}) and (\ref{IFT_add_assumption}),  the implicit function theorem
guarantees the existence of the function $r$ and the positive constant $\delta$.

{\it Step 3:} We want to show that the function $r$ determined by the implicit function theorem
coincides with $ s_{Y}(\cdot;\omega)$ on rational numbers in a neighborhood of $t_{0}$.
That is, there exists a constant $\delta'>0$ with  $\delta'\leq\delta$ such that
\begin{align}\label{IFT:g=s:Q}
r(t) = s_{Y}(t;\omega)\ \textup{for any}\ t\in (t_{0}-\delta',t_{0}+\delta')\cap \mathbb{Q}.
\end{align}
Suppose not. Then, (a) there exists a sequence $\{q_{n}\}_{n\in \mathbb{N}}\subset (t_{0}-\delta',t_{0}+\delta')\cap \mathbb{Q}\setminus\{t_{0}\}$
such that $q_{n}\rightarrow t_{0}$ as $n\rightarrow\infty$,
(b) $r(q_{n}) \neq s_{Y}(q_{n};\omega)$ for any $n\in \mathbb{N}$.
For the condition (b), note that if (\ref{IFT:g=s:Q}) does not hold for any $\delta'$, then
there exists $q_{1}\in (t_{0}-\delta',t_{0}+\delta')\cap \mathbb{Q}\setminus\{t_{0}\}$ such that
$r(q_{1}) \neq s_{Y}(q_{1};\omega)$. We can iteratively choose $q_{n}\in (t_{0}-\delta',t_{0}+\delta')\cap \mathbb{Q}\setminus\{t_{0}\}$  from the intervals
$(q_{n-1},t_{0})$ or $(t_{0},q_{n-1})$ for $n>1$ such that $r(q_{n}) \neq s_{Y}(q_{n};\omega)$.
Because
\begin{itemize}
\item the continuity of $r$ implies that $r(q_{n})\rightarrow r(t_{0})$ as $n\rightarrow\infty$,
\item the continuity
of $ s_{Y}(\cdot;\omega)$ on $\mathbb{Q}$ implies that $ s_{Y}(q_{n};\omega)\rightarrow  s_{Y}(t_{0};\omega)$ as $n\rightarrow\infty$, and
\item $ s_{Y}(t_{0};\omega) = r(t_{0})$,
\end{itemize}
we obtain that $r(q_{n})- s_{Y}(q_{n};\omega)\rightarrow0$ and
\begin{align}\label{difference_quotient:partialS}
\left[\frac{\partial S_{Y}}{\partial s}(q_{n},r(q_{n});\omega)
-\frac{\partial S_{Y}}{\partial s}(q_{n}, s_{Y}(q_{n});\omega)\right]\left[r(q_{n})- s_{Y}(q_{n})\right]^{-1}\rightarrow
\frac{\partial^{2} S_{Y}}{\partial s^{2}}(t_{0}, s_{Y}(t_{0});\omega)
\end{align}
as $n\rightarrow\infty$. Here, we used the two times continuous differentiability of $S_{Y}(\cdot,\cdot;\omega)$ to obtain the limit in (\ref{difference_quotient:partialS}).
However, (\ref{IFT1s}) and (\ref{IFT1}) imply that the left-hand side of (\ref{difference_quotient:partialS}) is equal to zero.
This leads to
\begin{equation*}
\frac{\partial^{2} S_{Y}}{\partial s^{2}}(t_{0}, s_{Y}(t_{0});\omega)=0,
\end{equation*}
which contradicts (\ref{IFT_add_assumption}). Therefore, (\ref{IFT:g=s:Q}) holds.

{\it Step 4:} We finally demonstrate that
\begin{equation}\label{g_Qc:max}
r(t)\in  s_{Y}(t;\omega)\ \textup{for}\ t\in (t_{0}-\delta',t_{0}+\delta')\cap \mathbb{Q}^{c}.
\end{equation}
Indeed, for any $t\in (t_{0}-\delta',t_{0}+\delta')\cap \mathbb{Q}^{c}$,
there exists a sequence $\{q_{n}\}_{n\in \mathbb{N}}\subset (t_{0}-\delta',t_{0}+\delta')\cap \mathbb{Q}$
such that $q_{n}\rightarrow t$ as $n\rightarrow\infty$. For any $s>0$,
\begin{align}\notag
S_{Y}(t,r(t);\omega) =& \underset{n\rightarrow\infty}{\lim} S_{Y}(q_{n},r(q_{n});\omega)
\\\notag\geq&  \underset{n\rightarrow\infty}{\lim} S_{Y}(q_{n},s;\omega)
\\\label{check:maximizerproperty}=& S_{Y}(t,s;\omega),
\end{align}
where the equalities follow from the continuity of
$r$ on the interval $(t_{0}-\delta',t_{0}+\delta')$ and the continuity of
$S_{Y}(\cdot,\cdot;\omega)$.
The observation (\ref{check:maximizerproperty}) implies that
\begin{align*}
S_{Y}(t,r(t);\omega)\geq \underset{s>0}{\sup}\ S_{Y}(t,s;\omega).
\end{align*}
Hence, (\ref{g_Qc:max}) holds.
\qed

\subsection{Proof of Theorem \ref{lemma:conditional}}\label{sec:proof:lemma:conditional}
By the definition of $ s_{Y}(t)$ in (\ref{proxy:s_Y}),
for any interval $I$ containing $s_{f}(t)$,
we have
\begin{align}\notag
&\mathbb{P}\left( s_{Y}(t)\in I\right)
=\mathbb{P}\left(\underset{s\in I}{\max}|W_{Y}(t,s)| >\underset{s\in I^{c}}{\max}|W_{Y}(t,s)|\right)
\\\notag\geq& \mathbb{P}\left(|W_{f}(t,s_{f}(t))|-|W_{\Phi}(t,s_{f}(t))| >\underset{s\in I^{c}}{\max}|W_{f}(t,s)|+\underset{s\in I^{c}}{\max}|W_{\Phi}(t,s)|\right)
\\\notag=&\mathbb{P}\left(|W_{f}(t,s_{f}(t))|-\underset{s\in I^{c}}{\max}|W_{f}(t,s)|>|W_{\Phi}(t,s_{f}(t))| +\underset{s\in I^{c}}{\max}|W_{\Phi}(t,s)|\right)
\\\notag=&\mathbb{P}\left(\triangle_{I}>|W_{\Phi}(t,s_{f}(t))| +\underset{s\in I^{c}}{\max}|W_{\Phi}(t,s)|\right),
 \end{align}
where $\triangle_{I}= |W_{f}(t,s_{f}(t))|-\underset{s\in I^{c}}{\max}|W_{f}(t,s)|$.
Hence, to prove (\ref{ineq:sYinI}), it suffices to show that
\begin{align}\label{complex_Borell_ineq}
\mathbb{P}\left(|W_{\Phi}(t,s_{f}(t))| +\underset{s\in I^{c}}{\max}|W_{\Phi}(t,s)|\geq \triangle_{I}\right)
\leq
\textup{exp}\left(-\frac{(\triangle_{I}-\mu_{I})^{2}}{\sigma_{I}^{2}}\right),
\end{align}
where $\mu_{I}$ and $\sigma_{I}$ are defined in (\ref{def:mu_ab}) and (\ref{def:sigma_ab}).
The proof of (\ref{complex_Borell_ineq}) consists of two steps.

{\it Step 1: Establish a discrete version of (\ref{complex_Borell_ineq}).}
For any $d\in \mathbb{N}$ and
$\mathbf{x}=[x_{1}\ x_{2}\ \cdots\ x_{d}]^{\top}\in \mathbb{C}^{d}$,
we define
$$g(\mathbf{x}) = |x_{1}|+\underset{2\leq \ell \leq d}{\max}|x_{\ell}|,$$
which is Lipschitz continuous.
For any $s_{1},s_{2},\ldots,s_{d}>0$, let
\begin{align*}
\mathbf{J}= \left[W_{\Phi}(t,s_{1})\ W_{\Phi}(t,s_{2})\ \cdots\ W_{\Phi}(t,s_{d})\right]^{\top}\in
\mathbb{C}^{d}.
\end{align*}
Because the wavelet $\psi$ is assumed to be analytic, $\mathbf{J}$ is a circularly symmetric complex Gaussian random vector, as shown in Proposition \ref{lemma:exp}(a).
By Lemma \ref{lemma:Lip},
for any $u>\mathbb{E}[g(\mathbf{J})]$,
\begin{align}\label{ineq:|W|+max|W|_v1}
\mathbb{P}\left(|W_{\Phi}(t,s_{1})|+\underset{2\leq \ell\leq d}{\max}|W_{\Phi}(t,s_{\ell})|>u\right)
= \mathbb{P}\left(g(\mathbf{J})>u\right)\leq  \textup{exp}\left(-\frac{(u-\mathbb{E}[g(\mathbf{J})])^{2}}{2\|g(\mathbf{A}\bullet)\|_{\textup{Lip}}^{2}}\right),
\end{align}
where $\mathbf{A}\in \mathbb{C}^{d\times d}$ satisfying
\begin{align}\label{EJJ=2AA}
\mathbf{A}\mathbf{A}^{*}= 2^{-1}\mathbb{E}\left[\mathbf{J}\mathbf{J}^{*}\right].
\end{align}
About $\|g(\mathbf{A}\bullet)\|_{\textup{Lip}}$,
for any $\mathbf{x},\mathbf{y}\in \mathbb{C}^{d}$ with $\mathbf{z} = \mathbf{x}-\mathbf{y}$,
\begin{align}\notag
\left|g(\mathbf{A}\mathbf{x})-g(\mathbf{A}\mathbf{y})\right|
\leq&
\left|e_{1}^{\top}\mathbf{A}\mathbf{z}\right|
+\underset{2\leq \ell \leq d}{\max}\left|e_{\ell}^{\top}\mathbf{A}\mathbf{z}\right|
\\\notag\leq&
|e_{1}^{\top}\mathbf{A}||\mathbf{z}|
+\underset{2\leq \ell \leq d}{\max}|e_{\ell}^{\top}\mathbf{A}||\mathbf{z}|
\\\label{comput_Lip_1}=&
\left(\sqrt{e_{1}^{\top}\mathbf{A}\mathbf{A}^{*}e_{1}}
+\underset{2\leq \ell \leq d}{\max}
\sqrt{e_{\ell}^{\top}\mathbf{A}\mathbf{A}^{*}e_{\ell}}\right)
|\mathbf{z}|,
\end{align}
where $e_{\ell}$ is the column vector in $\mathbb{R}^{d}$ with one in position $\ell$ and zeros elsewhere.
From (\ref{EJJ=2AA}), for $\ell\in\{1,2,\ldots,d\}$,
\begin{align}\label{comput_Lip_2}
e_{\ell}^{\top}\mathbf{A}\mathbf{A}^{*}e_{\ell}
=2^{-1} \mathbb{E}\left[|W_{\Phi}(t,s_{\ell})|^{2}\right].
\end{align}
Combining (\ref{comput_Lip_1}) and (\ref{comput_Lip_2}) yields
\begin{align}\label{upperbound_Lip}
\|g(\mathbf{A}\bullet)\|_{\textup{Lip}}\leq 2^{-\frac{1}{2}}\left\{\sqrt{\mathbb{E}\left[|W_{\Phi}(t,s_{1})|^{2}\right]}
+\underset{2\leq \ell \leq d}{\max}
\sqrt{\mathbb{E}\left[|W_{\Phi}(t,s_{\ell})|^{2}\right]}\right\}.
\end{align}

{\it Step 2: Extend the discrete version to the continuous case.}
By setting $s_{1}=s_{f}(t)$ and defining $T_{d}=\{s_{2},\ldots,s_{d}\}$ as a subset of $I^{c}$ such that
$T_{d}\subset T_{d+1}$ and $T_{d}$ increases to a dense subset of $I^{c}$, the inequality (\ref{ineq:|W|+max|W|_v1}) becomes
\begin{align*}
\mathbb{P}\left(|W_{\Phi}(t,s_{f}(t))|+\underset{s\in T_{d}}{\max}|W_{\Phi}(t,s)|>u\right)
\leq  \textup{exp}\left(-\frac{(u-\mathbb{E}[g(\mathbf{J})])^{2}}{2\|g(\mathbf{A}\bullet)\|_{\textup{Lip}}^{2}}\right).
\end{align*}
As $d$ increases to infinity, the continuity of the sample paths of $W_{\Phi}$ given in Proposition \ref{lemma:twice_differentiability},
implies that
\begin{align*}
|W_{\Phi}(t,s_{f}(t))|+\underset{s\in T_{d}}{\max}|W_{\Phi}(t,s)|\uparrow
|W_{\Phi}(t,s_{f}(t))|+\underset{s\in I^{c}}{\max}|W_{\Phi}(t,s)|
\end{align*}
almost surely, the upper bound of $\|g(\mathbf{A}\bullet)\|_{\textup{Lip}}$ in (\ref{upperbound_Lip}) increases to
\begin{align*}
 2^{-\frac{1}{2}}\left\{\sqrt{\mathbb{E}\left[|W_{\Phi}(t,s_{f}(t))|^{2}\right]}
+\underset{s\in I^{c}}{\max}
\sqrt{\mathbb{E}\left[|W_{\Phi}(t,s)|^{2}\right]}\right\}
\end{align*}
by the continuity of $s\mapsto \mathbb{E}\left[|W_{\Phi}(t,s)|^{2}\right]$, and
\begin{align}\label{finite_infinite_mean}
\mathbb{E}[g(\mathbf{J})] \uparrow \mathbb{E}\left[|W_{\Phi}(t,s_{f}(t))|\right]+\mathbb{E}\left[\underset{s\in I^{c}}{\max}|W_{\Phi}(t,s)|\right],
\end{align}
where the finiteness of the expectation in (\ref{finite_infinite_mean}) is guaranteed by Lemma \ref{lemma:Dudley}.
\qed

\subsection{Proof of Theorem \ref{mainresult:deviation}}\label{sec:proof:mainresult:deviation}

{\it Step 1.} For any $m\in\{1,2,\ldots,M\}$, we condition on the event $s_{Y,m}(t)\in B^{\circ}_{m}(t)\subset (\underline{\mathfrak{i}}_{m}(t),\overline{\mathfrak{i}}_{m}(t))$.
Under this condition, we derive the following formula:
\begin{align}\notag
 s_{Y,m}(t)-s_{f,m}(t) =&
-
\left[\frac{\partial^{2}S_{f}}{\partial s^{2}}(t,c(t))\right]^{-1}
\frac{\partial}{\partial s}\left[W_{f}(t,s)
\overline{W_{\Phi}(t,s)}+\overline{W_{f}(t,s)}
W_{\Phi}(t,s)\right]\Big|_{s= s_{Y,m}(t)}
\\\label{def:ridge_deviation}&-\left[\frac{\partial^{2}S_{f}}{\partial s^{2}}(t,c(t))\right]^{-1}
\frac{\partial S_{\Phi}}{\partial s}(t, s_{Y,m}(t)),
\end{align}
where $c(t)$ lies between $s_{f,m}(t)$ and $ s_{Y,m}(t)$.

Indeed, from the following relationship among the scalograms of the clean signal $f$, the noisy signal $Y$, and the noise $\Phi$:
\begin{align}\notag
S_{Y}(t,s) =& |W_{f}(t,s)+W_{\Phi}(t,s)|^{2}
\\\notag=& S_{f}(t,s)+
W_{f}(t,s)
\overline{W_{\Phi}(t,s)}+\overline{W_{f}(t,s)}
W_{\Phi}(t,s)
+S_{\Phi}(t,s),
\end{align}
we have
\begin{align}\label{label:copy3:scalogram}
\frac{\partial S_{Y}}{\partial s}(t,s)
=  \frac{\partial S_{f}}{\partial s}(t,s)+\frac{\partial}{\partial s}\left[W_{f}(t,s)
\overline{W_{\Phi}(t,s)}+\overline{W_{f}(t,s)}
W_{\Phi}(t,s)\right]
+\frac{\partial S_{\Phi}}{\partial s}(t,s).
\end{align}
Because Theorem \ref{prop:unique_argmax} shows that
$s_{Y,m}(t)$ is a singleton almost surely,
the condition $s_{Y,m}(t)\in B^{\circ}_{m}(t)$ implies that
$$
\frac{\partial S_{Y}}{\partial s}(t, s_{Y,m}(t))=0.
$$
By substituting $s= s_{Y,m}(t)$ into (\ref{label:copy3:scalogram}),
\begin{align}\label{label:copy3:scalogram_0=}
0
=  \frac{\partial S_{f}}{\partial s}(t,s_{Y,m}(t))+\frac{\partial}{\partial s}\left[W_{f}(t,s)
\overline{W_{\Phi}(t,s)}+\overline{W_{f}(t,s)}
W_{\Phi}(t,s)\right]
+\frac{\partial S_{\Phi}}{\partial s}(t,s)\Big|_{s= s_{Y,m}(t)}.
\end{align}
From the condition (R1) in Definition \ref{definition:ridgepoint},
we have
\begin{align*}
\frac{\partial S_{f}}{\partial s}(t,s_{f,m}(t))=0.
\end{align*}
By applying the mean value theorem to the first term on the right-hand side of (\ref{label:copy3:scalogram}),
there exists a random variable $c(t)$ between $s_{f,m}(t)$ and $ s_{Y,m}(t)$ such that
\begin{align}\notag
&\frac{\partial^{2}S_{f}}{\partial s^{2}}(t,c(t))\left( s_{Y,m}(t)-s_{f,m}(t)\right)
\\\label{label:copy4:scalogram}=&-
\frac{\partial}{\partial s}\left[W_{f}(t,s)
\overline{W_{\Phi}(t,s)}+\overline{W_{f}(t,s)}
W_{\Phi}(t,s)\right]\Big|_{s= s_{Y,m}(t)}
-\frac{\partial S_{\Phi}}{\partial s}(t, s_{Y,m}(t)).
\end{align}
Since the condition $B^{\circ}_{m}(t)\subset (\underline{\mathfrak{i}}_{m}(t),\overline{\mathfrak{i}}_{m}(t))$
 implies that
$$\frac{\partial^{2}S_{f}}{\partial s^{2}}(t,c(t))<0,$$
we obtain (\ref{def:ridge_deviation}) from (\ref{label:copy4:scalogram}).

{\it Step 2.} We derive an upper bound for $|s_{Y,m}(t)-s_{f,m}(t)|$.
Let $L_{m}$ be the constant defined in (\ref{def:lower_bound_L}).
On the event $\{s_{Y}(t)\in B^{\circ}_{m}(t)\}$,
$$
\left|\frac{\partial^{2}S_{f}}{\partial s^{2}}(t,c(t))\right|\geq L_{m}.
$$
Hence, (\ref{def:ridge_deviation}) implies that
\begin{align}\label{estimate:deviation1}
&L_{m}| s_{Y,m}(t)-s_{f,m}(t)|
\\\notag\leq &
\left|\frac{\partial}{\partial s}\left[W_{f}(t,s)
\overline{W_{\Phi}(t,s)}+\overline{W_{f}(t,s)}
W_{\Phi}(t,s)\right]\Big|_{s= s_{Y,m}(t)}\right|
+\left|\frac{\partial S_{\Phi}}{\partial s}(t, s_{Y,m}(t))\right|
\\\notag\leq& \left|\frac{\partial W_{f}}{\partial s}(t, s_{Y,m}(t))\right|\left|W_{\Phi}(t, s_{Y,m}(t))\right|
+\left|W_{f}(t, s_{Y,m}(t))\right|\left|\frac{\partial W_{\Phi}}{\partial s}(t, s_{Y,m}(t))\right|
\\\notag&+
\left|\frac{\partial W_{f}}{\partial s}(t, s_{Y,m})\right|\left|W_{\Phi}(t, s_{Y,m}(t))\right|
+\left|W_{f}(t, s_{Y,m}(t))\right|\left|\frac{\partial W_{\Phi}}{\partial s}(t, s_{Y,m}(t))\right|
+\left|\frac{\partial S_{\Phi}}{\partial s}(t, s_{Y,m}(t))\right|
\\\notag=&2\left|\frac{\partial W_{f}}{\partial s}(t, s_{Y,m}(t))\right|
\left|W_{\Phi}(t, s_{Y,m}(t))\right|
+2\left|W_{f}(t, s_{Y,m}(t))\right|\left|\frac{\partial W_{\Phi}}{\partial s}(t, s_{Y,m}(t))\right|
+\left|\frac{\partial S_{\Phi}}{\partial s}(t, s_{Y,m}(t))\right|.
\end{align}
Under the condition $B_{m}(t)\subset (\underline{\mathfrak{i}}_{m}(t),\overline{\mathfrak{i}}_{m}(t))$,
$L_{m}>0$.
By (\ref{estimate:deviation1}), for any $\varepsilon>0$,
\begin{align}\notag
&\mathbb{P}\left(| s_{Y,m}(t)-s_{f,m}(t)|>\varepsilon\ \textup{and}\  s_{Y,m}(t)\in B^{\circ}_{m}(t) \right)
\\\notag\leq
& \mathbb{P}\left(\left|\frac{\partial W_{f}}{\partial s}(t, s_{Y,m}(t))
W_{\Phi}(t, s_{Y,m}(t))\right|>
\frac{\varepsilon L_{m}}{6} \ \textup{and}\  s_{Y,m}\in B^{\circ}_{m}(t) \right)
\\\notag&+\mathbb{P}\left(\left|W_{f}(t, s_{Y,m})\frac{\partial W_{\Phi}}{\partial s}(t, s_{Y,m})\right|>\frac{\varepsilon L_{m}}{6}\ \textup{and}\  s_{Y,m}(t)\in B^{\circ}_{m}(t)  \right)
\\\label{estimate:deviation3}&+\mathbb{P}\left(\left|\frac{\partial S_{\Phi}}{\partial s}(t, s_{Y,m})\right|>\frac{\varepsilon L_{m}}{3}\ \textup{and}\  s_{Y,m}(t)\in B^{\circ}_{m}(t)  \right).
\end{align}

For the first term on the right-hand side of  (\ref{estimate:deviation3}),
by the first part of Lemma \ref{lemma:complex-Borell},
we have
\begin{align}\notag
&\mathbb{P}\left(\left|\frac{\partial W_{f}}{\partial s}(t, s_{Y,m}(t))W_{\Phi}(t, s_{Y,m}(t))\right|>
\frac{\varepsilon L_{m}}{6}\ \textup{and}\  s_{Y,m}(t)\in B^{\circ}_{m}(t)\right)
\\\notag\leq& \mathbb{P}\left(\underset{s\in B_{m}(t)}{\max}
\left|\frac{\partial W_{f}}{\partial s}(t,s)\right|\left|W_{\Phi}(t,s)\right|>
\frac{\varepsilon L_{m}}{6}\right)
\\\label{estimate:deviation4}\leq&
\textup{exp}\left[-\frac{1}{\sigma_{m,1}^{2}}\left(\frac{\varepsilon L_{m}}{6}-\mu_{m,1}\right)^{2}\right]
\end{align}
for any $\varepsilon>6\mu_{m,1}L_{m}^{-1}$, where $\mu_{m,1}$ and $\sigma_{m,1}^{2}$ are defined in (\ref{def:mu1_4}).

For the second term on the right-hand side of  (\ref{estimate:deviation3}), by the second part of Lemma \ref{lemma:complex-Borell},
we have
\begin{align}\notag
&\mathbb{P}\left(\left|W_{f}(t, s_{Y,m}(t))\frac{\partial W_{\Phi}}{\partial s}(t, s_{Y,m}(t))\right|>\frac{\varepsilon L_{m}}{6}
\ \textup{and}\  s_{Y,m}(t)\in B^{\circ}_{m}(t)\right)
\\\notag\leq&\mathbb{P}\left(\underset{s\in B_{m}(t)}{\max}
\left|W_{f}(t,s)\right|\left|\frac{\partial W_{\Phi}}{\partial s}(t,s)\right|>\frac{\varepsilon L_{m}}{6}
\right)
\\\label{estimate:deviation5}\leq&
\textup{exp}\left[-\frac{1}{\sigma_{m,2}^{2}}\left(\frac{\varepsilon L_{m}}{6}-\mu_{m,2}\right)^{2}\right]
\end{align}
for any $\varepsilon>6\mu_{2}L_{m}^{-1}$, where $\mu_{m,2}$ and $\sigma_{m,2}^{2}$ are defined in (\ref{def:mu1_4}).

For the third term in (\ref{estimate:deviation3}), because
\begin{align}\notag
\left|\frac{\partial S_{\Phi}}{\partial s}(t,s)\right|
=&\left|\frac{\partial}{\partial s}\left(W_{\Phi}(t,s)\overline{W_{\Phi}(t,s)}\right)\right|
\\\notag=&
\left|\left(\frac{\partial}{\partial s}W_{\Phi}(t,s)\right)\overline{W_{\Phi}(t,s)}
+W_{\Phi}(t,s)\left(\frac{\partial}{\partial s}\overline{W_{\Phi}(t,s)}\right)\right|
\\\notag\leq&2\left|W_{\Phi}(t,s)\right|
\left|\frac{\partial W_{\Phi}}{\partial s}(t,s)\right|,
\end{align}
we have
\begin{align}\notag
&\mathbb{P}\left(\left|\frac{\partial S_{\Phi}}{\partial s}(t, s_{Y,m}(t))\right|>\frac{\varepsilon L_{m}}{3}
\ \textup{and}\  s_{Y,m}(t)\in B^{\circ}_{m}(t)\right)
\\\notag\leq &
\mathbb{P}\left(\left|W_{\Phi}(t, s_{Y,m}(t))\right|>\sqrt{\frac{\varepsilon L_{m}}{6}}
\ \textup{and}\  s_{Y,m}(t)\in B^{\circ}_{m}(t)\right)
\\\notag&+\mathbb{P}\left(\left|\frac{\partial W_{\Phi}}{\partial s}(t, s_{Y,m}(t))\right|>\sqrt{\frac{\varepsilon L_{m}}{6}}
\ \textup{and}\  s_{Y,m}(t)\in B^{\circ}_{m}(t)\right).
\end{align}
By the first part of Lemma \ref{lemma:complex-Borell},
\begin{align}\notag
&\mathbb{P}\left(\left|W_{\Phi}(t, s_{Y,m}(t))\right|>\sqrt{\frac{\varepsilon L_{m}}{6}}\ \textup{and}\  s_{Y,m}(t)\in B^{\circ}_{m}(t)\right)
\\\label{estimate:deviation6_a}\leq&
\mathbb{P}\left(\underset{s\in B_{m}(t)}{\max}\left|W_{\Phi}(t,s)\right|>\sqrt{\frac{\varepsilon L_{m}}{6}}\right)
\leq  \textup{exp}\left[-\frac{1}{\sigma_{m,3}^{2}}
\left(\sqrt{\frac{\varepsilon L_{m}}{6}}-\mu_{m,3}\right)^{2}\right]
\end{align}
if
$\sqrt{\varepsilon L_{m}}>\mu_{m,3}/\sqrt{6}$,
where $\mu_{m,3}$ and $\sigma_{m,3}^{2}$ are defined in (\ref{def:mu1_4}).
Similarly, by the second part of Lemma \ref{lemma:complex-Borell},
\begin{align}\notag
&\mathbb{P}\left(\left|\frac{\partial W_{\Phi}}{\partial s}(t, s_{Y,m}(t))\right|>\sqrt{\frac{\varepsilon L_{m}}{6}}
\ \textup{and}\  s_{Y,m}\in B^{\circ}_{m}(t)\right)
\\\label{estimate:deviation6_b}\leq&
\mathbb{P}\left(\underset{s\in B_{m}(t)}{\max}\left|\frac{\partial W_{\Phi}}{\partial s}(t,s)\right|>\sqrt{\frac{\varepsilon L_{m}}{6}}\right)
\leq
\textup{exp}\left[-\frac{1}{\sigma_{m,4}^{2}}\left(\sqrt{\frac{\varepsilon L_{m}}{6}}-\mu_{m,4}\right)^{2}\right]
\end{align}
if $\sqrt{\varepsilon L_{m}}>\mu_{m,4}/\sqrt{6}$,
where $\mu_{4}$ and $\sigma_{4}^{2}$ are defined in (\ref{def:mu1_4}).
The proof of Theorem \ref{mainresult:deviation} is completed by summing (\ref{estimate:deviation4}), (\ref{estimate:deviation5}), (\ref{estimate:deviation6_a}), and (\ref{estimate:deviation6_b}).
\qed

\subsection{Proof of Lemma \ref{lemma:Lip}}\label{sec:proof:lemma:Lip}

Because the complex Gaussian random vector $\mathbf{J}$ is assumed to be circularly symmetric,
there exists a matrix $\mathbf{A}\in \mathbb{C}^{d\times d}$ such that
\begin{align}\label{J=AI_c}
\mathbf{J} \overset{\textup{Law}}{=} \mathbf{A}\mathbf{N},
\end{align}
where
$\mathbf{N}=[N_{1}+iN_{d+1}\ N_{2}+iN_{d+2}\ \cdots\ N_{d}+iN_{2d}]^{\top}\in
\mathbb{C}^{d}$
and $\{N_{1},N_{2},\ldots, N_{2d}\}$ is a set of independent standard normal random variables.
The matrix $\mathbf{A}$ in (\ref{J=AI_c}) satisfies
\begin{align*}
\mathbb{E}\left[\mathbf{J}\mathbf{J}^{*}\right] =
 \mathbf{A}\mathbb{E}\left[\mathbf{N}\mathbf{N}^{*} \right]\mathbf{A}^{*}
 =  2\mathbf{A}\mathbf{A}^{*}.
\end{align*}
By the following vector augmentation:
\begin{align*}
\mathbf{N}_{2d}= \begin{bmatrix}\textup{Re}\left(\mathbf{N}\right)
\\ \textup{Im}\left(\mathbf{N}\right) \end{bmatrix}\in \mathbb{R}^{2d},\
\mathbf{J}_{2d}= \begin{bmatrix}\textup{Re}\left(\mathbf{J}\right)
\\ \textup{Im}\left(\mathbf{J}\right) \end{bmatrix}\in \mathbb{R}^{2d}
\ \textup{and}\ \mathbf{A}_{2d}= \left[\begin{array}{lr}\textup{Re}\left(\mathbf{A}\right) & -\textup{Im}\left(\mathbf{A}\right)
\\ \textup{Im}\left(\mathbf{A}\right) & \textup{Re}\left(\mathbf{A}\right)\end{array}\right]\in \mathbb{R}^{2d\times 2d},
\end{align*}
we have $\mathbf{J}_{2d} \overset{\textup{Law}}{=} \mathbf{A}_{2d}\mathbf{N}_{2d}$ from (\ref{J=AI_c}),
and for any $u\in \mathbb{R}$,
\begin{align}\label{cited_ineq0}
\mathbb{P}\left(g(\mathbf{J})>u\right) =
\mathbb{P}\left(g_{2d}(\mathbf{J}_{2d})>u\right)
=\mathbb{P}\left(g_{2d}(\mathbf{A}_{2d}\mathbf{N}_{2d})>u\right),
\end{align}
where $g_{2d}: \mathbb{R}^{2d}\mapsto \mathbb{R}$ is defined by
\begin{align*}
g_{2d}\left(\begin{bmatrix}\textup{Re}(\mathbf{w}) \\ \textup{Im}(\mathbf{w})\end{bmatrix}\right)
= g(\mathbf{w}),\ \mathbf{w}\in \mathbb{C}^{d}.
\end{align*}
By the Gaussian concentration inequality \cite[Corollary 8.5]{kuhn2023maximal}, for any $u>\mathbb{E}[g_{2d}(\mathbf{A}_{2d}\mathbf{N}_{2d})]$,
\begin{align}\label{cited_ineq1}
\mathbb{P}\left(g_{2d}(\mathbf{A}_{2d}\mathbf{N}_{2d})>u\right)
\leq \textup{exp}\left(-\frac{(u-\mathbb{E}[g_{2d}(\mathbf{A}_{2d}\mathbf{N}_{2d})])^{2}}{2\|g_{2d}(\mathbf{A}_{2d}\bullet)\|_{\textup{Lip}}^{2}}\right),
\end{align}
where
\begin{align*}
\|g_{2d}(\mathbf{A}_{2d}\bullet)\|_{\textup{Lip}}
=
\underset{\begin{subarray}{c}\mathbf{x},\mathbf{y}\in \mathbb{R}^{2d}\\ \mathbf{x}\neq\mathbf{y} \end{subarray}}
{\sup} \frac{|g_{2d}(\mathbf{A}_{2d}\mathbf{x})-g_{2d}(\mathbf{A}_{2d}\mathbf{y})|}{|\mathbf{x}-\mathbf{y}|}=\|g(\mathbf{A}\bullet)\|_{\textup{Lip}}.
\end{align*}
The proof of Lemma \ref{lemma:Lip} is completed by substituting (\ref{cited_ineq1}) into
(\ref{cited_ineq0})
and observing that
$\mathbb{E}[g_{2d}(\mathbf{A}_{2d}\mathbf{N}_{2d})]= \mathbb{E}[g(\mathbf{J})]$.
\qed


\subsection{Proof of Lemma \ref{lemma:d_s1_s2}}\label{sec:proof:lemma:d_s1_s2}

(a) By the spectral representation of $W_{\Phi}(t,s) $ in (\ref{spect:WPhi}) and the orthogonal property (\ref{ortho}),
\begin{align}\label{decay:EWPhi2}
\mathbb{E}[|W_{\Phi}(0,s)|^{2}]=& \int_{0}^{\infty} |\widehat{\psi}(s\lambda)|^{2}p(\lambda)d\lambda
=s^{-1}\int_{0}^{\infty} |\widehat{\psi}(\lambda)|^{2}p\left(s^{-1}\lambda\right)d\lambda.
\end{align}
Under Assumption \ref{assumption:spectral_density}(b1), $p(\lambda) = L(|\lambda|^{-1})|\lambda|^{H-1}$, where $L$ is slowly varying at infinity. According to \cite[Theorem 1.5.3]{bingham1989regular}, for any slowly varying function $L$ and any $\delta>0$, $L(s)s^{-\delta}\rightarrow 0$ as $s\rightarrow\infty$.
Thus, (\ref{decay:EWPhi2}) implies that
\begin{align}\notag
&\underset{s\rightarrow\infty}{\lim} L(s)^{-1}s^{H}\mathbb{E}[|W_{\Phi}(0,s)|^{2}]
\\\notag=& \underset{s\rightarrow\infty}{\lim} \int_{0}^{\infty} |\widehat{\psi}(\lambda)|^{2}|\lambda|^{H-1}
\frac{L(s/|\lambda|)}{L(s)}d\lambda
\\\label{secondmoment_tail_limit}=&\int_{0}^{\infty} |\widehat{\psi}(\lambda)|^{2}|\lambda|^{H-1}d\lambda.
\end{align}
The inequality (\ref{decay_var_X1}) follows from (\ref{secondmoment_tail_limit}).

(b) By (\ref{spect:WPhi}),
\begin{equation*}
W_{\Phi}(0,s_{1})-W_{\Phi}(0,s_{2})
= \int_{0}^{\infty} \left[\ \overline{\widehat{\psi}(s_{1}\lambda)}-\overline{\widehat{\psi}(s_{2}\lambda)}\ \right]\sqrt{p(\lambda)}Z(d\lambda).
\end{equation*}
By the orthogonal property (\ref{ortho}),
\begin{align}\notag
d_{W_{\Phi}}(s_{1},s_{2}) =&\left\{\mathbb{E}\left[|W_{\Phi}(0,s_{1})-W_{\Phi}(0,s_{2})|^{2}\right]\right\}^{\frac{1}{2}}
\\\label{d_s1_s2_v2:X0}=&
\left(\int_{0}^{\infty}
|\widehat{\psi}(s_{1}\lambda)-\widehat{\psi}(s_{2}\lambda)|^{2}p(\lambda)
d\lambda\right)^{\frac{1}{2}}.
\end{align}
Because $\widehat{\psi}$ is assumed to be Lipschitz continuous on $(0,\infty)$
with Lipschitz constant $\|\widehat{\psi}\|_{\textup{Lip}}$,
 (\ref{d_s1_s2_v2:X0}) can be estimated as follows
\begin{align}\notag
d^{2}_{W_{\Phi}}(s_{1},s_{2})  \leq& \int_{0}^{\infty}\left(\min\left\{2\|\widehat{\psi}\|_{\infty},\|\widehat{\psi}\|_{\textup{Lip}}|s_{1}-s_{2}|\lambda\right\}\right)^{2} p(\lambda)d\lambda
\\\label{d_s1_s2_v3:X0}=&4\|\widehat{\psi}\|_{\infty}^{2}\int_{\Lambda}^{\infty}p(\lambda)d\lambda
+
\|\widehat{\psi}\|_{\textup{Lip}}^{2} |s_{1}-s_{2}|^{2}\int_{0}^{\Lambda}\lambda^{2} p(\lambda)d\lambda,
\end{align}
where $\Lambda = 2\|\widehat{\psi}\|_{\infty}\|\widehat{\psi}\|_{\textup{Lip}}^{-1}|s_{1}-s_{2}|^{-1}$.
Under Assumption \ref{assumption:spectral_density}(a),
\begin{align}\label{estimate_d_part1:X0}
\int_{\Lambda}^{\infty}p(\lambda)d\lambda\leq C_{1}\int_{\Lambda}^{\infty}\lambda^{-(1+\gamma)}d\lambda
=\frac{C_{1}}{\gamma}\left(\frac{2\|\widehat{\psi}\|_{\infty}}{\|\widehat{\psi}\|_{\textup{Lip}}|s_{1}-s_{2}|}\right)^{-\gamma}.
\end{align}
Let $\kappa=2\|\widehat{\psi}\|_{\infty}\|\widehat{\psi}\|_{\textup{Lip}}^{-1}$.
When $|s_{1}-s_{2}|\leq 1$, we have $\Lambda\geq \kappa$
and
\begin{align}\notag
\int_{0}^{\Lambda}\lambda^{2} p(\lambda)d\lambda
\leq& \int_{0}^{\kappa} \lambda^{2}p(\lambda)d\lambda +C_{1}\int_{\kappa}^{\Lambda}\lambda^{1-\gamma}d\lambda
\\\label{estimate_d_part2:X0}
\leq&2\frac{\|\widehat{\psi}\|_{\infty}^{2}}{\|\widehat{\psi}\|_{\textup{Lip}}^{2}}\textup{Var}\left(\Phi(0)\right)+
\frac{2^{2-\gamma}C_{1}}{2-\gamma}\left(\frac{\|\widehat{\psi}\|_{\infty}}{\|\widehat{\psi}\|_{\textup{Lip}}}\right)^{2-\gamma}
\left(|s_{1}-s_{2}|^{\gamma-2}
-1\right).
\end{align}
By substituting (\ref{estimate_d_part1:X0}) and (\ref{estimate_d_part2:X0}) into (\ref{d_s1_s2_v3:X0}),
\begin{align}\notag
d^{2}_{W_{\Phi}}(s_{1},s_{2})  \leq&
C_{1}\left(\frac{1}{\gamma}+\frac{1}{2-\gamma}\right)(2\|\widehat{\psi}\|_{\infty})^{2-\gamma}\|\widehat{\psi}\|_{\textup{Lip}}^{\gamma}|s_{1}-s_{2}|^{\gamma}
\\\notag&+
\left(2\|\widehat{\psi}\|_{\infty}^{2}\textup{Var}\left(\Phi(0)\right)
-\frac{C_{1}}{2-\gamma}(2\|\widehat{\psi}\|_{\infty})^{2-\gamma}\|\widehat{\psi}\|_{\textup{Lip}}^{\gamma}\right)|s_{1}-s_{2}|^{2}
\\\notag\leq &
\left[C_{1}\left(\frac{1}{\gamma}+\frac{2}{2-\gamma}\right)(2\|\widehat{\psi}\|_{\infty})^{2-\gamma}\|\widehat{\psi}\|_{\textup{Lip}}^{\gamma}
+2\|\widehat{\psi}\|_{\infty}^{2}\textup{Var}\left(\Phi(0)\right)\right]|s_{1}-s_{2}|^{\gamma}
\end{align}
for $|s_{1}-s_{2}|\leq 1$.
\qed

\subsection{Proof of Lemma \ref{lemma:Dudley}}\label{sec:proof:lemma:Dudley}
The proof idea can be traced back to the chaining argument in \cite{talagrand2014upper}.
We firstly create a sequence of sets, denoted as $\{\pi_{j}\}_{j\in \mathbb{N}\cup\{0\}}$, of finite points in $[0,\infty)$
such that for each $j\in \mathbb{N}\cup\{0\}$,
\begin{align}\label{nest_structure1}
\underset{s> 0}{\sup} \underset{s'\in \pi_{j}}{\inf} d_{W_{\Phi}}(s,s')\leq \mathcal{D}^{1-j},
\end{align}
where $\mathcal{D}$ is a constant greater than two, as defined in (\ref{def:D}).
Because $\underset{s\rightarrow 0+}{\lim}W_{\Phi}(\cdot,s)=0$ almost surely, as shown in Corollary \ref{lemma:continuity_extension}, we define $W_{\Phi}(\cdot,0)=0$.
Hence, for any $s\geq0$,
$$
d_{W_{\Phi}}(0,s) = \left\{\mathbb{E}\left[|W_{\Phi}(0,s)|^{2}\right]\right\}^{\frac{1}{2}}\leq\mathcal{D},
$$
where the inequality follows from the definition of $\mathcal{D}$ in (\ref{def:D}).
Hence, we can choose $\pi_{0} = \{0\}$.
For every $j\in \mathbb{N}$, denote
\begin{align}\label{def:Sj}
S_{j}=\max\{(2^{-1}\mathcal{D}^{1-j})^{-2/H^{-}},T\}.
\end{align}
By Lemma \ref{lemma:d_s1_s2}(a), for any $s\geq S_{j}$,
\begin{align}\label{pij:property1}
d_{W_{\Phi}}(s,S_{j}) \leq \left\{\mathbb{E}\left[|W_{\Phi}(0,s)|^{2}\right]\right\}^{\frac{1}{2}}+\left\{\mathbb{E}\left[|W_{\Phi}(0,S_{j})|^{2}\right]\right\}^{\frac{1}{2}}
\leq  s^{-H^{-}/2}+S_{j}^{-H^{-}/2}\leq \mathcal{D}^{1-j}.
\end{align}
Consider a set of partition points on $[0,S_{j})$ with lag size $2(C_{2}\mathcal{D}^{j-1})^{-2/\gamma}$ as follows
\begin{align}\label{partition_0_S}
\pi_{j,k} = 2k(C_{2}\mathcal{D}^{j-1})^{-\frac{2}{\gamma}},\ k=0,1,2,\ldots,\lfloor S_{j}2^{-1}(C_{2}\mathcal{D}^{j-1})^{\frac{2}{\gamma}}\rfloor,
\end{align}
where $C_{2}$ is a constant defined in (\ref{def:C2}).
For any $s\in [0,S_{j})$, one of the following situations occurs:
$|s-S_{j}|\leq (C_{2}\mathcal{D}^{j-1})^{-2/\gamma}\leq 1$,
or there exists an integer $k\in\{0,1,2,\ldots,\lfloor S_{j}2^{-1}(C_{2}\mathcal{D}^{j-1})^{2/\gamma}\rfloor\}$ such that
$|s-\pi_{j,k}|\leq (C_{2}\mathcal{D}^{j-1})^{-2/\gamma}\leq 1$.
If the first situation occurs, by Lemma \ref{lemma:d_s1_s2}(b),
\begin{align}\label{pij:property2}
d_{W_{\Phi}}(s,S_{j}) \leq C_{2}|s-S_{j}|^{\gamma/2}\leq \mathcal{D}^{1-j}.
\end{align}
If the second situation occurs, by Lemma \ref{lemma:d_s1_s2}(b) again, we have
\begin{align}\label{pij:property3}
d_{W_{\Phi}}(s,\pi_{j,k}) \leq |s-\pi_{j,k}|^{\gamma/2}\leq \mathcal{D}^{1-j}.
\end{align}
Hence, the set
\begin{align*}
\pi_{j} = \{\pi_{j,0},\pi_{j,1},\ldots\}\cup \{S_{j}\}
\end{align*}
fulfills the condition (\ref{nest_structure1}).
Denote the cardinality of  $\pi_{j}$ by $|\pi_{j}|$.
By (\ref{def:Sj}) and  (\ref{partition_0_S}), for $j\in\mathbb{N}$,
\begin{align}\label{def:cardinality}
|\pi_{j}| \leq  C_{2}^{\frac{2}{\gamma}}\left[2^{\frac{2}{H^{-}}}\mathcal{D}^{\left(\frac{2}{\gamma}+\frac{2}{H^{-}}\right)(j-1)}\vee T\mathcal{D}^{\frac{2}{\gamma}(j-1)}\right].
\end{align}

For any $s>0$ and $j\in \mathbb{N}\cup\{0\}$, we choose a proxy for $s$ from the set $\pi_{j}$ as follows
\begin{align*}
\pi_{j}(s)  = \textup{arg} \underset{s'\in \pi_{j}}{\min}|s-s'|.
\end{align*}
According to (\ref{pij:property1}), (\ref{pij:property2}), and  (\ref{pij:property3}),
\begin{align*}
d_{W_{\Phi}}(s,\pi_{j}(s))\leq \mathcal{D}^{1-j}
\end{align*}
Because $\pi_{0}=\{0\}$ and $W_{\Phi}(0,0)=0$ by definition, we can represent $|W_{\Phi}(0,s)|$ via a telescoping sum as follows
\begin{align}\label{telescoping}
|W_{\Phi}(0,s)| = \overset{\infty}{\underset{j=0}{\sum}} |W_{\Phi}(0,\pi_{j+1}(s))|-|W_{\Phi}(0,\pi_{j}(s))|.
\end{align}
Let $\{a_{j}\}_{j=0}^{\infty}\subset [0,\infty)$ satisfying $A := \sum_{j=0}^{\infty}a_{j}<\infty$,
where $a_{j}$ will be chosen later.
 By (\ref{telescoping}), for any $u>0$,
\begin{align}\notag
\mathbb{P}\left(\underset{s> 0}{\sup}\ |W_{\Phi}(0,s)|\geq u\right)
=&\mathbb{P}\left(\underset{s> 0}{\sup}\ \overset{\infty}{\underset{j=0}{\sum}} |W_{\Phi}(0,\pi_{j+1}(s))|-|W_{\Phi}(0,\pi_{j}(s))|\geq uA^{-1}\overset{\infty}{\underset{j=0}{\sum}} a_{j}\right)
\\\label{ineq:supX}\leq& \overset{\infty}{\underset{j=0}{\sum}}\mathbb{P}\left(\underset{s>0}{\sup}\  |W_{\Phi}(0,\pi_{j+1}(s))|-|W_{\Phi}(0,\pi_{j}(s))|\geq u A^{-1}a_{j}\right).
\end{align}
For any $s>0$ and $j\in \mathbb{N}\cup\{0\}$,
\begin{align}\notag
&\mathbb{P}\left(|W_{\Phi}(0,\pi_{j+1}(s))|-|W_{\Phi}(0,\pi_{j}(s))|\geq u A^{-1} a_{j}\right)
\\\notag\leq&\mathbb{P}\left(|W_{\Phi}(0,\pi_{j+1}(s))-W_{\Phi}(0,\pi_{j}(s))|\geq u A^{-1} a_{j}\right)
\\\label{tail_prob_Xs_v1}=&\mathbb{P}\left(\frac{|W_{\Phi}(0,\pi_{j+1}(s))-W_{\Phi}(0,\pi_{j}(s))|^{2}}{d^{2}_{W_{\Phi}}(\pi_{j+1}(s),\pi_{j}(s))}
\geq \frac{u^{2} A^{-2}a_{j}^{2}}{d^{2}_{W_{\Phi}}(\pi_{j+1}(s),\pi_{j}(s))}\right),
\end{align}
where
\begin{align*}
d^{2}_{W_{\Phi}}(\pi_{j+1}(s),\pi_{j}(s))= \mathbb{E}\left[\left|W_{\Phi}(0,\pi_{j+1}(s))-W_{\Phi}(0,\pi_{j}(s))\right|^{2}\right].
\end{align*}
By Proposition \ref{lemma:exp}(b), the random variable
$$\frac{|W_{\Phi}(0,\pi_{j+1}(s))-W_{\Phi}(0,\pi_{j}(s))|^{2}}{d^{2}_{W_{\Phi}}(\pi_{j+1}(s),\pi_{j}(s))}$$
follows an exponential distribution with a mean of one. Continuing from (\ref{tail_prob_Xs_v1}),
for any $s>0$ and $j\in \mathbb{N}\cup\{0\}$,
\begin{align}\notag
&\mathbb{P}\left(|W_{\Phi}(0,\pi_{j+1}(s))|-|W_{\Phi}(0,\pi_{j}(s))|\geq u A^{-1}a_{j}\right)
\\\label{tail_prob_Xs_v2}\leq& \textup{exp}\left(-\frac{u^{2} A^{-2} a_{j}^{2}}{d^{2}_{W_{\Phi}}(\pi_{j+1}(s),\pi_{j}(s))}\right)
\leq \textup{exp}\left(-u^{2} A^{-2} a_{j}^{2}\mathcal{D}^{2(j-1)}\right),
\end{align}
where the last inequality follows from that $d_{W_{\Phi}}\left(\pi_{j+1}(s),\pi_{j}(s)\right)\leq \mathcal{D}^{1-j}$.
Because the number of pairs $\{(\pi_{j}(s),\pi_{j+1}(s))\mid s\geq 0\}$ is less than $|\pi_{j}| |\pi_{j+1}|$,
(\ref{ineq:supX}) and (\ref{tail_prob_Xs_v2}) imply that
\begin{align}\notag
\mathbb{P}\left(\underset{s>0}{\sup}\ |W_{\Phi}(0,s)|\geq u\right)
\leq \overset{\infty}{\underset{j=0}{\sum}}|\pi_{j}| |\pi_{j+1}|
\textup{exp}\left(-u^{2} A^{-2} a_{j}^{2}\mathcal{D}^{2(j-1)}\right)
\end{align}
for any $u\geq0$, and
\begin{align}\notag
\mathbb{E}\left[\underset{s>0}{\sup}\ |W_{\Phi}(0,s)|\right] =&
\int_{0}^{\infty}\mathbb{P}\left(\underset{s>0}{\sup}\ |W_{\Phi}(0,s)|\geq u\right)du
\\\label{ineq:supX:v2}\leq& A+\int_{A}^{\infty}\overset{\infty}{\underset{j=0}{\sum}}|\pi_{j}| |\pi_{j+1}|
\textup{exp}\left(-u^{2} A^{-2} a_{j}^{2}\mathcal{D}^{2(j-1)}\right)du.
\end{align}
By choosing
\begin{align*}
a_{j} = \mathcal{D}^{1-j}\sqrt{\ln\left(2^{j+1}|\pi_{j}| |\pi_{j+1}|\right)},
\end{align*}
the integrand function in (\ref{ineq:supX:v2}) can be estimated as follows
\begin{align}\label{estimate:integrand}
\overset{\infty}{\underset{j=0}{\sum}}|\pi_{j}| |\pi_{j+1}|
\textup{exp}\left(-u^{2} A^{-2} a_{j}^{2}\mathcal{D}^{2(j-1)}\right)
= \overset{\infty}{\underset{j=0}{\sum}}|\pi_{j}| |\pi_{j+1}|
\left(2^{j+1}|\pi_{j}| |\pi_{j+1}|\right)^{-u^{2}A^{-2}}\leq 2^{-u^{2}A^{-2}}
\end{align}
for any $u\geq A$. By substituting (\ref{estimate:integrand}) into (\ref{ineq:supX:v2}),
\begin{align*}
\mathbb{E}\left[\underset{s> 0}{\sup}\ |W_{\Phi}(0,s)|\right]
\leq A\left(1+\int_{1}^{\infty}2^{-v^{2}}dv\right)\leq\frac{4}{3}A.
\end{align*}
Because
\begin{align}\label{def:normalizer}
A = \overset{\infty}{\underset{j=0}{\sum}}\mathcal{D}^{1-j}\sqrt{\ln\left(2^{j+1}|\pi_{j}| |\pi_{j+1}|\right)},
\end{align}
by substituting (\ref{def:cardinality}) into (\ref{def:normalizer}) and performing elementary calculation,
\begin{align}\notag
A \leq \sqrt{3\ln\mathcal{D}}\frac{\mathcal{D}^{3}}{(\mathcal{D}-1)^{2}}\sqrt{\frac{1}{\gamma}+\frac{1}{H^{-}}}
+\left(\sqrt{\frac{4}{\gamma}\ln C_{2}}+\sqrt{2\ln T}\right)\frac{\mathcal{D}^{2}}{\mathcal{D}-1}.
\end{align}
\qed

\subsection{Proof of Lemma \ref{lemma:complex-Borell}}\label{sec:proof:lemma:complex-Borell}
Since the proofs of the first and second parts of Lemma \ref{lemma:complex-Borell} are similar, with the first being simpler, we provide only the proof of the second part below.
Consider a complex-valued Gaussian process as follows
\begin{align*}
V(s) = w(s)\frac{\partial W_{\Phi}}{\partial s}(0,s),\ s\in [a,b].
\end{align*}
As shown in Proposition \ref{lemma:twice_differentiability}, $W_{\Phi}$ is two times continuously differentiable,
which implies that $V$ is continuous on $[a,b]$ almost surely, and
\begin{align}\label{finite:sup:V:as}
\mathbb{P}\left(\|V\|_{[a,b]}<\infty\right) = 1,\ \textup{where}\ \|V\|_{[a,b]}:= \underset{s\in [a,b]}{\max} |V(s)|.
\end{align}
Using (\ref{spect:WPhi}) and denoting the derivative of $\widehat{\psi}$ by $D\widehat{\psi}$, we have
\begin{align}\notag
V(s) =& w(s)\int_{0}^{\infty} \overline{[D\widehat{\psi}](s\lambda)}\lambda\sqrt{p(\lambda)}Z(d\lambda)
\\\label{spect:gWPhi}=& w(s)\int^{0}_{-\infty} \overline{[D\widehat{\psi}](-s\lambda)}\lambda\sqrt{p(\lambda)}\ \overline{Z(d\lambda)},
\end{align}
which implies that for any $s_{1},s_{2}\in [a,b]$,
\begin{align}\label{check:circular}
\mathbb{E}[V(s_{1})V(s_{2})] = 0.
\end{align}
Let $\{T_{d}\}_{d\in \mathbb{N}}$ be a sequence of subsets of $[a,b]$ such that $T_{d}$ has the cardinality $d$, $T_{d}\subset T_{d+1}$,
and $T_{d}$ increases to a dense subset of $[a,b]$.
By (\ref{check:circular}), $\mathbf{V}_{d}:=\{V(s)\mid s\in T_{d}\}$ is a circularly symmetric complex Gaussian random vector.
By Lemma \ref{lemma:Lip}, for any $u>0$,
\begin{align}\label{Lip_ineq:Va}
\mathbb{P}\left(\pm\left(\|\mathbf{V}_{d}\|_{\ell^{\infty}}
-\mathbb{E}\left[\|\mathbf{V}_{d}\|_{\ell^{\infty}}\right]\right)> u\right) \leq
\textup{exp}\left(-\frac{u^{2}}{2\|\|\mathbf{A}_{d}\bullet\|_{\ell^{\infty}}\|_{\textup{Lip}}^{2}}\right),
\end{align}
where
$\mathbf{A}_{d}$ is a matrix in $\mathbb{C}^{d\times d}$
satisfying
\begin{align}\label{relation:AAVV}
\mathbf{A}_{d}\mathbf{A}_{d}^{*}=2^{-1}\mathbb{E}\left[\mathbf{V}_{d}\mathbf{V}_{d}^{*}\right].
\end{align}
For any $\mathbf{x},\mathbf{y}\in \mathbb{C}^{d\times 1}$ with $\mathbf{x}\neq\mathbf{y}$,
\begin{align}\notag
|\|\mathbf{A}_{d}\mathbf{x}\|_{\ell^{\infty}}-\|\mathbf{A}_{d}\mathbf{y}\|_{\ell^{\infty}}|
\leq& \underset{1\leq k\leq d}{\max} \left|e_{k}^{\top}\mathbf{A}_{d}\left(\mathbf{x}-\mathbf{y}\right)\right|
\\\notag\leq& \underset{1\leq k\leq d}{\max} \left|e_{k}^{\top}\mathbf{A}_{d}\right| \left|\mathbf{x}-\mathbf{y}\right|
\\\label{Lip_ineq:Vb}=&\left(\underset{1\leq k\leq d}{\max} e_{k}^{\top}\mathbf{A}_{d}\mathbf{A}_{d}^{*}e_{k}\right) \left|\mathbf{x}-\mathbf{y}\right|,
\end{align}
where $e_{k}$ is the column vector in $\mathbb{R}^{d}$ with one in position $k$ and zeros elsewhere.
If we denote the elements of $T_{d}$ by $s_{d,1},s_{d,2},\ldots,s_{d,d}$, then (\ref{spect:gWPhi}) and (\ref{relation:AAVV}) imply that
\begin{align}\label{comput_Lip_eAde}
e_{k}^{\top}\mathbf{A}_{d}\mathbf{A}_{d}^{*}e_{k}
=2^{-1}w^{2}(s_{d,k})\int_{0}^{\infty} |[D\widehat{\psi}](s_{d,k}\lambda)|^{2}\lambda^{2}p(\lambda)d\lambda.
\end{align}
By combining (\ref{Lip_ineq:Va}), (\ref{Lip_ineq:Vb}), and (\ref{comput_Lip_eAde}),
we obtain that for any $u>0$,
\begin{align}\label{Lip_ineq:summary1}
\mathbb{P}\left(\pm \left(\|\mathbf{V}_{d}\|_{\ell^{\infty}}-\mathbb{E}\left[\|\mathbf{V}_{d}\|_{\ell^{\infty}}\right]\right)> u\right) \leq
\textup{exp}\left(-\frac{u^{2}}{\sigma_{T_{d}}^{2}}\right),
\end{align}
where
\begin{align*}
\sigma_{T_{d}}^{2}=\underset{1\leq \ell\leq d}{\max}\ w^{2}(s_{d,\ell})\int_{0}^{\infty} |[D\widehat{\psi}](s_{d,\ell}\lambda)|^{2}\lambda^{2}p(\lambda)d\lambda.
\end{align*}

Because $T_{d}$ increases to a dense subset of $[a,b]$ as $d\rightarrow\infty$,
the continuity of $V$ on $[a,b]$ implies that
\begin{align}\label{gdVd:;conv}
\|\mathbf{V}_{d}\|_{\ell^{\infty}}\uparrow\|V\|_{[a,b]}<\infty
\end{align}
almost surely, and under Assumption \ref{assumption:boundedness:psi},
\begin{align}\label{sigmaTd_conv}
\sigma_{T_{d}}^{2}\uparrow \sigma_{[a,b]}^{2}=\underset{s\in [a,b]}{\max}\ w^{2}(s)\int_{0}^{\infty} |[D\widehat{\psi}](s\lambda)|^{2}\lambda^{2}p(\lambda)d\lambda<\infty.
\end{align}
On the other hand, when $d$ increases,
\begin{align}\label{EVTd:increasing}
\mathbb{E}\left[\|\mathbf{V}_{d}\|_{\ell^{\infty}}\right]\uparrow \mathbb{E}\left[\|V\|_{[a,b]}\right].
\end{align}
If the limit $\mathbb{E}[\|V\|_{[a,b]}]$ is finite,
then for any $u>0$, (\ref{gdVd:;conv}) and (\ref{EVTd:increasing}) imply that
\begin{align}\notag
\mathbb{P}\left(\pm\left(\|V\|_{[a,b]}-\mathbb{E}\left[\|V\|_{[a,b]}\right]\right)> u\right)
=&\underset{d\rightarrow\infty}{\lim} \mathbb{P}\left(\pm\left(\|\mathbf{V}_{d}\|_{\ell^{\infty}}-\mathbb{E}\left[\|\mathbf{V}_{d}\|_{\ell^{\infty}}\right]\right)> u\right)
\\\notag\overset{(\ref{Lip_ineq:summary1})}{\leq}&
\underset{d\rightarrow\infty}{\lim}\textup{exp}\left(-\frac{u^{2}}{\sigma_{T_{d}}^{2}}\right)
\overset{(\ref{sigmaTd_conv})}{=}
\textup{exp}\left(-\frac{u^{2}}{\sigma_{[a,b]}^{2}}\right).
\end{align}
Therefore, to complete the proof, it suffices to show that $\mathbb{E}[\|V\|_{[a,b]}]<\infty$.

We proceed by contradiction. Thus, assume that $\mathbb{E}[\|V\|_{[a,b]}]=\infty$, and choose $u_{0}>0$ such that
\begin{align*}
\exp\left(-\frac{u_{0}^{2}}{\sigma_{[a,b]}^{2}}\right)< \frac{1}{3}\ \textup{and}\ \mathbb{P}\left(\|V\|_{[a,b]}<u_{0}\right)\geq \frac{1}{2},
\end{align*}
which is attainable due to (\ref{finite:sup:V:as}).
Now choose $d\geq1$ such that $\mathbb{E}[\|\mathbf{V}_{d}\|_{\ell^{\infty}}]>2u_{0}$, which is attainable due to (\ref{EVTd:increasing}) and $\mathbb{E}[\|V\|_{[a,b]}]=\infty$. By (\ref{Lip_ineq:summary1}),
\begin{align}\notag
\frac{1}{3}>& \exp\left(-\frac{u_{0}^{2}}{\sigma_{[a,b]}^{2}}\right)\geq \exp\left(-\frac{u_{0}^{2}}{\sigma_{T_{d}}^{2}}\right)
\\\notag\geq& \mathbb{P}\left(-\left(\|\mathbf{V}_{d}\|_{\ell^{\infty}}-\mathbb{E}[\|\mathbf{V}_{d}\|_{\ell^{\infty}}]\right)> u_{0}\right)
\\\notag=&
 \mathbb{P}\left(\|\mathbf{V}_{d}\|_{\ell^{\infty}}<\mathbb{E}[\|\mathbf{V}_{d}\|_{\ell^{\infty}}]- u_{0}\right)
\\\notag\geq&
\mathbb{P}\left(\|\mathbf{V}_{d}\|_{\ell^{\infty}}< u_{0}\right)
\\\notag\geq& \mathbb{P}\left(\|V\|_{[a,b]}<u_{0}\right)
\geq \frac{1}{2},
\end{align}
which is contradictory.
\qed
%
%
%
%
%
%
%
%
%
%
%
%
%
\bibliographystyle{abbrvnat}
\bibliography{ref_nonnull202412_journal_abbrev}

\vskip 20 pt
\newpage
{\center{\Large  \textbf{Appendix\\}}}

\hspace{-6cm}
\numberwithin{equation}{section}

\begin{table}[h]
\begin{center}\begin{tabular}{p{2.5cm}p{13cm}}
\toprule
{\bf Symbol} & {\bf Description}
\\
\midrule
$f$  & Adapted harmonic signal with $M$ components (see (\ref{AMFM_signal})).\\
$\phi'_{1},\ldots,\phi'_{M}$ & IFs of $f$ with $\phi'_{1}(t)<\cdots<\phi'_{M}(t)$ for any $t\in \mathbb{R}$.\\
$A_{1},\ldots,A_{M}$ & Amplitude functions of harmonic components of $f$.\\
$\Phi$      & Stationary Gaussian process (see Assumption \ref{assumption:Gaussian}).\\
$p$    & Spectral density function of $\Phi$.\\
$Y$      & Noisy signal modeled as $Y=f+\Phi$.\\
$Z$      & Orthogonally scattered Gaussian random measure on $\mathbb{R}$.\\
$\psi$, $\widehat{\psi}$ & Analytic mother wavelet and its Fourier transform.\\
$\psi_{R}$, $\psi_{I}$ & Real and imaginary parts of the wavelet $\psi$.\\
$D\widehat{\psi}, D^{2}\widehat{\psi}$ & First and second derivatives of $\widehat{\psi}$.\\
$W_{f},W_{\Phi},W_{Y}$       & Continuous wavelet transform of $f$, $\Phi$, and $Y$.\\
$W^{R}_{\Phi},W^{I}_{\Phi}$       & Real and imaginary parts of $W_{\Phi}$.\\
$S_{f}, S_{\Phi},S_{Y}$   & Scalograms of $f$, $\Phi$, and $Y$, i.e., the squared modulus of $W_{f}$, $W_{\Phi}$, and $W_{Y}$.\\
$s_{f}(t)$ & Argmax of the scalogram of $f$ at time $t$, i.e., $\arg \underset{s>0}{\max}\ S_{f}(t,s)$. \\
$s_{f,m}(t)$ &  Local maxima of $f$'s scalogram at time $t$, ordered as  $s_{f,1}>\cdots>s_{f,M}$.\\
$s_{Y}(t)$ & Set of scales maximizing $S_{Y}(t,\cdot)$, i.e., $s_{Y}(t) = \arg \underset{s>0}{\max}\ S_{Y}(t,s)$.\\
$\omega_{\psi}$ & Argmax of  $|\widehat{\psi}|$.\\
$B_{m}(t)$  & IF information of $f$: $\omega_{\psi}(\phi'_{m}(t))^{-1}\in B_{m}(t)$ (see (\ref{prior_information})).\\
$s_{Y,m}(t)$   & Scales in $B_{m}(t)$ maximizing $S_{Y}(t,\cdot)$, i.e., $s_{Y,m}(t) = \arg \underset{s\in B_{m}(t)}{\max}\ S_{Y}(t,s)$.\\
\bottomrule
\end{tabular}
\end{center}
\caption{List of frequently used symbols.}
\label{List_symbols}
\end{table}

\begin{Lemma}[Berge's maximum theorem]\label{lemma:Berge}
Let $\mathbb{X}$ and $\mathbb{Y}$ be Hausdorff topological spaces. If $u: \mathbb{X}\times \mathbb{Y}\rightarrow \mathbb{R}$ is a continuous function
and $F: \mathbb{X}\rightarrow  \{K\subseteq \mathbb{Y}:\ K\ \textup{is a nonempty compact set}\}$ is a continuous multifunction,
then
\begin{align*}
\arg \underset{y\in F(x)}{\sup} u(x,y) : = \left\{z\in F(x):\ u(x,z) = \underset{y\in F(x)}{\sup} u(x,y) \right\}
\end{align*}
is an upper hemicontinuous multifunction.
\end{Lemma}

\begin{Lemma}\cite[Lemma (4.3.16)]{MR1039321}\label{lemma:continuous_dense}
If $\mathbb{X}$ is a topological space, $\mathbb{Y}$ a complete metric space and
$r : \mathbb{D} \mapsto \mathbb{Y}$  a continuous mapping defined on a dense subset $\mathbb{D}$ of the space $\mathbb{X}$,
then the mapping $r$ is extendable to a continuous mapping  defined
on the set consisting of all points of $\mathbb{X}$ at which the oscillation of $r$ is equal
to zero. Here, we say that the
oscillation of the mapping $r$ at a point $x \in \mathbb{X}$ is equal to zero if for every
$\varepsilon > 0$
there exists a neighborhood $U$ of the point $x$ such that
\begin{align*}
\underset{x',x''\in \mathbb{D}\cap U}{\sup} r(x')-r(x'')<\varepsilon.
\end{align*}
The set of all points at which the oscillation of $r$ is equal to zero is a $G_\delta$-set containing
the dense set $\mathbb{D}$.
\end{Lemma}

\begin{figure}[hbt!]
\centering
\includegraphics[width=0.99\textwidth]{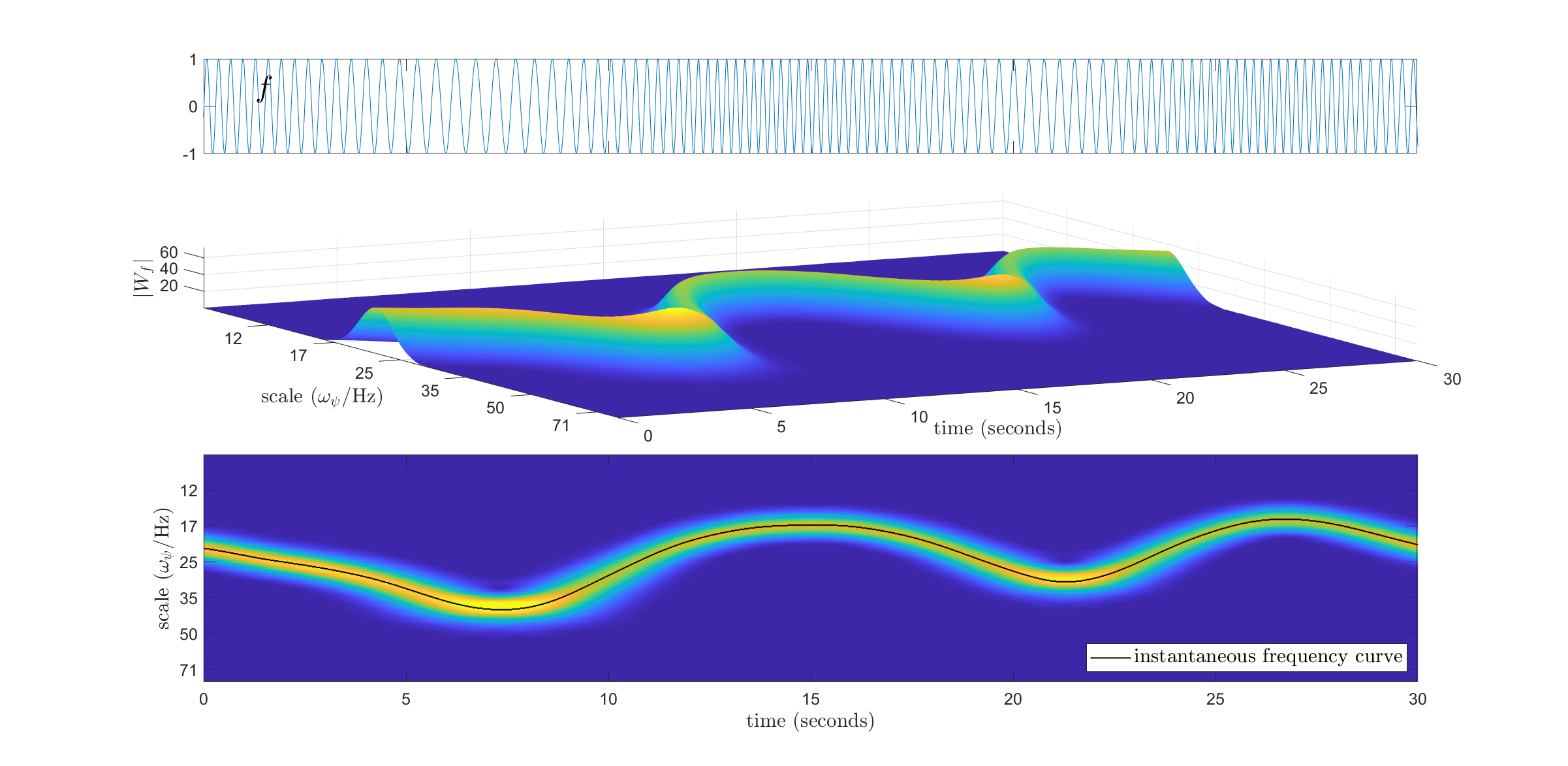}
\hspace{-1.5cm}
\includegraphics[width=0.99\textwidth]{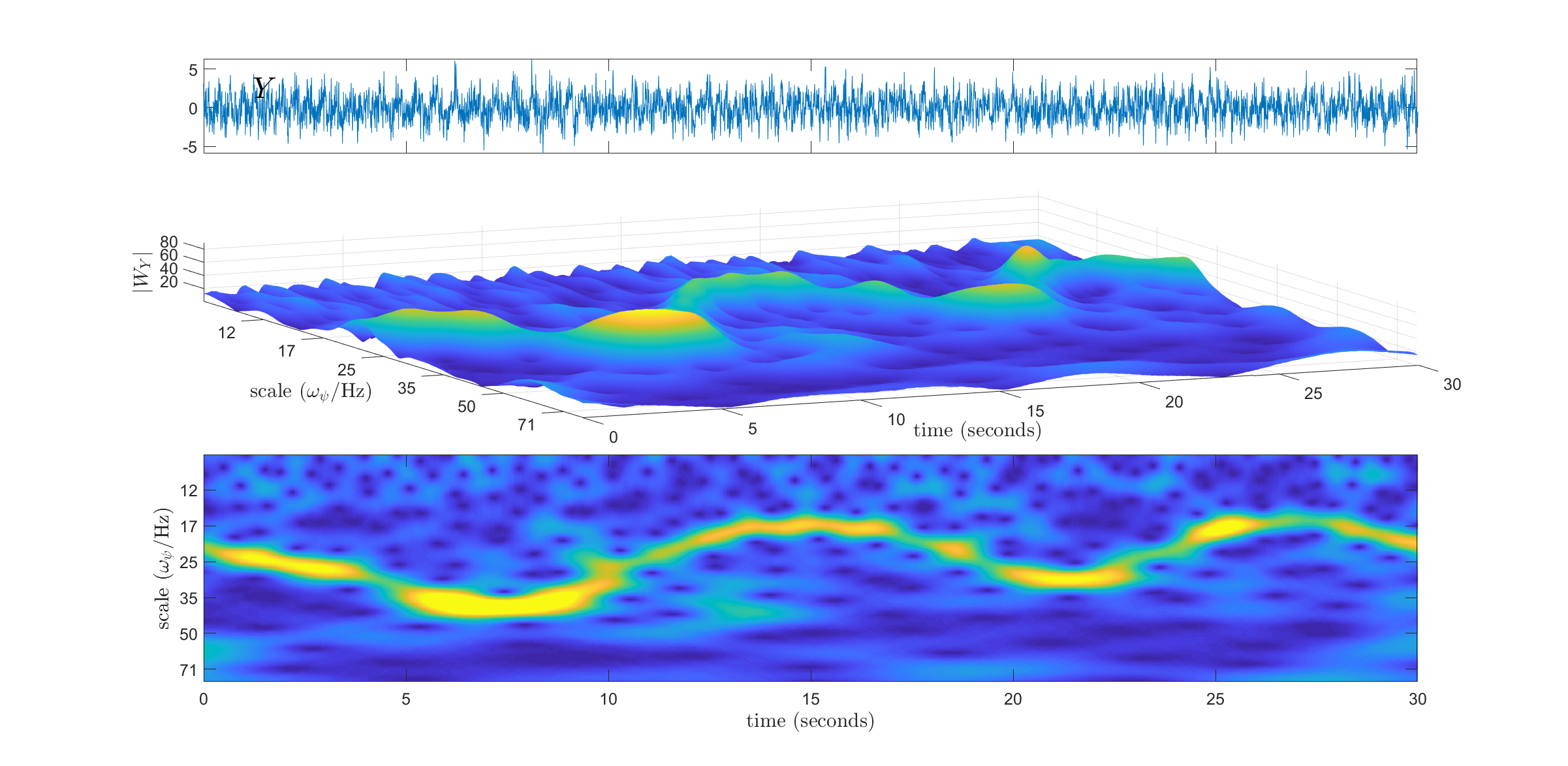}
 \caption{The top section illustrates a frequency-modulated signal $f$ along with its 3-dimensional and image representations in the time-scale domain.
The bottom section presents the same for its noise-affected counterpart $Y$, with a signal-to-noise ratio of -6.72 dB. Here, $Y$ is the sum of $f$ and a sample path of the Gaussian process $\Phi$.
  \label{fig:wavelet_potential_ridge_Surface plot:one}}
\end{figure}

\begin{figure}[hbt!]
  \centering
  \includegraphics[width=0.99\textwidth]{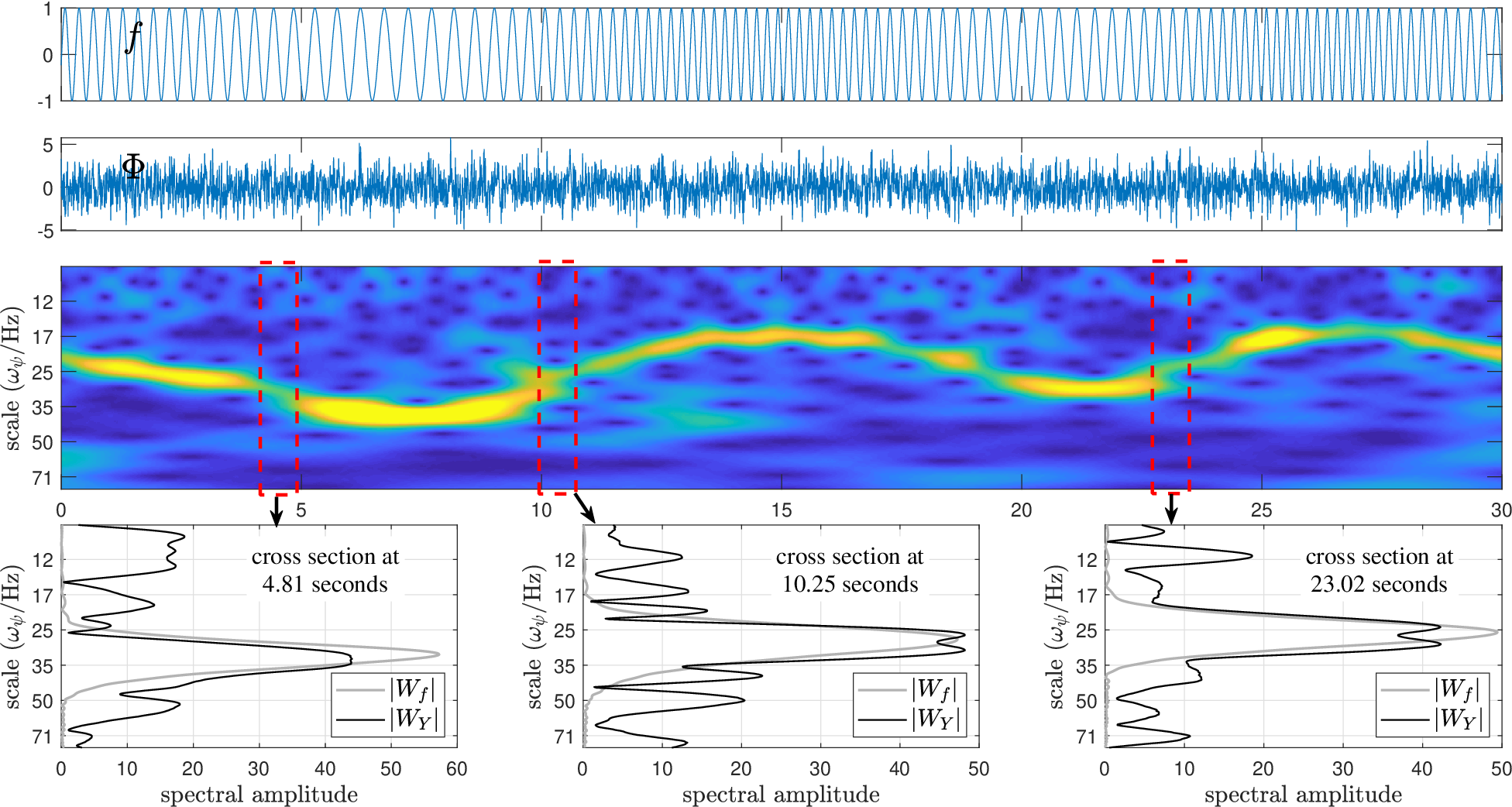}
 \caption{Example of the complex modulus of the AWT for a 30-second noise-contaminated and frequency-modulated signal.
The first row displays the frequency-modulated signal $f$, while the second row shows the noise, a sample from the stationary Gaussian process $\Phi$, identical to that in Figure \ref{fig:wavelet_potential_ridge_Surface plot}. The third row illustrates the complex modulus of the AWT, i.e., the spectral amplitude, of the noise-contaminated signal $Y$, where $Y=f+\Phi$. In the last three columns, the black curves represent the spectral amplitude of $Y$ at specific times, corresponding to cross-sections of the image in the third row. For comparison, the spectral amplitudes of $f$ are shown in gray.}\label{fig:wavelet_ridge:one}
\end{figure}

\begin{figure}[hbt!]
\centering
\subfigure[][]{\includegraphics[width=0.99\textwidth]{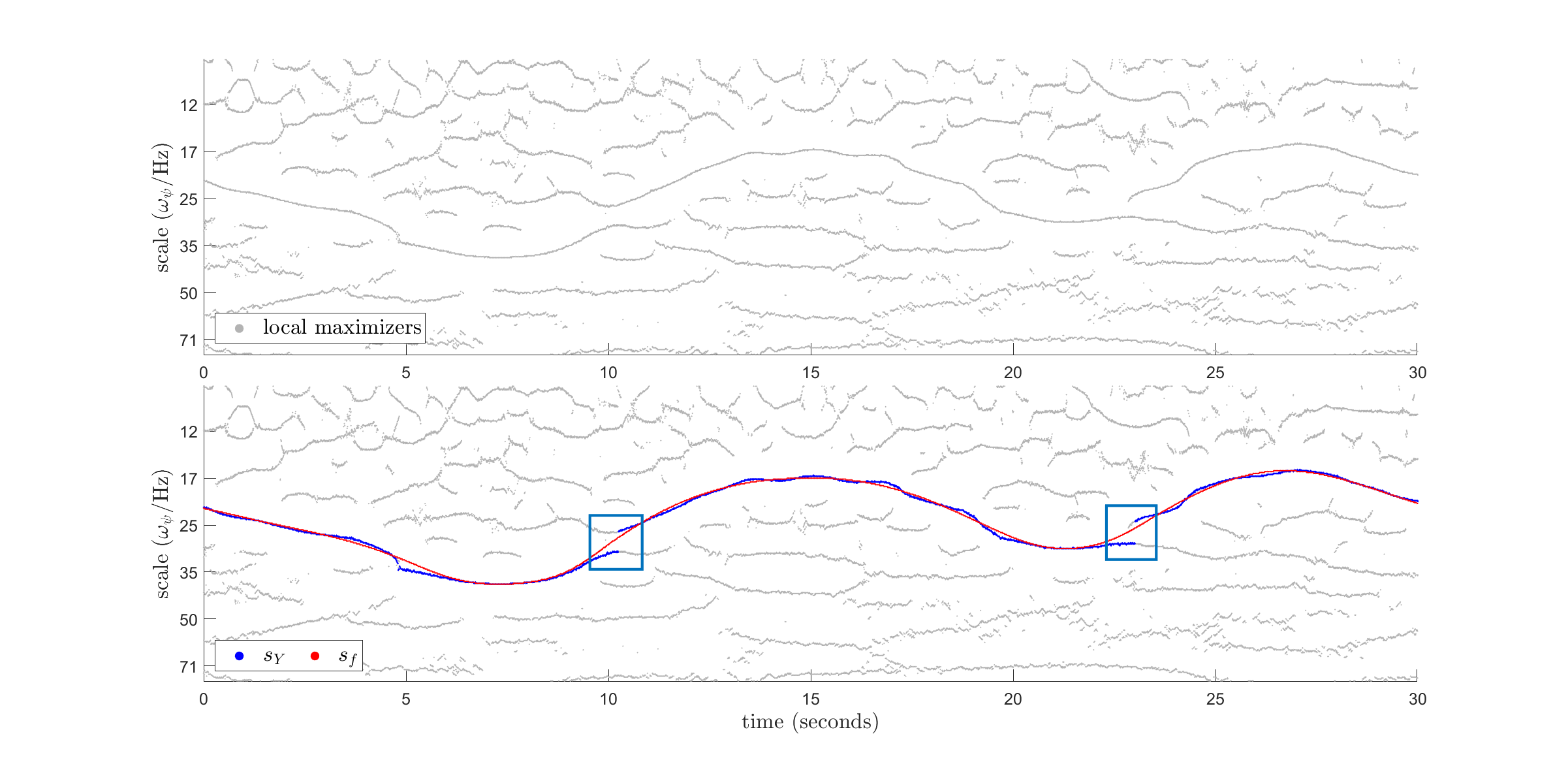}}
\subfigure[][]{\includegraphics[width=0.99\textwidth]{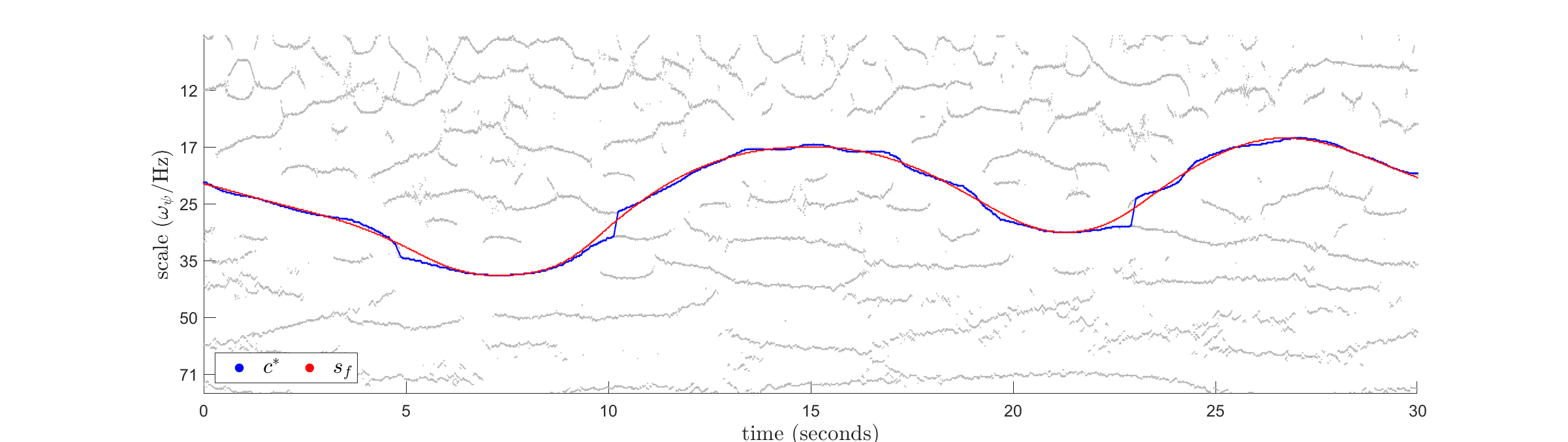}}
\caption{The top panel shows the ridge points corresponding to local maxima of the scalogram of $Y$ along the scale axis,
where $Y = f+\Phi$. The signal $f$ and the realization of $\Phi$ are identical to those in Figure \ref{fig:wavelet_potential_ridge_Surface plot:one}.
The middle panel displays the ridge curves $s_{Y}$ and $s_{f}$.
The blue boxes highlight the discontinuity of the ridge in the presence of noise.
The bottom panel shows the ridge curve extracted by the ridge extraction algorithm \eqref{eq ridge ex}.
  \label{fig:wavelet_potential_ridge:one}}
\end{figure}

\end{document}